\documentclass{article}

\usepackage{PRIMEarxiv}

\usepackage[utf8]{inputenc} 
\usepackage[T1]{fontenc}    
\usepackage{url}            
\usepackage{booktabs}       
\usepackage{amsfonts}       
\usepackage{nicefrac}       
\usepackage{microtype}      
\usepackage{fancyhdr}       
\usepackage{graphicx}       
\graphicspath{{media/}}     

\usepackage{amsmath,amssymb,amsthm}
\usepackage{tikz}
\usepackage{algorithm}
\usepackage{algpseudocode}
\usepackage{subcaption}
\usepackage{placeins}
\usepackage{hyperref}       

\newtheorem{theorem}{Theorem}[section]
\newtheorem{proposition}[theorem]{Proposition}
\newtheorem{corollary}[theorem]{Corollary}
\newtheorem{assumption}[theorem]{Assumption}
\newtheorem{lemma}[theorem]{Lemma}
\theoremstyle{definition}
\newtheorem{remark}[theorem]{Remark}

\title{A GPU Interior-Point Solver with Macro-Element Condensation for Equilibrium-Based Plasticity
}

\author{
	Youssef Essousy \\
	College of Computing \\
	Mohammed VI Polytechnic University (UM6P) \\
	Ben Guerir\\
	\texttt{youssef.essoussy@um6p.ma} \\
	\And
	Mohammed-Khalil Ferradi \\
	College of Computing \\
	Mohammed VI Polytechnic University (UM6P) \\
	Ben Guerir\\
	\texttt{khalil.ferradi@um6p.ma} \\
}

\begin{document}
	\maketitle

\begin{abstract}
Equilibrium-based conic formulations of elastoplasticity and limit analysis provide statically admissible stress
fields but generate large Newton systems whose sparse factorization limits three-dimensional computations.
This work develops and analyzes a macro-element condensation tailored to the Newton systems arising from primal--dual
interior-point methods. Exploiting four-tetrahedron macro-elements, the method eliminates stress directions and
macro-internal traction-continuity multipliers, yielding reduced systems that can be solved using positive-definite
factorizations, with estimated reductions of approximately \(4\times\) in dimension and structural nonzeros relative
to stress-only elimination. The validity of this elimination is established through global and macro-local
nonsingularity conditions. Local nonsingularity is independent of the yield criterion in elastoplasticity, whereas in limit
analysis it holds directly for Rankine and is ensured for von Mises through trace-bound augmentation. Numerical results
show that, relative to stress-only elimination, the proposed condensation achieves CPU speedups of up to approximately \(2\times\)
and reduces peak linear-solver memory by up to \(66.5\%\). We further develop a CUDA Fortran implementation that exploits
uniform macro-local blocks through batched routines and uses cuDSS for the reduced sparse solve. This implementation
achieves speedups of up to \(13.5\times\) over its CPU counterpart and up to \(87.5\times\) over MOSEK, while successfully
solving the largest tested 3D problems that exceed the device-memory capacity of the GPU solver CuClarabel. For elastoplasticity,
we also propose a preconditioned conjugate-gradient solver that reduces device-memory requirements relative to the direct
GPU solver, at the expense of increased solution time.
\end{abstract}



\section{Introduction}
\label{sec:introduction}

Assessing structural behavior beyond the elastic range and determining the ultimate load-bearing capacity are central
problems in computational solid mechanics. Elastoplastic analysis describes the redistribution of stresses as yielding
develops under prescribed loading~\cite{kanno2011nonsmooth}, whereas limit analysis determines the collapse load
directly without resolving the complete intermediate response~\cite{krabbenhoft2007formulation}. These two analyses
therefore provide complementary information for assessing and safely designing structures operating near their material
limits. Both admit kinematic, stress-based, mixed, and hybrid variational formulations~\cite{simo1989complementary,krabbenhoft2007interior}.
In limit analysis, the kinematic and static principles yield upper and lower bounds on the collapse load,
respectively~\cite{makrodimopoulos2006lower,bleyer2021automated}. In the equilibrium-based description, stresses are
treated as primary unknowns and constrained to satisfy equilibrium and yield admissibility directly.

Equilibrium finite elements provide a natural discretization of this stress-based description~\cite{ferradi2023elastoplastic}.
Many convex yield criteria admit second-order cone (SOC) or positive semidefinite (PSD) cone
representations~\cite{makrodimopoulos2006lower,bisbos2007second,krabbenhoft2012computational,krabbenhoft2008three}. With
these conic yield formulations and the linear equilibrium constraints, equilibrium-based elastoplasticity and limit analysis
can be expressed within a common conic optimization framework~\cite{alizadeh2003second,vandenberghe1996semidefinite}. Within
this equilibrium discretization, however, an individual tetrahedron may exhibit nonphysical zero-energy modes in addition to
the six rigid-body modes~\cite{maunder1996general}. The super-element construction of Kempeneers et al.~\cite{kempeneers2010pure}
eliminates these modes by subdividing a tetrahedron into four sub-elements connected at a barycentric node; the stability of the
resulting four-tetrahedron patch is also confirmed by the star-patch analysis of Maunder et al.~\cite{maunder2016stability}.
However, the equilibrium description still retains both the stress unknowns and the traction-continuity constraints, so the
resulting algebraic systems are typically larger than those arising from continuous kinematic formulations, particularly in
3D~\cite{makrodimopoulos2007upper,el2020elastoplastic}.

Large conic programs are commonly solved using first-order methods~\cite{tu2014practical} or primal--dual
interior-point methods (IPMs)~\cite{mehrotra1992implementation,nesterov1997self}. First-order methods offer inexpensive
iterations and low memory requirements, but attaining tight feasibility and optimality tolerances can require many
iterations, particularly for poorly conditioned large-scale problems. By contrast, IPMs typically achieve high accuracy
within a few tens of iterations, although each iteration requires forming and solving a large sparse Newton system
obtained by linearizing the Karush--Kuhn--Tucker conditions~\cite{nesterov1997self,nesterov1998primal}. Their robustness
and accuracy make IPMs particularly suitable for computational mechanics, where accurate stress fields and load multipliers
are often required, and underpin established general-purpose conic solvers such as MOSEK~\cite{mosek112},
Clarabel~\cite{goulart2026clarabel}, and CVXOPT~\cite{vandenberghe2010cvxopt}. Because the linear solve dominates the
cost of an IPM, the dimension, sparsity, and numerical structure of the Newton system largely determine its computational
and memory requirements, independently of the computing architecture.

Graphics processing units (GPUs) are increasingly used in scientific computing because their high memory bandwidth and
parallelism are well suited to data-parallel numerical operations. These capabilities make GPUs attractive for the
solution of large conic models. Existing GPU implementations have therefore focused largely on first-order and matrix-free
algorithms, whose dominant sparse matrix--vector products and cone projections map naturally to GPU
hardware~\cite{lin2025pdcs,han2024accelerating,lu2025overview}. Although conic IPMs also benefit from parallel cone
operations, their efficient GPU implementation remains more challenging because the dominant computational cost lies in
solving these sparse Newton systems. CPU implementations commonly rely on mature sparse direct solvers and can draw on
substantially larger host-memory capacities. These solvers provide robust Newton directions, but fill-in during
factorization can still require substantial storage~\cite{duff1983multifrontal,benzi2005numerical}. On GPUs, the smaller
available device memory makes this strategy particularly restrictive for large 3D problems. To reduce this storage
requirement, GPU-accelerated matrix-free and preconditioned iterative strategies have been investigated within
IPMs~\cite{smith2012gpu,kardani2015application}. However, as an IPM approaches convergence, the Newton systems become
increasingly ill-conditioned, and robust iterative convergence requires strong, often problem-dependent preconditioners.
Sparse direct solvers therefore remain the standard choice when sufficient factorization memory is available, whereas
iterative methods provide a memory-oriented alternative for problems beyond that limit~\cite{gondzio2012interior}.

Recent general-purpose GPU conic IPMs offering sparse-direct KKT solution paths include
QOCO-GPU~\cite{chari2026qocogpu}, CuClarabel~\cite{chen2026cuclarabel}, and Moreau~\cite{moreau2026}.
QOCO-GPU supports conic blocks only from the nonnegative orthant and SOCs. By contrast, CuClarabel, the GPU
implementation of Clarabel~\cite{goulart2026clarabel}, and Moreau, a GPU-native differentiable solver whose interior-point
algorithm derives from Clarabel and CuClarabel, support broader cone families that include PSD cones.
All three solvers provide general sparse KKT solution paths based on the NVIDIA cuDSS sparse direct
solver~\cite{cudss080}. When such a path is used, sparse-factor storage can limit the problem size. Beyond conic optimization,
recent condensed-space IPM strategies for nonlinear programming exploit positive-definite reformulations to enable GPU-based
Cholesky factorization without numerical pivoting~\cite{pacaud2026condensed}.

In parallel with these developments in conic optimization, GPU acceleration has been widely investigated in finite element
analysis, including structural elasticity~\cite{cecka2011assembly,georgescu2013gpu,dick2011real}, topology
optimization~\cite{wu2015system,traff2023simple,qi2025efficient}, and elastoplasticity~\cite{kiran2023gpu,kiran2024development}.
These implementations accelerate element calculations and sparse assembly or employ matrix-free iterative solvers to
reduce storage. For conic finite-element limit analysis, Kardani et al.~\cite{kardani2015application} used a GPU-accelerated
hybrid preconditioned conjugate-gradient (PCG) method to solve the sparse linear systems arising at each IPM iteration.
These studies demonstrate the potential of GPUs for computational mechanics, including memory-oriented iterative strategies for
large-scale problems.

Algebraic condensation is classical in equilibrium finite-element formulations. For linear elasticity, Fraeijs de Veubeke
and Pian derived element stiffness relations by eliminating local assumed-stress parameters from complementary-energy
formulations \cite{fraeijs1965displacement,pian1964derivation}. Fraeijs de Veubeke additionally considered
assemblies of primitive triangular equilibrium elements and eliminated unloaded internal generalized displacement variables,
anticipating the topology-aware equilibrium macro-element methodology later systematized by Maunder and co-authors
\cite{maunder1996general,pereira1999adaptive,de2016equilibrium}. In stress-based plasticity, nonlinear-programming and interior-point
algorithms have similarly exploited the block-diagonal local stress curvature to eliminate stress directions and obtain a
global system in the equilibrium multipliers, both in lower-bound limit analysis \cite{lyamin2002lower,krabbenhoft2003general}
and in incremental elastoplasticity \cite{krabbenhoft2007interior}.
This stress-only reduction requires the local curvature blocks to be invertible; in limit analysis with a pressure-insensitive
criterion such as von Mises, the unrestricted hydrostatic stress mode makes these blocks singular. Existing approaches have
addressed this difficulty through small diagonal regularizations \cite{lyamin2002lower,krabbenhoft2003general} or an exact change
of stress coordinates that separates the bounded deviatoric variables from an explicitly retained unrestricted stress
component \cite{christiansen1999computation}.

Building on the equilibrium macro-element framework \cite{maunder1996general} and the four-tetrahedron construction of Kempeneers
et al.~\cite{kempeneers2010pure}, we propose a condensation of the Newton systems associated with equilibrium-based conic
elastoplasticity and limit analysis. The method jointly eliminates stress directions and macro-internal traction-continuity
multipliers, leaving a smaller system in the variables that remain globally coupled. This reduction requires the coupled local
blocks to remain invertible at every IPM iteration, a property that global Newton nonsingularity alone does not guarantee.
We therefore establish conditions for both global and macro-local nonsingularity, with elastic compliance providing local
invertibility in elastoplasticity irrespective of the yield criterion. For limit analysis, the corresponding analysis verifies
local nonsingularity for Rankine and introduces a bound on the stress trace to remove the hydrostatic null mode associated with
von Mises. Under these conditions, the reduced elastoplastic system admits a positive-definite factorization, while the
load-factor border in limit analysis can be handled through a scalar Schur complement of a positive-definite multiplier block.
For this macro-element partition, the structural reductions are estimated at approximately \(8.67\times\) in dimension and
\(4.85\times\) in structural nonzeros relative to the original Newton system, and approximately \(4\times\) in both quantities
relative to stress-only condensation. CPU comparisons with stress-only elimination then quantify the computational benefit of
the additional condensation: speedups reach \(1.97\times\), and the peak memory required by the linear solver for the condensed
Newton system decreases by up to \(66.5\%\).

The independent local eliminations expose parallel work, while the smaller condensed Newton system reduces the storage
requirements of the global sparse solve. Our CUDA Fortran implementation of the full primal--dual IPM exploits both properties,
with batched CUDA routines processing the uniform macro-local blocks and cuDSS solving the condensed sparse system. Numerical
assessment on a single NVIDIA A100 GPU with \(80\,\mathrm{GB}\) of device memory covers three 3D benchmarks: elastoplastic
Cook's membrane, a steel beam under torsion, and a reinforced-concrete beam under three-point
bending~\cite{ferradi2023elastoplastic,ferradi2024steel}. End-to-end comparisons show that the GPU solver is up to \(13.5\times\)
faster than its macro-condensed CPU counterpart (executed on \(24\) threads), while achieving speedups of up to \(87.5\times\)
over MOSEK and approximately \(31\times\) over CuClarabel. It completes all three benchmarks at their largest tested mesh sizes,
for which CuClarabel exhausts the available device memory. The storage demands of sparse direct factorization also motivate an
additional PCG solution path for elastoplasticity, equipped with a preconditioner based on the macro-element structure. This
iterative alternative lowers device memory usage, although the direct cuDSS path remains faster in the reported tests.

The nonsingularity analysis provides a basis for extending the proposed condensation beyond the formulations considered here.
Potential extensions include other conic yield criteria, such as Drucker--Prager~\cite{makrodimopoulos2006lower} and
frictional or truncated Mohr--Coulomb models~\cite{krabbenhoft2008three}, as well as equilibrium-based shakedown
analysis~\cite{zouain2002algorithm,bisbos2005second}, yield design~\cite{krabbenhoft2007formulation,bleyer2021automated},
and topology optimization~\cite{mourad2021topology,ferradi2024steel}. Each extension requires a formulation-specific
verification of the local and global nonsingularity conditions.

The remainder of the paper is organized as follows. The equilibrium finite-element discretization, conic representation,
and associated interior-point Newton system are introduced in Section~\ref{sec:equilibrium-conic-formulation}.
Section~\ref{sec:structure-exploiting-condensation} establishes global regularity conditions, presents stress and macro-element
condensation, and compares the resulting system dimensions and sparsity. Local regularity requirements for both condensations
are analyzed in Section~\ref{sec:local-condensation-block-nonsingularity}, which also introduces trace-bound augmentation for
von Mises limit analysis. The GPU implementation of the IPM, including direct solution and the elastoplastic PCG alternative,
is described in Section~\ref{sec:gpu-condensed-solver}, followed by numerical comparisons of the CPU and GPU implementations
in Section~\ref{sec:numerical-results}. Finally, Section~\ref{sec:conclusions} summarizes the findings and discusses future work.


\section{Equilibrium-Based Conic Finite-Element Formulation}
\label{sec:equilibrium-conic-formulation}

This section formulates the equilibrium-based elastoplasticity and lower-bound limit-analysis problems and
introduces their common finite-element discretization. Both problems are then cast in a unified conic form,
from which the primal--dual interior-point Newton system condensed in the next section is derived.

\noindent\textbf{Notation.}
Following common conventions in numerical optimization, vectors, matrices, tensors, and their Voigt representations are
written without boldface; their roles and dimensions are specified when introduced. Boldface is reserved for unit normal vectors.

\subsection{Continuous Problems and Equilibrium Discretization}
\label{subsec:equilibrium-problem-discretization}

Let \(V\subset\mathbb{R}^3\) denote a connected structural domain, with boundary \(\partial V=\Gamma_u\cup\Gamma_t\) and outward unit
normal \(\boldsymbol n\). A displacement \(\bar u\) is prescribed on \(\Gamma_u\), and a reference traction \(\bar t\) is prescribed
on \(\Gamma_t\). The prescribed-displacement conditions are assumed to eliminate all global rigid-body motions. Body forces are omitted
from the displayed problems for simplicity. Using the Cauchy stress tensor \(\sigma\) as the primary unknown, the continuous elastoplastic
problem is written as follows~\cite{ferradi2023elastoplastic}:
\begin{equation}
	\begin{aligned}
		\min_{\sigma}\quad
		& \frac12 \int_V \sigma:\mathbb C^{-1}:\sigma {\,\mathrm dV} -\int_{\Gamma_u}(\sigma \boldsymbol n)\cdot \bar u {\,\mathrm d\Gamma} ,\\
		\text{s.t.}\quad
		& \begin{alignedat}[t]{2}
			& \operatorname{div}\sigma=0  &&\qquad\text{in }V, \\
			& \sigma \boldsymbol n=\bar t &&\qquad\text{on }\Gamma_t, \\
			& f(\sigma)\le0 &&\qquad\text{in }V .
		\end{alignedat}
	\end{aligned}
	\label{eq:continuous-elastoplastic}
\end{equation}
Here \(\mathbb C\) is the symmetric positive-definite elasticity tensor, and \(f(\sigma)\le0\) denotes the material admissibility condition.

In lower-bound limit analysis, the elastic energy is omitted and the objective is to determine the largest
statically admissible load multiplier~\cite{ferradi2023elastoplastic}:
\begin{equation}
	\begin{aligned}
		\max_{\sigma,\lambda}\quad
		& \lambda, \\
		\text{s.t.}\quad
		& \operatorname{div}\sigma=0 &&\text{in }V, \\
		& \sigma \boldsymbol n=\lambda\bar t &&\text{on }\Gamma_t, \\
		& f(\sigma)\le0 &&\text{in }V .
	\end{aligned}
	\label{eq:continuous-limit-analysis}
\end{equation}
The two problems share the same equilibrium and admissibility conditions and differ in their objective and loading.

To construct the equilibrium finite-element discretization common to both problems, the domain \(V\) is conformingly partitioned into
\(n_{\rm vol}\) nondegenerate tetrahedral sub-elements \(V_i\), with \(V=\bigcup_{i=1}^{n_{\rm vol}}V_i\). Following the super-element
construction of Kempeneers et al.~\cite{kempeneers2010pure}, four sub-elements are assembled into each super element to suppress the
spurious kinematic modes exhibited by individual equilibrium elements~\cite{maunder1996general}.

Using the Voigt notation, the stress tensor is represented by
\begin{equation}
	\sigma_v =
	\begin{bmatrix}
		\sigma_{xx} & \sigma_{yy} & \sigma_{zz} & \sigma_{xy} & \sigma_{yz} & \sigma_{xz}
	\end{bmatrix}^\top \in\mathbb{R}^6.
	\label{eq:stress-voigt-vector}
\end{equation}
The inverse Voigt mapping is \(\operatorname{mat}_v:\mathbb{R}^6\to\mathbb{S}^3\), where \(\mathbb{S}^n\)
denotes the space of \(n\times n\) symmetric matrices.

On each sub-element \(V_i\), we use an affine pure-equilibrium stress interpolation whose modes satisfy local
equilibrium by construction~\cite{kempeneers2010pure}:
\begin{equation}
	\sigma_v=S_i q_i, \qquad S_i\in\mathbb{R}^{6\times21}, \qquad q_i\in\mathbb{R}^{21}.
	\label{eq:stress-interpolation}
\end{equation}
The construction extends directly to a nonzero body-force density \(b_V\): a particular stress field \(\sigma^p\) satisfying
\(\operatorname{div}\sigma^p+b_V=0\) can be added to the homogeneous interpolation, giving \(\sigma_v=\sigma_v^p+S_i q_i\).

The stress field is allowed to be discontinuous across sub-element faces, while continuity of the normal traction is enforced through the
face equilibrium equations. On an internal face \(F=\partial V_i\cap\partial V_j\), let \(\sigma_i\) and \(\sigma_j\) denote the stress
traces from the adjacent sub-elements, and let \(\boldsymbol n_{F,i}\) and \(\boldsymbol n_{F,j}\) denote the outward unit normals to
\(V_i\) and \(V_j\), respectively. Defining
\(\left[\!\left[\sigma\right]\!\right]\boldsymbol n:=\sigma_i \boldsymbol n_{F,i}+\sigma_j \boldsymbol n_{F,j}\),
traction continuity requires \(\left[\!\left[\sigma\right]\!\right]\boldsymbol n=0 \;\,\text{on }F\).

The Lagrange multipliers associated with these face equilibrium equations are interpreted as face displacement
fields. With the linear face interpolation represented by \(N_F\in\mathbb{R}^{3\times9}\), each triangular face
carries nine displacement degrees of freedom~\cite{kempeneers2010pure,ferradi2023elastoplastic}.

Since local equilibrium is embedded in the interpolation \eqref{eq:stress-interpolation}, the discrete form of
\eqref{eq:continuous-elastoplastic} imposes traction continuity, prescribed boundary tractions, and vertexwise yield admissibility:
\begin{subequations}
	\label{eq:discrete-elastoplastic-sigma}
	\begin{alignat}{2}
		\min_{\sigma}\quad
		& \makebox[0pt][l]{\(\displaystyle
			\frac12 \int_V \sigma:\mathbb C^{-1}:\sigma {\,\mathrm dV} - \int_{\Gamma_u}(\sigma \boldsymbol n)\cdot \bar u {\,\mathrm d\Gamma}\)},
		\label{eq:discrete-elastoplastic-objective}\\
		\text{s.t.}\quad & \left[\!\left[\sigma\right]\!\right]\boldsymbol n=0 &&\qquad\text{on all internal faces}, \label{eq:discrete-traction-continuity}\\
		& \sigma \boldsymbol n=\bar t &&\qquad\text{on }\Gamma_t, \label{eq:discrete-traction-boundary} \\
		& f(\sigma)\le0 &&\qquad\text{for all vertices of every }V_i. \label{eq:discrete-yield-vertices}
	\end{alignat}
\end{subequations}
The stress interpolation is affine and the yield domains considered here are convex; therefore, enforcing admissibility at the vertices
enforces it throughout each sub-element. The same discrete equilibrium and admissibility conditions apply to the limit-analysis problem
\eqref{eq:continuous-limit-analysis}.

The two convex yield criteria considered below are von Mises and Rankine. Denoting by \(\sigma_y\) the uniaxial yield stress, the von Mises
criterion is written as
\begin{equation}
	\left\|\operatorname{dev}\sigma\right\|_F \le \sqrt{\frac23}\,\sigma_y ,
	\label{eq:vm-continuous}
\end{equation}
where \(\operatorname{dev}\sigma\) denotes the deviatoric part of \(\sigma\), and \(\lVert\cdot\rVert_F\) denotes the Frobenius norm.

Denoting by \(\sigma_{\min}\) and \(\sigma_{\max}\) the minimum and maximum principal stresses, respectively, the Rankine criterion is
\begin{equation}
	-f_c\le\sigma_{\min}, \qquad \sigma_{\max}\le f_t,
	\label{eq:rankine-continuous}
\end{equation}
where \(f_c>0\) and \(f_t>0\) are the compressive and tensile strengths.

To obtain the global algebraic formulation, we collect the coefficient vectors of all sub-elements, define \(n_{\rm dof}:=21n_{\rm vol}\), and set
\begin{equation}
	q = 
	\begin{bmatrix}
		q_1^\top & \cdots & q_{n_{\rm vol}}^\top
	\end{bmatrix}^\top \in\mathbb{R}^{n_{\rm dof}} .
	\label{eq:assembled-stress-vector}
\end{equation}
Let \(n_{\rm face}\) denote the number of triangular faces entering the assembly of equations \eqref{eq:discrete-traction-continuity} and
\eqref{eq:discrete-traction-boundary}. Since each face contributes nine scalar equations, the total number of equality constraints is
\begin{equation}
	n_{\rm eq}=9n_{\rm face}.
	\label{eq:number-equality-constraints}
\end{equation}
The resulting assembled system is \(A_q q=d_q\), where \(A_q\in\mathbb{R}^{n_{\rm eq}\times n_{\rm dof}}\) and
\(d_q\in\mathbb{R}^{n_{\rm eq}}\). The vector \(d_q\) collects the corresponding right-hand sides.

\subsection{Conic Formulation}
\label{subsec:conic-formulation}

For \(a=1,\ldots,4\), let \(S_{i,a}\) denote the interpolation matrix \(S_i\) from \eqref{eq:stress-interpolation} evaluated at vertex \(a\)
of sub-element \(V_i\). In terms of the generalized stress vector, the elastoplastic problem \eqref{eq:discrete-elastoplastic-sigma} becomes
\begin{equation}
	\begin{aligned}
		\min_q\quad
		& \frac12 q^\top P_q q+c_q^\top q, \\
		\text{s.t.}\quad
		& A_q q=d_q, \\
		& f_v(S_{i,a} q_i)\le0, \qquad i=1,\dots,n_{\rm vol},\quad a=1,\dots,4,
	\end{aligned}
	\label{eq:discrete-elastoplastic-q}
\end{equation}
where \(f_v\) denotes the yield function expressed in Voigt stress coordinates.

Denoting by \(\operatorname{blkdiag}(\cdot)\) the block-diagonal concatenation of a sequence of matrices, the
assembled compliance matrix \(P_q\) is given by
\begin{equation}
	P_q = \operatorname{blkdiag}_{i=1}^{n_{\rm vol}}(P_i), \qquad
	P_i=\int_{V_i}S_i^\top\,\mathbb C_i^{-1}\,S_i{\,\mathrm dV} \in\mathbb{R}^{21\times21},
	\label{eq:assembled-compliance}
\end{equation}
where \(\mathbb C_i^{-1}\) is the positive definite compliance matrix on \(V_i\) expressed in the adopted Voigt
convention. Since the stress modes are linearly independent and \(\mathbb C_i^{-1}\succ0\), each \(P_i\succ0\).

The vector \(c_q\) is the assembled linear contribution associated with the prescribed displacement data on \(\Gamma_u\), defined by
\begin{equation*}
	c_q^\top q = -\int_{\Gamma_u}(\sigma \boldsymbol n)\cdot\bar u\,{\,\mathrm d\Gamma} .
\end{equation*}
The equality matrix \(A_q\) is assembled from the local face matrices \(A_{F,i}\). On an internal face, the matrices associated with the
two adjacent sub-elements enter the same block row. Each local face matrix is given by
\begin{equation}
	A_{F,i} = \int_F N_F^\top \mathcal{E}(\boldsymbol n_{F,i})S_i\,{\,\mathrm d\Gamma}
	\in\mathbb{R}^{9\times21}, \qquad F\subset\partial V_i\setminus\Gamma_u,
	\label{eq:local-face-equilibrium-block}
\end{equation}
where \(\mathcal{E}(\boldsymbol n)\) is the stress-to-traction extraction operator defined by
\begin{equation*}
	\mathcal{E}(\boldsymbol n) =
	\begin{bmatrix}
		n_x & 0   & 0   & n_y & 0   & n_z\\
		0   & n_y & 0   & n_x & n_z & 0\\
		0   & 0   & n_z & 0   & n_y & n_x
	\end{bmatrix}, \qquad \sigma \boldsymbol n=\mathcal{E}(\boldsymbol n)\sigma_v .
\end{equation*}
For purely proportional lower-bound limit analysis, let \(\ell\in\mathbb{R}^{n_{\rm eq}}\) denote the assembled reference-load vector.
The equilibrium equations are \(A_q q=\lambda\ell\). Appending \(\lambda\) to the stress variables and defining
\begin{equation}
	A_\lambda =
	\begin{bmatrix}
		A_q & -\ell
	\end{bmatrix} \in\mathbb{R}^{n_{\rm eq}\times(n_{\rm dof}+1)}
	\label{eq:limit-analysis-equality-matrix}
\end{equation}
gives the discrete problem
\begin{equation}
	\begin{aligned}
		\min_{q,\lambda}\quad
		& -\lambda, \\ 
		\text{s.t.}\quad
		& A_\lambda
		\begin{bmatrix}
			q\\ \lambda
		\end{bmatrix} = 0, \\
		& f_v(S_{i,a} q_i)\le0, \qquad i=1,\dots,n_{\rm vol},\quad a=1,\dots,4 .
	\end{aligned}
	\label{eq:discrete-limit-analysis-q}
\end{equation}
To express the vertexwise admissibility conditions in standard conic form, we use the nonnegative orthant,
second-order cone, and positive-semidefinite cone, defined, respectively, by
\begin{equation*}
	\begin{aligned}
		\mathbb{R}_+^m
		&:= \left\{ v\in\mathbb{R}^m : v_j\ge0,\quad j=1,\ldots,m \right\}, \\
		\mathcal Q_m
		&:= \left\{ (t,v)\in\mathbb{R}^m :\; t\in\mathbb{R}_+,\; v\in\mathbb{R}^{m-1},\; t\ge\lVert v\rVert_2 \right\}, \\
		\mathbb{S}^n_+
		&:= \left\{ X\in\mathbb{S}^n : X\succeq0 \right\},
	\end{aligned}
\end{equation*}
where \(X\succeq0\) indicates that \(X\) is positive semidefinite, whereas \(X\succ0\) indicates that \(X\) is positive definite.

For von Mises, let \(D\in\mathbb{R}^{5\times6}\) denote the norm-preserving deviatoric map in the adopted Voigt
convention~\cite{bruno2020return}. Condition \eqref{eq:vm-continuous} becomes
\begin{equation*}
	\left\|DS_{i,a} q_i\right\|_2 \le \sqrt{\frac23}\,\sigma_y.
\end{equation*}
It admits the following local SOC representation~\cite{kanno2011nonsmooth, bisbos2007second}:
\begin{equation*}
	G_{i,a} q_i+s_{i,a} = h_{i,a}, \qquad s_{i,a}\in\mathcal Q_6,
\end{equation*}
with
\begin{equation}
	G_{i,a} =
	\begin{bmatrix}
		0_{1\times21}\\
		DS_{i,a}
	\end{bmatrix}, \qquad
	h_{i,a} =
	\begin{bmatrix}
		\sqrt{2/3}\,\sigma_y\\
		0_5
	\end{bmatrix}.
	\label{eq:von-mises-vertex-data}
\end{equation}
For Rankine, condition \eqref{eq:rankine-continuous} is equivalently expressed by the semidefinite constraints~\cite{larsen2012limit}:
\begin{equation*}
	\operatorname{mat}_v\!\left(S_{i,a} q_i\right) +f_c{\mathbf I}\in\mathbb{S}^3_+, \qquad
	-\operatorname{mat}_v\!\left(S_{i,a} q_i\right) +f_t{\mathbf I}\in\mathbb{S}^3_+.
\end{equation*}
For their algebraic representation, define the symmetric vectorization by
\begin{equation*}
	T_v = \operatorname{diag}\!\left(1,1,1,\sqrt2,\sqrt2,\sqrt2\right), \qquad \operatorname{svec}(X)=T_vX_v,
\end{equation*}
where \(X_v\) is the Voigt vector of \(X\in\mathbb{S}^3\) in the ordering of \eqref{eq:stress-voigt-vector}.
This scaling preserves the inner product:
\[
\operatorname{svec}(X)^\top\operatorname{svec}(Y)=\operatorname{tr}(XY)=X:Y, \qquad X,Y\in\mathbb{S}^3.
\]
The vectorized identity is denoted by
\begin{equation}
	I_v=\operatorname{svec}({\mathbf I}) =
	\begin{bmatrix}
		1&1&1&0&0&0
	\end{bmatrix}^{\top}.
	\label{eq:svec-identity}
\end{equation}
With this convention, the Rankine constraints take the local conic form
\begin{equation*}
	G_{i,a} q_i+s_{i,a} = h_{i,a},
\end{equation*}
where the vertexwise slack is partitioned as
\begin{equation*}
	s_{i,a} =
	\begin{bmatrix}
		s_{i,a}^{\,c}\\
		s_{i,a}^{\,t}
	\end{bmatrix}, \qquad
	s_{i,a}^{\,c},s_{i,a}^{\,t} \in\operatorname{svec}(\mathbb{S}^3_+).
\end{equation*}
Thus, \(s_{i,a}\in\mathbb{R}^{12}\) contains two vectorized \(3\times3\) semidefinite blocks, and the
corresponding local matrix and data vector are
\begin{equation}
	G_{i,a} =
	\begin{bmatrix}
		-T_v S_{i,a}\\
		T_v S_{i,a}
	\end{bmatrix}, \qquad
	h_{i,a} =
	\begin{bmatrix}
		f_c I_v\\
		f_t I_v
	\end{bmatrix}.
	\label{eq:rankine-vertex-data}
\end{equation}
Since each sub-element has four vertices, the von Mises criterion gives a total of \(4n_{\rm vol}\) SOC constraints. The Rankine criterion
gives two \(3\times3\) PSD cones at each vertex, and therefore a total of \(8n_{\rm vol}\) PSD cones.

For either discretized problem, let \(x\) denote the vector of primal variables, and write the assembled equality constraints as \(Ax=b\).
After collecting all vertexwise yield constraints, both problems can be written in the unified conic form
\begin{equation}
	\begin{aligned}
		\min_{x,s}\quad
		& \frac12 x^\top P x+c^\top x, \\
		\text{s.t.}\quad
		& Ax=b, \\
		& Gx+s=h, \qquad s\in\mathcal{K} .
	\end{aligned}
	\label{eq:unified-primal-conic}
\end{equation}
Here \(G\) and \(h\) are assembled using local \(G_{i,a}\) and \(h_{i,a}\), \(s\) denotes the global slack vector, and \(\mathcal{K}\) is
a Cartesian product of the nonnegative orthants, SOCs, and PSD cones present in the problem. These cone families are self-dual under the
adopted inner products; hence, \(\mathcal{K}^*=\mathcal{K}\).

For elastoplasticity, \(x=q\), \(P=P_q\), \(c=c_q\), \(A=A_q\), and \(b=d_q\). For limit analysis, the corresponding conic data are
\begin{equation*}
	x=
	\begin{bmatrix}
		q\\
		\lambda
	\end{bmatrix}, \qquad
	P=0, \qquad c=
	\begin{bmatrix}
		0_{n_{\rm dof}}\\
		-1
	\end{bmatrix}, \qquad
	A=A_\lambda, \qquad b=0_{n_{\rm eq}}.
\end{equation*}
The column of \(G\) associated with \(\lambda\) is zero because the yield constraints act only on the stress coefficients.

\subsection{Primal-Dual Interior-Point Method and Newton System}
\label{subsec:ipm-newton-system}

We briefly recall the standard primal--dual predictor--corrector interior-point framework~\cite{mehrotra1992implementation,bleyer2018advances}
applied to \eqref{eq:unified-primal-conic}, following the conic implementation in CVXOPT~\cite{vandenberghe2010cvxopt}. Here, \(P\succeq0\)
is the objective Hessian, \(y\) is the multiplier associated with \(Ax=b\), and \(z\in\mathcal K\) is the dual variable associated with
\(Gx+s=h\). Using Nesterov--Todd scaling for the complementarity relation, the perturbed KKT conditions are
\begin{equation}
	\begin{aligned}
		Px+c+A^\top y+G^\top z &=0,\\
		Ax-b &=0,\\
		Gx+s-h &=0,\\
		(W^{-\top}s)\circ(Wz) &=\mu e,\\
		s,z&\in\operatorname{int}\mathcal{K}.
	\end{aligned}
	\label{eq:ipm-perturbed-kkt-system}
\end{equation}
Here \(\mu=s^\top z/m\) is the duality measure, where \(m\) is the degree of the product cone \(\mathcal{K}\); \(W\) is the Nesterov--Todd
scaling operator computed from the current cone iterates~\cite{nesterov1997self,nesterov1998primal}; \(\circ\) denotes the Euclidean Jordan
product associated with the cone factors of \(\mathcal{K}\); and \(e\) is the identity element of the corresponding product
algebra~\cite{vandenberghe2010cvxopt}. The scaled cone variable satisfies \(\xi:=Wz=W^{-\top}s\).

At each IPM iteration, the scaled system \eqref{eq:ipm-perturbed-kkt-system} is linearized first for the affine
predictor and then for the combined corrector.

The associated primal and dual residuals are given by
\begin{equation}
	\begin{aligned}
		r_{\rm d} &= Px+c+A^\top y+G^\top z,\\
		r_{\rm e} &= Ax-b,\\
		r_{\rm c} &= Gx+s-h.
	\end{aligned}
	\label{eq:ipm-residuals}
\end{equation}
Convergence is assessed using the normalized primal and dual residuals and the absolute- and relative-gap
criteria adopted in CVXOPT~\cite{vandenberghe2010cvxopt}. Let \(v:=(x,s,y,z)\) denote the current iterate and
\(\Delta:=(\Delta x,\Delta s,\Delta y,\Delta z)\) a search direction. The linearized scaled KKT system is
\begin{equation}
	\begin{aligned}
		P\Delta x+A^\top\Delta y+G^\top\Delta z &=-r_{\rm d},\\
		A\Delta x &=-r_{\rm e},\\
		G\Delta x+\Delta s &=-r_{\rm c},\\
		\xi\circ\left(W\Delta z+W^{-\top}\Delta s\right) &=d_\mu .
	\end{aligned}
	\label{eq:ipm-linearized-kkt-system}
\end{equation}
Only the complementarity right-hand side \(d_\mu\) changes between the two solve phases. For the affine predictor
direction \(\Delta_{\rm a}\), the right-hand side is \(d_\mu^{\rm a} =-\xi\circ\xi\).

To ensure the updated iterates remain within the cone interior, the maximum admissible step length for an
arbitrary direction \(\Delta\) is defined as
\begin{equation}
	\alpha_{\max}(\Delta) :=
	\sup\left\{ \alpha\in[0,1]: \left( s+\alpha\Delta s, z+\alpha\Delta z \right)\in\mathcal{K}\times\mathcal{K}
	\right\}.
	\label{eq:ipm-admissible-step}
\end{equation}
Based on the affine step length \(\alpha_{\rm a}=\alpha_{\max}(\Delta_{\rm a})\), a centering parameter
\(\sigma_{\rm cen}\in[0,1]\) is determined using the predictor-corrector heuristic of~\cite{vandenberghe2010cvxopt}.
This parameter is then used to form the right-hand side for the combined corrector step:
\begin{equation*}
	d_\mu^{\rm c} =-\xi\circ\xi+\sigma_{\rm cen}\mu e -(W^{-\top}\Delta s_{\rm a}) \circ(W\Delta z_{\rm a}),
\end{equation*}
which yields the combined search direction \(\Delta_{\rm c}\).

For both the predictor and corrector steps, the cone directions \(\Delta s\) and \(\Delta z\) are algebraically
eliminated from the linearized system \eqref{eq:ipm-linearized-kkt-system}. This yields the global
symmetric-indefinite Newton system for \((\Delta x,\Delta y)\), together with explicit recovery relations for the cone directions:
\begin{subequations}
	\label{eq:ipm-reduced-system}
	\begin{equation}
		\begin{bmatrix}
			\Omega & A^\top\\
			A      & 0
		\end{bmatrix}
		\begin{bmatrix}
			\Delta x\\
			\Delta y
		\end{bmatrix} =
		\begin{bmatrix}
			b_x\\
			b_y
		\end{bmatrix},
		\label{eq:ipm-newton-system}
	\end{equation}
	\begin{equation}
		\begin{aligned}
			\Delta s &=-r_{\rm c}-G\Delta x,\\
			\Delta z &=W^{-1}\left(\xi\backslash d_\mu\right) +\Theta\left(r_{\rm c}+G\Delta x\right),
		\end{aligned}
		\label{eq:ipm-cone-direction-recovery}
	\end{equation}
\end{subequations}
Here, \(\Theta:=W^{-1}W^{-\top}\) and \(\Omega:=P+G^\top\Theta G\), with \(\xi\backslash w\) denoting the unique
\(u\) satisfying \(\xi\circ u=w\). The right-hand side is
\begin{equation*}
	\begin{aligned}
		b_x &:=-r_{\rm d} -G^\top\left( \Theta r_{\rm c} + W^{-1}\left(\xi\backslash d_\mu\right) \right),\\
		b_y &:=-r_{\rm e}.
	\end{aligned}
\end{equation*}
Since the iterates remain in the cone interior, the Nesterov--Todd scaling operator \(W\) is nonsingular.
Consequently, \(\Theta\succ0\) and \(\Omega\succeq0\).

\begin{assumption}[Global Newton nonsingularity~\cite{vandenberghe2010cvxopt}]
	\label{ass:global-newton-regularity}
	The assembled conic data satisfy
	\begin{equation*}
		\operatorname{rank}(A)=n_{\rm eq}, \qquad \ker(P)\cap\ker(G)\cap\ker(A)=\{0\}.
	\end{equation*}
\end{assumption}
Because \(P\succeq0\) and \(\Theta\succ0\), the definition of \(\Omega\) gives \(\ker(\Omega)=\ker(P)\cap\ker(G)\).
Hence Assumption~\ref{ass:global-newton-regularity} is equivalently written as
\begin{equation*}
	\operatorname{rank}(A)=n_{\rm eq}, \qquad \ker(\Omega)\cap\ker(A)=\{0\},
\end{equation*}
and ensures that the coefficient matrix in \eqref{eq:ipm-newton-system} is nonsingular at every interior-point
iteration~\cite{benzi2005numerical}. While this algebraic property guarantees that the Newton steps are uniquely
defined, it does not imply the existence of a solution to the underlying conic optimization problem.

Within a single IPM iteration, the predictor and corrector systems share the coefficient matrix in
\eqref{eq:ipm-newton-system}, which is therefore factorized once and applied to both right-hand sides. Despite
this reuse, repeated factorization and solution of the Newton system across IPM iterations remain the dominant
computational and memory costs of the method. The next section therefore develops a structure-exploiting
condensation that reduces the size of this global Newton system and, consequently, the cost and memory
requirements of its solution.


\section{Structure-Exploiting Condensation of the Newton System}
\label{sec:structure-exploiting-condensation}

At a fixed IPM iteration, we consider the Newton system \eqref{eq:ipm-newton-system} generated by the
equilibrium finite-element discretization introduced in the preceding section. In this discretization,
four sub-elements form each super element, which is treated here as a macro element. The resulting
local/global partition exposes a bordered block-diagonal structure that permits the macro-local variables
to be eliminated independently, leaving a smaller globally coupled system.

\subsection{Macro Element Partition and Block Structure}
\label{subsec:mechanics-induced-structure}

Denoting the macro elements by \(\mathcal M_e\), the domain partition is written as
\begin{equation*}
	V=\bigcup_{e=1}^{N_m}\mathcal M_e, \qquad \mathcal M_e =\bigcup_{r=1}^{4}V_{e,r},
\end{equation*}
with mutually disjoint interiors. Here, \(V_{e,r}\) denotes sub-element \(r\) of macro element \(e\). Let
\(i(e,r)\) denote its global sub-element index, so that \(V_{e,r}=V_{i(e,r)}\) and \(q_{e,r}=q_{i(e,r)}\),
where \(q_i\) is the \(21\)-component stress coefficient vector in \eqref{eq:stress-interpolation}.
Henceforth, \(i\) indexes sub-elements, whereas \(e\) indexes macro elements. The four stress vectors are collected in
\begin{equation*}
	q_e = \begin{bmatrix}
		q_{e,1}\\
		q_{e,2}\\
		q_{e,3}\\
		q_{e,4} \end{bmatrix} \in\mathbb{R}^{84}.
\end{equation*}
Thus, \(q_e\) is the block of the assembled vector \(q\) in \eqref{eq:assembled-stress-vector} associated
with \(\mathcal M_e\).
\begin{figure}[htbp]
	\centering
	\includegraphics[width=0.40\textwidth]{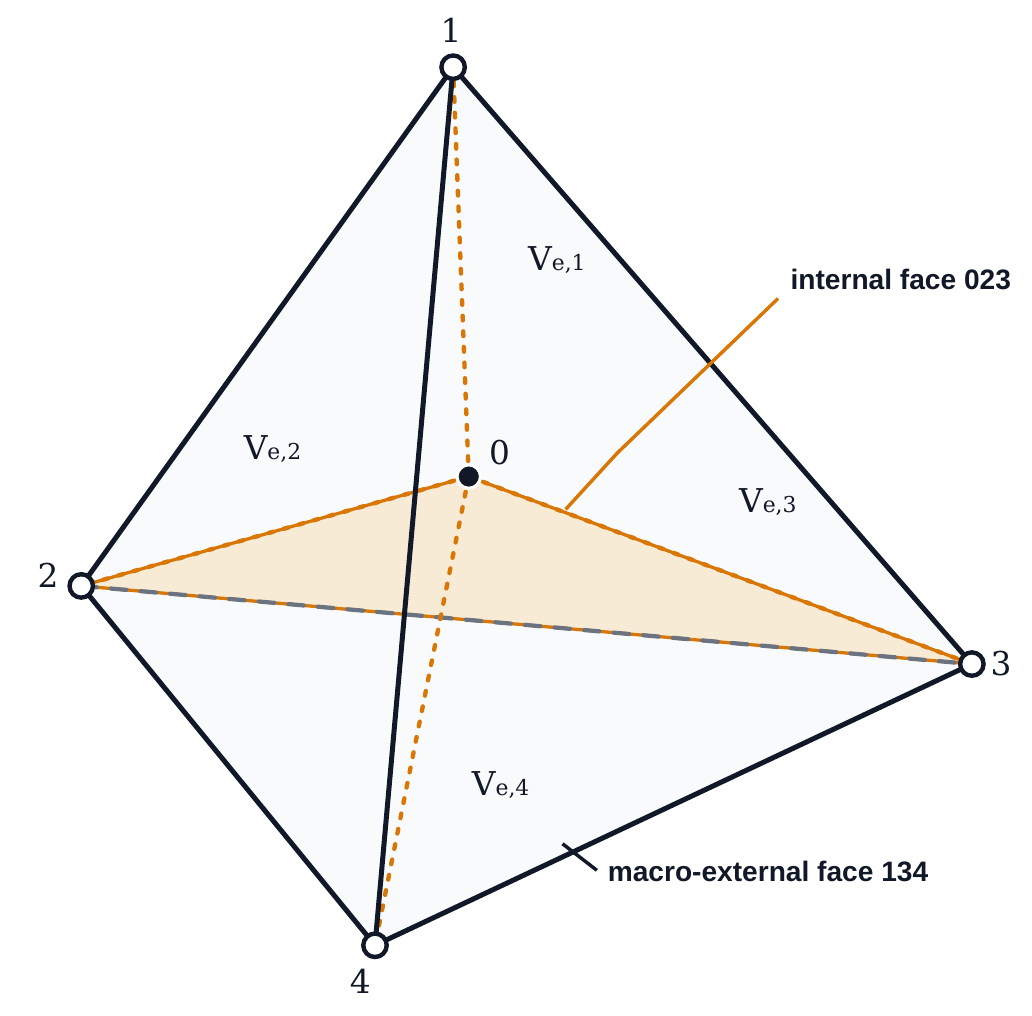}
	\caption{Topology of the four-sub-element macro element \(\mathcal M_e\). Connecting the interior point
		\(0\) to the outer vertices \(1,\ldots,4\) defines the four tetrahedral sub-elements
		\(V_{e,r}\), \(r=1,\ldots,4\). The partition contains 6 macro-internal faces, including the shaded face
		\(023\), and 4 macro-external faces.}
	\label{fig:macro-element-topology}
\end{figure}
As illustrated in Figure~\ref{fig:macro-element-topology}, each macro element has 6 macro-internal faces,
shared by two of its own sub-elements, and 4 macro-external faces. A macro-external face is either shared
with a neighboring macro element or lies on the physical boundary. The macro-element mesh is assumed
face-connected: for any two macro elements, there exists a sequence of macro elements in which each consecutive
pair shares a complete macro-external face. Faces on \(\Gamma_u\) generate no equality rows, whereas faces shared
by neighboring macro elements and faces on \(\Gamma_t\) enter the assembled equality system. This distinction
induces the following partition of the equality-multiplier vector \(y\):
\begin{equation*}
	y= \begin{bmatrix}
		y_{\mathrm I}\\
		y_{\mathrm E} \end{bmatrix}, \qquad
	y_{\mathrm I}= \begin{bmatrix}
		y_{\mathrm I,1}\\
		\vdots\\
		y_{\mathrm I,N_m} \end{bmatrix}.
\end{equation*}
The vector \(y_{\mathrm I,e}\in\mathbb{R}^{54}\) contains the 9 multipliers associated with each of the 6
macro-internal faces of \(\mathcal M_e\). The vector \(y_{\mathrm E}\) contains the multipliers associated with the
contributing macro-external faces. We denote its dimension by \(n_{\mathrm E}:=n_{\rm eq}-54N_m\).

The elastoplastic and lower-bound limit-analysis formulations given in \eqref{eq:discrete-elastoplastic-q}
and \eqref{eq:discrete-limit-analysis-q}, respectively, constitute the two base problems considered in this
work. To accommodate extensions such as reinforced-concrete formulations featuring additional
reinforcement-stress variables~\cite{ferradi2023elastoplastic,ferradi2024steel}, the condensation
is developed for an extended formulation with the primal partition
\begin{equation}
	x= \begin{bmatrix}
		q_1\\
		\vdots\\
		q_{N_m}\\
		x_{\mathrm g} \end{bmatrix},
	\label{eq:primal-macro-global-partition}
\end{equation}
where \(x_{\mathrm g}\in\mathbb{R}^{n_{\mathrm g}}\) collects the remaining primal variables outside the macro-element stress
blocks. This block is absent in the elastoplastic formulation \eqref{eq:discrete-elastoplastic-q}, whereas
\(x_{\mathrm g}=\lambda\) in the lower-bound limit-analysis problem \eqref{eq:discrete-limit-analysis-q}.

We consider extensions in which the variables in \(x_{\mathrm g}\) contribute only to the existing macro-external equilibrium equations,
without introducing additional equality equations. With the row ordering \((y_{\mathrm I},y_{\mathrm E})\) and the column ordering in
\eqref{eq:primal-macro-global-partition}, the equality matrix \(A\) in \eqref{eq:unified-primal-conic} has the block structure
\begin{equation}
	A= \begin{bmatrix}
		A_{\mathrm I,1} &        &               & 0\\
		& \ddots &               & \vdots\\
		&        & A_{\mathrm I,N_m}     & 0\\
		A_{\mathrm E,1} & \cdots & A_{\mathrm E,N_m}     & B_{\mathrm E}
	\end{bmatrix},
	\qquad
	A_{\mathrm g}=
	\begin{bmatrix}0\\B_{\mathrm E}\end{bmatrix}.
	\label{eq:equality-local-global-structure}
\end{equation}
Equivalently, \(A=[\,A_q\;A_{\mathrm g}\,]\), where \(A_{\mathrm g}\in\mathbb{R}^{n_{\rm eq}\times n_{\mathrm g}}\) denotes the last block
column in \eqref{eq:equality-local-global-structure}. Here, \(A_{\mathrm I,e}\in\mathbb{R}^{54\times84}\) contains the macro-internal
traction equations, while \(A_{\mathrm E,e}\in\mathbb{R}^{n_{\mathrm E}\times84}\) contains the contribution of macro element \(e\) to
the macro-external equations. Although written in the global macro-external space, \(A_{\mathrm E,e}\) has nonzero rows only for the
contributing macro-external faces of \(\mathcal M_e\). Thus, an interior macro element contributes \(36\) such rows, whereas a macro
element with one macro-external face on \(\Gamma_u\) contributes \(27\). The block \(B_{\mathrm E}\in\mathbb{R}^{n_{\mathrm E}\times n_{\mathrm g}}\)
collects the contributions of \(x_{\mathrm g}\) to the macro-external equilibrium equations. For the standard elastoplastic formulation,
the last block column is absent (\(n_{\mathrm g}=0\)), and the equality matrix reduces to \(A=A_q\). For standard lower-bound limit analysis,
the reference loading acts only through the prescribed-traction equations. Hence, under the \(\mathrm I/\mathrm E\) row partition,
\(\ell=[\,0_{54N_m}^\top\ \ell_{\mathrm E}^\top\,]^\top\), so that \(B_{\mathrm E}=-\ell_{\mathrm E}\), where \(\ell_{\mathrm E}\in\mathbb{R}^{n_{\mathrm E}}\)
has nonzero entries only in rows associated with faces contained in \(\Gamma_t\).

We next describe the macro-local part of the primal block \(\Omega=P+G^\top\Theta G\). For each sub-element
\(V_i\), collect the four vertex blocks \(G_{i,a}\) defined in \eqref{eq:von-mises-vertex-data} or \eqref{eq:rankine-vertex-data} as:
\begin{equation*}
	G_i := \begin{bmatrix}
		G_{i,1}^\top & G_{i,2}^\top &
		G_{i,3}^\top & G_{i,4}^\top \end{bmatrix}^\top.
\end{equation*}
Let \(P_i\in\mathbb{R}^{21\times21}\) and \(\Theta_i\succ0\) denote the corresponding sub-element blocks of
\(P\) and \(\Theta\), respectively. For each macro element \(\mathcal M_e\), define
\begin{equation}
	\begin{aligned}
		\mathcal P_e &:= \operatorname{blkdiag} \left( P_{i(e,1)},\ldots,P_{i(e,4)} \right),\\
		\mathcal G_e &:= \operatorname{blkdiag} \left( G_{i(e,1)},\ldots,G_{i(e,4)} \right),\\
		\Psi_e &:= \operatorname{blkdiag} \left( \Theta_{i(e,1)},\ldots,\Theta_{i(e,4)} \right).
	\end{aligned}
	\label{eq:macro-primal-components}
\end{equation}
The corresponding macro element primal block is
\begin{equation}
		\Omega_e := \mathcal P_e+\mathcal G_e^\top\Psi_e\mathcal G_e
		= \operatorname{blkdiag}_{r=1}^{4} \left( P_{i(e,r)}+G_{i(e,r)}^\top\Theta_{i(e,r)}G_{i(e,r)}
		\right) \in\mathbb{R}^{84\times84}.
	\label{eq:macro-primal-block}
\end{equation}
We restrict attention to extended formulations with no mixed compliance or yield terms between \(q_e\) and \(x_{\mathrm g}\). Hence,
the corresponding off-diagonal blocks of \(\Omega\) vanish. We denote by \(\Omega_{\mathrm g}\in\mathbb{S}^{n_{\mathrm g}}_+\) the
positive semidefinite block associated with \(x_{\mathrm g}\), without prescribing any further internal structure. Coupling with the
macro element stress variables enters through the equality matrix. Consequently,
\begin{equation}
\Omega=
\begin{bmatrix}
\Omega_q & 0\\
0 & \Omega_{\mathrm g}
\end{bmatrix}
=
\begin{bmatrix}
\Omega_1 & & & 0\\
& \ddots & & \vdots\\
& & \Omega_{N_m} & 0\\
0 & \cdots & 0 & \Omega_{\mathrm g}
\end{bmatrix},
\label{eq:omega-local-retained-partition}
\end{equation}
where \(\Omega_q\in\mathbb S_+^{84N_m}\) denotes the symmetric
positive semidefinite stress block.

The block row and column corresponding to \(x_{\mathrm g}\) in \eqref{eq:omega-local-retained-partition} are absent
for the standard elastoplastic formulation \eqref{eq:discrete-elastoplastic-q}. In standard lower-bound
limit analysis \eqref{eq:discrete-limit-analysis-q}, \(x_{\mathrm g}=\lambda\), \(P=0\), and the \(\lambda\)-column
of \(G\) is zero; hence, \(\Omega_{\mathrm g}=0\). In other extended formulations, \(\Omega_{\mathrm g}\) has a problem-dependent
positive-semidefinite structure determined by the roles of the variables in \(x_{\mathrm g}\) in the objective and cone
constraints. For example, when \(x_{\mathrm g}\) includes reinforcement stresses, their corresponding contributions to
the primal Newton block are contained in \(\Omega_{\mathrm g}\).

The corresponding Newton directions are partitioned as
\begin{equation}
	\Delta x= \begin{bmatrix}
		\Delta q_1\\
		\vdots\\
		\Delta q_{N_m}\\
		\Delta x_{\mathrm g} \end{bmatrix}, \qquad
	\Delta y= \begin{bmatrix}
		\Delta y_{\mathrm I,1}\\
		\vdots\\
		\Delta y_{\mathrm I,N_m}\\
		\Delta y_{\mathrm E} \end{bmatrix}.
	\label{eq:newton-direction-split}
\end{equation}
We collect the macro-local and globally coupled directions as
\begin{equation}
	w_e= \begin{bmatrix}
		\Delta q_e\\
		\Delta y_{\mathrm I,e} \end{bmatrix}\in\mathbb{R}^{138}, \qquad
	\Delta\chi= \begin{bmatrix}
		\Delta y_{\mathrm E}\\
		\Delta x_{\mathrm g} \end{bmatrix}\in\mathbb{R}^{n_{\mathrm{cond}}},
	\label{eq:local-retained-variables}
\end{equation}
where \(n_{\mathrm{cond}}:=n_{\mathrm E}+n_{\mathrm g}\). Thus, the macro-local directions \(w_e\) couple to the retained
directions only through \(\Delta y_{\mathrm E}\), providing the block structure required for their independent elimination.

\subsection{Rank Properties and Global Newton Nonsingularity}
\label{subsec:equilibrium-rank-global-nonsingularity}

Before deriving this elimination, we establish the rank properties of the equilibrium operators and the resulting nonsingularity conditions
for the global Newton system \eqref{eq:ipm-newton-system}. The assembled equilibrium matrix has the following structural rank property.

\begin{proposition}[Rank of the assembled equilibrium matrix]
	\label{prop:assembled-equilibrium-full-rank}
	Under the prescribed-displacement condition stated in Subsection~\ref{subsec:equilibrium-problem-discretization},
	the assembled equilibrium matrix \(A_q\) has full row rank, with
	\(\operatorname{rank}(A_q)=n_{\rm eq}\); equivalently, \(\ker(A_q^\top)=\{0\}\).
\end{proposition}
\begin{proof}
	See Appendix~\ref{app:proof-assembled-equilibrium-full-rank}.
\end{proof}

\begin{corollary}[Rank preservation under primal extension]
	\label{cor:augmented-equilibrium-full-rank}
	The assembled equality matrix \(A=[\,A_q\;A_{\mathrm g}\,]\) of the extended formulation \eqref{eq:unified-primal-conic}
	has full row rank, with \(\operatorname{rank}(A)=n_{\rm eq}\); equivalently, \(\ker(A^\top)=\{0\}\).
\end{corollary}
\begin{proof}
	If \(v\in\ker(A^\top)\), then \(A_q^\top v=0\). Proposition~\ref{prop:assembled-equilibrium-full-rank}
	therefore gives \(v=0\).
\end{proof}

The rank of the macro-internal connection blocks, needed for the macro-element condensation, follows
immediately from the full rank of \(A_q\).

\begin{corollary}[Rank of the macro-internal connection matrix]
	\label{cor:internal-connection-rank}
	For every macro element \(\mathcal M_e\), the macro-internal connection matrix \(A_{\mathrm I,e}\) has full row rank,
	with \(\operatorname{rank}(A_{\mathrm I,e})=54\); equivalently, \(\ker(A_{\mathrm I,e}^\top)=\{0\}\).
\end{corollary}

\begin{proof}
	Let \(v_e\in\ker(A_{\mathrm I,e}^\top)\), and form a global multiplier vector \(\widetilde v\in\mathbb{R}^{n_{\rm eq}}\)
	by placing \(v_e\) in the internal rows of \(\mathcal M_e\) and zeros in all other rows. These rows have
	support only on \(q_e\), so the block structure \eqref{eq:equality-local-global-structure} gives
	\(A_q^\top\widetilde v=0\). Proposition~\ref{prop:assembled-equilibrium-full-rank} therefore implies
	\(\widetilde v=0\), and hence \(v_e=0\).
\end{proof}

The following standard criterion links these rank properties and kernel conditions directly to matrix nonsingularity.

\begin{theorem}[Saddle-point nonsingularity criterion~{\cite[Theorem~3.2]{benzi2005numerical}}]
	\label{thm:saddle-point-nonsingularity}
	Let \(H\in\mathbb{S}^n_+\), and suppose \(B\in\mathbb{R}^{m\times n}\) has full row rank. Then
	\begin{equation*}
		\begin{bmatrix}
			H & B^\top\\
			B & 0
		\end{bmatrix} \text{ is nonsingular} \quad\Longleftrightarrow\quad \ker(H)\cap\ker(B)=\{0\}.
	\end{equation*}
	In particular, under the stated rank condition on \(B\), the saddle-point matrix is nonsingular whenever \(H\succ0\).
\end{theorem}

Corollary~\ref{cor:augmented-equilibrium-full-rank} shows that adding \(A_{\mathrm g}\) preserves the full row rank of
the equilibrium matrix. Unlike the full-row-rank property, the kernel-intersection condition is not preserved
automatically when the primal space is extended and depends on the problem-specific blocks \(\Omega_{\mathrm g}\) and
\(A_{\mathrm g}\). Accordingly, Theorem~\ref{thm:saddle-point-nonsingularity} identifies the remaining requirement for
global Newton nonsingularity, which we state as follows.

\begin{assumption}[Non-degeneracy of the extended formulation]
	\label{ass:column-augmentation-regularity}
	For any extended formulation considered here, we require Assumption~\ref{ass:global-newton-regularity} to hold;
	equivalently, \(\ker(\Omega)\cap\ker(A)=\{0\}\).
\end{assumption}

We now specialize the global regularity requirement to the base formulations. For elastoplasticity, the
elastic compliance ensures \(P\succ0\). Consequently, \(\Omega\succ0\) and the kernel-intersection condition
is satisfied unconditionally. For limit analysis, however, \(P=0\) and \(\Omega\) is only positive semidefinite.
Preventing its kernel from intersecting \(\ker(A)\) requires an additional admissibility condition on the reference loading.

\begin{assumption}[Reference-loading admissibility]
	\label{ass:reference-loading-regularity}
	For Rankine limit analysis, the reference load vector satisfies \(\ell\ne0\). For von Mises limit analysis,
	let \(q_h\ne0\) denote the assembled stress-coefficient vector representing the uniform unit hydrostatic stress
	field. The reference load vector \(\ell\) and the resulting traction pattern \(A_q q_h\) are required to be
	linearly independent.
\end{assumption}

\begin{proposition}[Global Newton nonsingularity of the base formulations]
	\label{prop:equilibrium-formulations-regularity}
	For each of the Rankine and von Mises yield criteria, the elastoplasticity formulation
	\eqref{eq:discrete-elastoplastic-q} satisfies Assumption~\ref{ass:global-newton-regularity}. Under
	Assumption~\ref{ass:reference-loading-regularity}, the corresponding limit-analysis formulation
	\eqref{eq:discrete-limit-analysis-q} also satisfies Assumption~\ref{ass:global-newton-regularity}.
\end{proposition}
\begin{proof}
	See Appendix~\ref{app:proof-equilibrium-formulations-regularity}.
\end{proof}

The preceding proposition establishes global Newton nonsingularity for the two base formulations under the stated conditions. For the
reinforced-concrete extension, Assumption~\ref{ass:column-augmentation-regularity} is verified in Section~\ref{sec:numerical-results}.

Within this globally regular setting, a stress condensation of the Newton system is available whenever its local stress blocks are
invertible. This reduction is well established in lower-bound limit analysis, where small diagonal perturbations have been used when
necessary to enforce local invertibility \cite{lyamin2002lower,krabbenhoft2003general}. No such diagonal regularization is introduced in
the mathematical formulation considered here. In limit analysis, the positive definiteness of \(\Omega_q\) depends on the yield criterion.
Section~\ref{sec:local-condensation-block-nonsingularity} establishes this property for Rankine and, through the trace augmentation
introduced there, for von Mises. Since each \(\Omega_e\) is itself block diagonal over its four sub-elements, applications of \(\Omega_q^{-1}\)
reduce to independent solves with \(21\times21\) positive-definite blocks.

If \(\Omega_q\succ0\), partition \(b_x=[\,b_q^\top\;b_{x,\mathrm g}^\top\,]^\top\), where \(b_q\) and \(b_{x,\mathrm g}\) are the components
associated with \(q\) and \(x_{\mathrm g}\), respectively. Eliminating the stress directions from~\eqref{eq:ipm-newton-system} gives
\begin{equation}
	\begin{bmatrix}
		-\Sigma_{\mathrm s} & A_{\mathrm g}\\
		A_{\mathrm g}^\top & \Omega_{\mathrm g}
	\end{bmatrix}
	\begin{bmatrix}
		\Delta y\\
		\Delta x_{\mathrm g}
	\end{bmatrix}
	=
	\begin{bmatrix}
		b_y-A_q\Omega_q^{-1}b_q\\
		b_{x,\mathrm g}
		\end{bmatrix},
		\label{eq:stress-only-condensed-system}
	\end{equation}
	where \(\Sigma_{\mathrm s}:=A_q\Omega_q^{-1}A_q^\top\succ0\). The eliminated stress directions are recovered from
	\begin{equation*}
		\Delta q=
		\Omega_q^{-1}\left(b_q-A_q^\top\Delta y\right).
\end{equation*}
The coefficient matrix in \eqref{eq:stress-only-condensed-system} is the Schur complement of \(\Omega_q\) in the global Newton matrix and
is therefore nonsingular under the global regularity conditions established above.

\subsection{Macro-Element Condensation}
\label{subsec:condensed-newton-system}

For macro element \(e\), define the local saddle-point block
\begin{equation}
	K_e= \begin{bmatrix}
		\Omega_e & A_{\mathrm I,e}^\top\\
		A_{\mathrm I,e}  & 0 \end{bmatrix} \in\mathbb{R}^{138\times138},
	\label{eq:local-saddle-point-block}
\end{equation}
and its coupling with the macro-external multiplier direction
\(\Delta y_{\mathrm E}\) by
\begin{equation}
	C_e:=
	\begin{bmatrix}
		A_{\mathrm E,e}^\top\\
		0_{54\times n_{\mathrm E}}
	\end{bmatrix}
	\in\mathbb{R}^{138\times n_{\mathrm E}}.
	\label{eq:macro-external-coupling-block}
\end{equation}
Partitioning \(b_x\) and \(b_y\) according to \eqref{eq:newton-direction-split}, we define
\(b_e:=[\,b_{x,e}^{\top}\;b_{y,\mathrm I,e}^{\top}\,]^{\top}\in\mathbb R^{138}\). The Newton system \eqref{eq:ipm-newton-system} then reads
\begin{equation}
	\begin{aligned}
	K_e w_e+C_e\Delta y_{\mathrm E}
	&=b_e,
	&&e=1,\ldots,N_m,\\
	\sum_{e=1}^{N_m}C_e^\top w_e
	+B_{\mathrm E}\Delta x_{\mathrm g}
	&=b_{y,\mathrm E},\\
	B_{\mathrm E}^{\top}\Delta y_{\mathrm E}
	+\Omega_{\mathrm g}\Delta x_{\mathrm g}
	&=b_{x,\mathrm g}.
	\end{aligned}
	\label{eq:permuted-newton-system}
\end{equation}
When \(x_{\mathrm g}\) is absent, the last equation and the terms involving \(\Delta x_{\mathrm g}\) are omitted.

\begin{assumption}[Local block nonsingularity]
	\label{ass:local-block-regularity}
	Every macro-local saddle-point block \(K_e\) in \eqref{eq:local-saddle-point-block} is nonsingular.
\end{assumption}
For elastoplasticity, \(\Omega_e\succ0\), together with the established full row rank of \(A_{\mathrm I,e}\), ensures
this assumption. In limit analysis, its validity depends on the yield criterion. Section~\ref{sec:local-condensation-block-nonsingularity}
proves that it holds for Rankine, identifies the hydrostatic null mode that makes the original von Mises blocks singular, and introduces
a trace augmentation that ensures local nonsingularity.

Under Assumption~\ref{ass:local-block-regularity}, the first equation of \eqref{eq:permuted-newton-system} gives
\begin{equation}
	w_e=K_e^{-1}\left(b_e-C_e\Delta y_{\mathrm E}\right),
	\qquad e=1,\ldots,N_m.
	\label{eq:local-solution}
\end{equation}
Substitution into the macro-external equilibrium row introduces the
condensed matrix and right-hand side
\begin{equation}
	\Sigma_{\mathrm E}:=
	\sum_{e=1}^{N_m}C_e^\top K_e^{-1}C_e,\qquad
	\widehat\beta_{\mathrm E}:=
	-b_{y,\mathrm E}+\sum_{e=1}^{N_m}C_e^\top K_e^{-1}b_e.
	\label{eq:macro-external-schur-complement}
\end{equation}
Together with the retained primal equations, this gives the condensed system
\begin{equation}
	\Sigma\Delta\chi=\widehat\beta,
	\label{eq:compact-condensed-system}
\end{equation}
where
\begin{equation}
	\Sigma=
	\begin{bmatrix}
		\Sigma_{\mathrm E} & -B_{\mathrm E}\\
		-B_{\mathrm E}^\top & -\Omega_{\mathrm g}
	\end{bmatrix},
	\qquad
	\widehat\beta=
	\begin{bmatrix}
		\widehat\beta_{\mathrm E}\\
		-b_{x,\mathrm g}
	\end{bmatrix}.
	\label{eq:compact-schur-complement}
\end{equation}

\begin{proposition}[Nonsingularity of the condensed system for the extended formulation]
	\label{prop:general-condensed-system-nonsingular}
	Under Assumption~\ref{ass:local-block-regularity}, the condensed matrix \(\Sigma\) in \eqref{eq:compact-schur-complement} is
	nonsingular if and only if Assumption~\ref{ass:column-augmentation-regularity} holds.
\end{proposition}

\begin{proof}
	Eliminating the nonsingular blocks \(K_e\) from \eqref{eq:permuted-newton-system} gives the Schur complement \(-\Sigma\), so the global
	Newton matrix is nonsingular if and only if \(\Sigma\) is nonsingular. By Corollary~\ref{cor:augmented-equilibrium-full-rank} and
	Theorem~\ref{thm:saddle-point-nonsingularity}, this is equivalent to Assumption~\ref{ass:column-augmentation-regularity}.
\end{proof}

For the standard elastoplastic formulation, \(x_{\mathrm g}\) is absent, and \eqref{eq:compact-condensed-system} reduces to
\begin{equation}
	\Sigma_{\mathrm E}\Delta y_{\mathrm E}=\widehat\beta_{\mathrm E}.
	\label{eq:elastoplastic-condensed-system}
\end{equation}
In lower-bound limit analysis, \(x_{\mathrm g}=\lambda\) and
\(\Omega_{\mathrm g}=0\). Denoting by \(b_\lambda\) the component of \(b_x\) associated with \(\lambda\) and using
\(B_{\mathrm E}=-\ell_{\mathrm E}\), the condensed system is
\begin{equation}
	\begin{bmatrix}
		\Sigma_{\mathrm E}         & \ell_{\mathrm E}\\
		\ell_{\mathrm E}^\top & 0
	\end{bmatrix}
	\begin{bmatrix}
		\Delta y_{\mathrm E}\\
		\Delta\lambda
	\end{bmatrix}
	=
	\begin{bmatrix}
		\widehat\beta_{\mathrm E}\\
		-b_\lambda
	\end{bmatrix}.
	\label{eq:limit-analysis-condensed-system}
\end{equation}
The load multiplier therefore adds one row and column without modifying the local blocks \(K_e\) or their factorizations. For both base
formulations, the common macro-external block \(\Sigma_{\mathrm E}\) satisfies the stronger positive-definiteness property stated next.

\begin{proposition}[Positive definiteness of the macro-external condensed matrix]
	\label{prop:external-schur-positive-definite}
	For the two base formulations, Assumption~\ref{ass:local-block-regularity} implies that the matrix \(\Sigma_{\mathrm E}\)
	defined in \eqref{eq:macro-external-schur-complement} is symmetric positive definite.
\end{proposition}

\begin{proof}
	Since \(\Omega_e\succeq0\) and \(A_{\mathrm I,e}\) has full row rank by Corollary~\ref{cor:internal-connection-rank},
	Assumption~\ref{ass:local-block-regularity} and Theorem~\ref{thm:saddle-point-nonsingularity} give
	\(\ker(\Omega_e)\cap\ker(A_{\mathrm I,e})=\{0\}\). If \(q\in\ker(\Omega_q)\cap\ker(A_q)\), then
	\(\Omega_e q_e=0\) and \(A_{\mathrm I,e}q_e=0\) for every \(e\). The local kernel condition therefore gives
	\(q_e=0\) for every macro element, and hence \(q=0\). Together with
	Proposition~\ref{prop:assembled-equilibrium-full-rank}, the same theorem shows that the global
	stress--multiplier saddle-point matrix defined by \(\Omega_q\) and \(A_q\) is nonsingular. After the
	local/external permutation, its leading diagonal block is \(\operatorname{blkdiag}(K_1,\ldots,K_{N_m})\),
	and its exact Schur complement is \(-\Sigma_{\mathrm E}\). Therefore, \(\Sigma_{\mathrm E}\) is nonsingular.
	
	Since every \(K_e\) is symmetric and nonsingular, \(\Sigma_{\mathrm E}\) is symmetric. It
	remains to determine its sign. For \(v\in\mathbb{R}^{n_{\mathrm E}}\), define
	\begin{equation*}
		\begin{bmatrix}
			\zeta_e\\
			\eta_e
		\end{bmatrix} :=K_e^{-1}C_e v.
	\end{equation*}
	The second local block equation gives \(A_{\mathrm I,e}\zeta_e=0\). Consequently,
	\begin{equation*}
		v^\top \Sigma_{\mathrm E}v = \sum_{e=1}^{N_m}
		\begin{bmatrix}
			\zeta_e\\ \eta_e
		\end{bmatrix}^{\!\top}
		K_e
		\begin{bmatrix}
			\zeta_e\\ \eta_e
		\end{bmatrix}
		= \sum_{e=1}^{N_m}\zeta_e^\top\Omega_e\zeta_e \ge0.
	\end{equation*}
	Hence, \(\Sigma_{\mathrm E}\succeq0\). Its nonsingularity then implies \(\Sigma_{\mathrm E}\succ0\).
\end{proof}

The elastoplastic system~\eqref{eq:elastoplastic-condensed-system} therefore has a symmetric positive-definite coefficient matrix.
For limit analysis, the decomposition \(\ell=[\,0_{54N_m}^\top\ \ell_{\mathrm E}^\top\,]^\top\) and the nonzero reference loading
give \(\ell_{\mathrm E}\ne0\). The Schur complement of \(\Sigma_{\mathrm E}\) in the coefficient matrix of
\eqref{eq:limit-analysis-condensed-system} is the negative scalar \(-\ell_{\mathrm E}^\top \Sigma_{\mathrm E}^{-1}\ell_{\mathrm E}<0.\)
Consequently, the limit-analysis coefficient matrix is symmetric, nonsingular, and indefinite.

\begin{remark}[Positive-definite factorization for limit analysis]
	\label{rem:positive-definite-limit-analysis-solve}
	Although the bordered matrix in \eqref{eq:limit-analysis-condensed-system} is indefinite, it need not be factorized directly.
	With \(u:=\Sigma_{\mathrm E}^{-1}\widehat\beta_{\mathrm E}\) and \(v:=\Sigma_{\mathrm E}^{-1}\ell_{\mathrm E}\), we have
	\(\ell_{\mathrm E}^{\top}v=\ell_{\mathrm E}^{\top}\Sigma_{\mathrm E}^{-1}\ell_{\mathrm E}>0\) since \(\Sigma_{\mathrm E}\succ0\) and
	\(\ell_{\mathrm E}\ne0\). Block elimination then gives \(\Delta\lambda=(\ell_{\mathrm E}^{\top}u+b_\lambda)/(\ell_{\mathrm E}^{\top}v)\)
	and \(\Delta y_{\mathrm E}=u-v\Delta\lambda\). Thus, one factorization of the positive-definite matrix \(\Sigma_{\mathrm E}\), reused
	for the two right-hand sides \(\widehat\beta_{\mathrm E}\) and \(\ell_{\mathrm E}\), can replace a factorization of the full indefinite
	matrix.
\end{remark}

Macro-element condensation further reduces the dimension of the global Newton system while preserving or improving the spectral
conditioning of the positive-definite multiplier matrix.

\begin{proposition}[Conditioning of the macro-external condensed matrix]
	\label{prop:macro-external-conditioning}
	For the two base formulations, suppose \(\Omega_q\succ0\). The matrices \(\Sigma_{\mathrm s}\) and \(\Sigma_{\mathrm E}\),
	defined in \eqref{eq:stress-only-condensed-system} and \eqref{eq:macro-external-schur-complement}, respectively, satisfy
	\[\kappa_2(\Sigma_{\mathrm E})\leq\kappa_2(\Sigma_{\mathrm s}),\]
	where \(\kappa_2(M):=\lVert M\rVert_2\lVert M^{-1}\rVert_2\) denotes the spectral condition number of a nonsingular matrix \(M\).
\end{proposition}

\begin{proof}
	Eliminating the macro-internal multiplier directions from the stress-condensed system identifies \(\Sigma_{\mathrm E}\) as the
	Schur complement of the internal-multiplier block of \(\Sigma_{\mathrm s}\). The block-inverse identity therefore shows that
	\(\Sigma_{\mathrm E}^{-1}\) is a principal submatrix of \(\Sigma_{\mathrm s}^{-1}\). Since \(\Sigma_{\mathrm s}\succ0\), the
	eigenvalue bounds for positive-definite principal submatrices give \(\kappa_2(\Sigma_{\mathrm E}^{-1})\leq\kappa_2(\Sigma_{\mathrm s}^{-1})\).
	The condition number is unchanged by inversion, so the claim follows.
\end{proof}

The terms \(C_e^\top K_e^{-1}C_e\) and \(C_e^\top K_e^{-1}b_e\) provide the macro-local contributions to \(\Sigma_{\mathrm E}\) and
\(\widehat\beta_{\mathrm E}\) in \eqref{eq:macro-external-schur-complement}, respectively. After solving the condensed system,
the local directions \(w_e\) are recovered from \eqref{eq:local-solution}. Both the formation and recovery stages consist
of independent macro-element operations, exposing the coarse-grained parallelism exploited by the GPU implementation presented later.
Together, the recovered directions \(w_e\) and \(\Delta\chi\) yield the full Newton step \((\Delta x,\Delta y)\).

Both the stress condensation introduced above and the macro-element condensation avoid direct factorization of the original Newton matrix
\eqref{eq:ipm-newton-system}, which is a symmetric-indefinite saddle-point matrix with a zero multiplier block of dimension
\(n_{\rm eq}\). Direct factorization of this matrix generally requires a pivoted \(LDL^\top\) factorization, and its zero diagonal block can
increase pivoting and fill \cite{duff1983multifrontal,benzi2005numerical}. Stress condensation eliminates the stress directions but
retains all traction-continuity multipliers in the global system. The macro-element condensation additionally eliminates the
macro-internal multipliers, leaving only \(\Delta\chi=(\Delta y_{\mathrm E},\Delta x_{\mathrm g})\) in the global system.

Although Schur-complement condensation can, in general, introduce more fill-in, the induced couplings in the two base formulations
are confined to macro-external faces incident on the same macro element. The resulting matrix therefore remains sparse while having a
substantially smaller dimension than both the original and stress-condensed systems. The corresponding reductions in dimension and
structural nonzeros generally lower the memory requirements and factorization cost of the global solve and improve its suitability for
GPU-resident factorization.

\subsection{Dimension and Sparsity Estimates}
\label{subsec:sparsity-size-parallelism}

To quantify the impact of the proposed macro-element condensation, this subsection compares the dimensions and structural-nonzero counts
of the original Newton system, the stress-condensed system \eqref{eq:stress-only-condensed-system}, and the macro-condensed system
\eqref{eq:compact-condensed-system}. Physical-boundary contributions are retained in the exact dimension relations but neglected in the
leading-order estimates. The two base formulations have \(n_{\mathrm g}=0\) and \(n_{\mathrm g}=1\), respectively, a difference
that does not affect these estimates.

Let \(n_{\mathrm{full}}\) and \(n_{\mathrm{stress}}\) denote the dimensions of the original and stress-condensed Newton systems,
respectively. Using the local/global partition established in \eqref{eq:newton-direction-split}--\eqref{eq:local-retained-variables},
the three systems have the exact dimensions
\begin{equation}
	n_{\mathrm{full}}=84N_m+n_{\rm eq}+n_{\mathrm g},\qquad
	n_{\mathrm{stress}}=n_{\rm eq}+n_{\mathrm g},\qquad
	n_{\mathrm{cond}}=n_{\rm eq}-54N_m+n_{\mathrm g}.
	\label{eq:full-condensed-dimensions}
\end{equation}
Thus, stress condensation removes the \(84N_m\) stress directions, whereas macro-element condensation additionally removes the \(54N_m\)
macro-internal multipliers and therefore eliminates a total of \(138N_m\) unknowns from the globally coupled system.

An estimate of \(n_{\rm eq}\) follows directly from face counting. Let \(N_{\mathrm{sh}}^{\mathrm f}\) be the number of faces
shared by two macro elements, and let \(N_{\Gamma_t}^{\mathrm f}\) and \(N_{\Gamma_u}^{\mathrm f}\) be the numbers of
physical boundary faces on \(\Gamma_t\) and \(\Gamma_u\), respectively. Counting the four outer-face incidences of each macro element gives
\begin{equation*}
	4N_m = 2N_{\mathrm{sh}}^{\mathrm f}+N_{\Gamma_t}^{\mathrm f}+N_{\Gamma_u}^{\mathrm f}.
\end{equation*}
Combining \eqref{eq:number-equality-constraints} with this relation and the contributing-face convention in \eqref{eq:local-face-equilibrium-block} gives
\begin{equation*}
	n_{\rm eq} = 9\left(6N_m+N_{\mathrm{sh}}^{\mathrm f}+N_{\Gamma_t}^{\mathrm f}\right)
	= 72N_m +\frac{9}{2} \left(N_{\Gamma_t}^{\mathrm f}-N_{\Gamma_u}^{\mathrm f}\right).
\end{equation*}
Neglecting the physical-boundary contribution gives
\begin{equation}
	n_{\rm eq}\approx72N_m .
	\label{eq:equality-count-estimate}
\end{equation}
Equations \eqref{eq:full-condensed-dimensions} and \eqref{eq:equality-count-estimate} yield
\begin{equation*}
	n_{\mathrm{full}}\approx156N_m, \qquad
	n_{\mathrm{stress}}\approx72N_m, \qquad
	n_{\mathrm{cond}}\approx18N_m.
\end{equation*}
Consequently,
\begin{equation}
	\frac{n_{\mathrm{full}}}{n_{\mathrm{stress}}}
	\approx\frac{13}{6}\approx2.17, \qquad
	\frac{n_{\mathrm{stress}}}{n_{\mathrm{cond}}}
	\approx4, \qquad
	\frac{n_{\mathrm{full}}}{n_{\mathrm{cond}}}
	\approx\frac{26}{3}\approx8.67.
	\label{eq:dimension-reduction-ratios}
\end{equation}
The estimates involving stress condensation apply whenever the local stress blocks are invertible.
Section~\ref{sec:local-condensation-block-nonsingularity} shows that this requirement is more restrictive than
the macro-local nonsingularity condition used by the proposed condensation.

We next estimate the structural-nonzero counts \(\operatorname{nnz}_{\mathrm{full}}\), \(\operatorname{nnz}_{\mathrm{stress}}\), and
\(\operatorname{nnz}_{\mathrm{cond}}\) of the three system matrices, with each count taken over the upper triangular part, including the
diagonal. The primal block contributes \(4N_m\) dense symmetric \(21\times21\) blocks, or \(4N_m\,21(21+1)/2\) entries. The equality matrix
contributes approximately \(16N_m\) local face blocks: \(12N_m\) from the two blocks associated with each of the 6 macro-internal faces
and approximately \(4N_m\) from macro-external face incidences. For the interpolation used here, direct evaluation of
\eqref{eq:local-face-equilibrium-block} gives exactly \(13\) structural nonzeros per row for a generic face away from the physical boundary.
Hence, \(\operatorname{nnz}(A_{F,i})=9\times13=117\), compared with \(189\) entries for a dense \(9\times21\) block. Storing the off-diagonal
equality block once therefore gives
\begin{equation}
	\operatorname{nnz}_{\mathrm{full}} \approx 4N_m\frac{21(21+1)}{2}+16\times117N_m = 2796N_m.
	\label{eq:full-structural-nnz}
\end{equation}
For the stress-condensed matrix, each contributing face generates a dense \(9\times9\) diagonal block. Since
\(n_{\rm eq}/9\approx8N_m\), these blocks contribute approximately \(8\times45N_m\) entries to the stored
triangle. Each of the four sub-elements of a macro element couples its four face-multiplier blocks and therefore
generates \(\binom{4}{2}=6\) dense off-diagonal \(9\times9\) blocks. No two distinct sub-elements share the same
pair of faces, so these off-diagonal blocks are not double-counted. Hence,
\begin{equation*}
	\operatorname{nnz}_{\mathrm{stress}}
	\approx 8\times45N_m+4\times6\times9^2N_m
	=2304N_m.
\end{equation*}
For the macro-condensed matrix \(\Sigma\) defined in \eqref{eq:compact-schur-complement}, the number of contributing
macro-external face blocks is
\begin{equation*}
	N_{\mathrm E}^{\mathrm f} := N_{\mathrm{sh}}^{\mathrm f}+N_{\Gamma_t}^{\mathrm f}
	= \frac{n_{\mathrm E}}{9} = \frac{n_{\rm eq}-54N_m}{9}.
\end{equation*}
Equation~\eqref{eq:equality-count-estimate} then gives \(N_{\mathrm E}^{\rm f}\approx2N_m\). Each contributing macro-external face
contributes \(9(9+1)/2=45\) entries to the triangularly stored diagonal, while the four outer faces of an interior macro element generate
\(\binom{4}{2}=6\) dense off-diagonal \(9\times9\) blocks. Because two distinct conforming macro elements cannot share the same pair of
complete faces, these off-diagonal blocks are not double-counted. Under the same refined-mesh approximation,
\begin{equation*}
	\operatorname{nnz}_{\mathrm{cond}} \approx 45N_{\mathrm E}^{\rm f}+6\times 9^2 N_m\approx 576N_m.
\end{equation*}
The corresponding leading-order ratios are
\begin{equation}
	\frac{\operatorname{nnz}_{\mathrm{full}}}
	     {\operatorname{nnz}_{\mathrm{stress}}}
	\approx\frac{2796}{2304}\approx1.21, \qquad
	\frac{\operatorname{nnz}_{\mathrm{stress}}}
	     {\operatorname{nnz}_{\mathrm{cond}}}
	\approx\frac{2304}{576}=4, \qquad
	\frac{\operatorname{nnz}_{\mathrm{full}}}
	     {\operatorname{nnz}_{\mathrm{cond}}}
	\approx\frac{2796}{576}\approx4.85.
	\label{eq:nnz-reduction-estimate}
\end{equation}
Thus, stress condensation reduces the global dimension by a factor of approximately \(2.17\) but the structural-nonzero count by only
\(1.21\), whereas eliminating the macro-internal multipliers provides a further factor of approximately \(4\) in both quantities.

The estimate in \eqref{eq:full-structural-nnz} excludes explicitly stored zero entries. If the \(n_{\rm eq}\) zero diagonal entries of the
multiplier block are stored, they must be added to \(\operatorname{nnz}_{\mathrm{full}}\); this convention does not alter the mathematical sparsity.

Finally, these estimates remain representative for the extended formulation defined by \eqref{eq:primal-macro-global-partition}
and \eqref{eq:equality-local-global-structure}, provided that the number of additional retained variables and the structural
nonzeros contributed by \(B_{\mathrm E}\) and \(\Omega_{\mathrm g}\) are small relative to the corresponding dimensions and
structural-nonzero counts of the base formulations. The accuracy of these estimates is assessed in Section~\ref{sec:numerical-results}.


\section{Nonsingularity of the Macro-Local Condensation Blocks}
\label{sec:local-condensation-block-nonsingularity}

The stress and macro-element condensations introduced in the preceding section impose related but distinct local regularity requirements.
Stress condensation requires the local stress block \(\Omega_e\) to be nonsingular for every macro element; since \(\Omega_e\succeq0\) by
\eqref{eq:macro-primal-block}, this is equivalent to \(\Omega_e\succ0\). Macro-element condensation instead requires the macro-local
saddle-point block \(K_e\) to be nonsingular, as formalized in Assumption~\ref{ass:local-block-regularity}.

Corollary~\ref{cor:internal-connection-rank} establishes that \(A_{\mathrm I,e}\) has full row rank, so the remaining requirement for
macro-element condensation concerns its interaction with the kernel of \(\Omega_e\). Theorem~\ref{thm:saddle-point-nonsingularity} applied
to \eqref{eq:local-saddle-point-block} therefore gives
\begin{equation}
	K_e\ \text{is nonsingular} \quad\Longleftrightarrow\quad \ker(\Omega_e)\cap\ker(A_{\mathrm I,e})=\{0\}.
	\label{eq:local-newton-kernel-condition}
\end{equation}
Thus, positive definiteness of \(\Omega_e\) guarantees the local validity of both condensations, whereas the weaker kernel-intersection
condition may suffice only for macro-element condensation. The following analysis establishes these properties for elastoplasticity and
Rankine limit analysis, identifies the hydrostatic null mode of the original von Mises formulation, and shows how trace-bound augmentation
restores the required local regularity.

For the elastoplastic formulation \eqref{eq:discrete-elastoplastic-q}, each \(P_i\succ0\) by \eqref{eq:assembled-compliance}, and hence
\(\mathcal P_e\succ0\) by \eqref{eq:macro-primal-components}. Equation~\eqref{eq:macro-primal-block} then gives \(\Omega_e\succ0\),
independently of the yield criterion, and \eqref{eq:local-newton-kernel-condition} establishes the nonsingularity of \(K_e\).

In limit analysis, \(\mathcal P_e=0\), so the local kernel condition depends on the macro-level yield operator and its
interaction with the internal traction equations. For any criterion-specific or augmented macro-level yield
operator \(\mathcal G_e\) considered below, let \(\Psi_e\succ0\) denote its corresponding scaling block.
The associated primal block is given by \(\Omega_e=\mathcal G_e^\top\Psi_e\mathcal G_e\).

For every \(q_e\in\mathbb{R}^{84}\),
\[ q_e^\top\Omega_eq_e = \left\|\Psi_e^{1/2}\mathcal G_eq_e\right\|_2^2. \]
Combining the resulting kernel identity with \eqref{eq:local-newton-kernel-condition} gives
\begin{equation}
	\begin{aligned}
		\ker(\Omega_e) &=\ker(\mathcal G_e),\\
		K_e\ \text{is nonsingular} &\quad\Longleftrightarrow\quad \ker(\mathcal G_e)\cap\ker(A_{\mathrm I,e})=\{0\}.
	\end{aligned}
	\label{eq:limit-analysis-macro-kernel-criterion}
\end{equation}
Thus, the IPM scaling does not affect the relevant nullspace. In particular, full column rank of \(\mathcal G_e\)
is sufficient, but not necessary, for the nonsingularity of \(K_e\).

\subsection{Rankine Limit Analysis}
\label{subsec:rankine-local-coercivity}

For the Rankine formulation \eqref{eq:rankine-vertex-data}, the yield operator of sub-element \(V_i\) is
\begin{equation}
	G_i^{\rm R} =
	\begin{bmatrix}
		-T_vS_{i,1}\\
		T_vS_{i,1}\\
		\vdots\\
		-T_vS_{i,4}\\
		T_vS_{i,4}
	\end{bmatrix} \in\mathbb{R}^{48\times21}.
	\label{eq:rankine-subelement-operator}
\end{equation}
For a macro element \(\mathcal M_e\), these operators assemble as
\begin{equation}
	\mathcal G_e^{\rm R} := \operatorname{blkdiag}_{r=1}^{4} \left(G_{i(e,r)}^{\rm R}\right).
	\label{eq:rankine-macro-operator}
\end{equation}

\begin{proposition}[Local nonsingularity for Rankine limit analysis]
	\label{prop:rankine-operator-full-rank}
	For every sub-element \(V_i\),
	\begin{equation}
		\operatorname{rank}(G_i^{\rm R})=21.
		\label{eq:rankine-subelement-full-rank}
	\end{equation}
	Consequently, \(\ker(\mathcal G_e^{\rm R})=\{0\}\). In Rankine limit analysis, \(\Omega_e\succ0\), and \(K_e\)
	is nonsingular for every macro element \(\mathcal M_e\).
\end{proposition}

\begin{proof}
	Let \(q_i\in\ker(G_i^{\rm R})\). Since \(T_v\) is nonsingular, \eqref{eq:rankine-subelement-operator} implies
	\(S_{i,a}q_i=0\) for \(a=1,\ldots,4\). The stress field \(S_iq_i\) is affine and vanishes at the four vertices
	of the nondegenerate tetrahedron. Affine unisolvency therefore gives \(S_iq_i\equiv0\) on \(V_i\). Because the
	columns of \(S_i\) form a basis of the local stress space, \(q_i=0\). Hence \(\ker(G_i^{\rm R})=\{0\}\), which
	proves \eqref{eq:rankine-subelement-full-rank}. The stated macro-level conclusions follow from
	\eqref{eq:rankine-macro-operator} and \eqref{eq:limit-analysis-macro-kernel-criterion}.
\end{proof}

\subsection{Von Mises Limit Analysis: Hydrostatic Null Mode and Trace-Bound Augmentation}
\label{subsec:von-mises-local-coercivity}

For the von Mises representation \eqref{eq:von-mises-vertex-data}, the sub-element yield operator is
\begin{equation}
	G_i^{\rm VM} =
	\begin{bmatrix}
		0_{1\times21}\\DS_{i,1}\\
		\vdots\\
		0_{1\times21}\\DS_{i,4}
	\end{bmatrix} \in\mathbb{R}^{24\times21}.
	\label{eq:von-mises-subelement-operator}
\end{equation}
Unlike the Rankine map, this operator is insensitive to hydrostatic stress.

\begin{proposition}[Kernel of the von Mises operator]
	\label{prop:von-mises-hydrostatic-kernel}
	For every sub-element \(V_i\),
	\begin{equation*}
		\ker(G_i^{\rm VM}) = \left\{ q_i\in\mathbb{R}^{21}: S_iq_i=pI_v \text{ for some }p\in\mathbb{R} \right\},
	\end{equation*}
	and hence
	\begin{equation*}
		\dim\ker(G_i^{\rm VM})=1, \qquad \operatorname{rank}(G_i^{\rm VM})=20.
	\end{equation*}
\end{proposition}

\begin{proof}
	Fix a sub-element \(V_i\), and let \(q_i\in\ker(G_i^{\rm VM})\). By \eqref{eq:von-mises-subelement-operator},
	\(DS_{i,a}q_i=0\) for \(a=1,\ldots,4\). The field \(DS_iq_i\) is affine and vanishes at the four vertices of
	the nondegenerate tetrahedron. Affine unisolvency therefore gives \(DS_iq_i\equiv0\). Since
	\(\ker(D)=\operatorname{span}\{I_v\}\), there is an affine scalar field \(p\) such that \(S_iq_i=p(x)I_v\).
	The associated tensor field is \(p(x){\mathbf I}\). Since the interpolation satisfies local equilibrium
	identically, \(0=\operatorname{div}(p{\mathbf I})=\nabla p\), so \(p\) is constant.
	
	Conversely, every constant hydrostatic field \(pI_v\) belongs to the complete \(21\)-mode affine equilibrium
	space and is mapped to zero by \(D\). Injectivity of the interpolation makes its coefficient representation
	unique. Therefore, \(\ker(G_i^{\rm VM})\) is one-dimensional, and rank--nullity gives \(\operatorname{rank}(G_i^{\rm VM})=20\).
\end{proof}

\begin{lemma}[Propagation of hydrostatic pressure]
	\label{lem:hydrostatic-pressure-propagation}
	Let a connected collection of sub-elements carry constant hydrostatic stress fields \(p_i{\mathbf I}\). If the
	internal traction-continuity equations are satisfied on every face shared by two sub-elements in the collection,
	then all \(p_i\) are equal.
\end{lemma}

\begin{proof}
	Let \(F\) be a face shared by sub-elements \(V_i\) and \(V_j\), and orient its unit normal \(\boldsymbol n_F\)
	outward from \(V_i\). The outward normal from \(V_j\) is then \(-\boldsymbol n_F\), so the discrete
	traction-continuity equation on \(F\) becomes
	\begin{equation*}
		\int_F N_F^\top (p_i-p_j)\boldsymbol n_F\,\mathrm d\Gamma=0.
	\end{equation*}
	The traction jump \((p_i-p_j)\boldsymbol n_F\) is constant and belongs to the linear face space represented by
	\(N_F\). Exact integration therefore implies \((p_i-p_j)\boldsymbol n_F=0\), and hence \(p_i=p_j\). Connectivity
	propagates this equality throughout the collection.
\end{proof}

For a macro element \(\mathcal M_e\), define
\begin{equation*}
	\mathcal G_e^{\rm VM} := \operatorname{blkdiag}_{r=1}^{4} \left(G_{i(e,r)}^{\rm VM}\right).
\end{equation*}
For pure von Mises limit analysis, its four sub-elements therefore contribute four independent hydrostatic null
modes to \(\Omega_e=(\mathcal G_e^{\rm VM})^\top\Psi_e\mathcal G_e^{\rm VM}\).
Lemma~\ref{lem:hydrostatic-pressure-propagation} shows that the internal traction equations remove their relative
variations but not their common value. More precisely, let \(q_{h,e}\ne0\) represent the same constant
hydrostatic stress \(p{\mathbf I}\), \(p\ne0\), on all four sub-elements of \(\mathcal M_e\). The outward normals
on each shared face are opposite, and consequently
\begin{equation*}
	\mathcal G_e^{\rm VM}q_{h,e}=0, \qquad A_{\mathrm I,e}q_{h,e}=0.
\end{equation*}
The local block is therefore singular:
\begin{equation*}
	K_e \begin{bmatrix} q_{h,e}\\0_{54} \end{bmatrix} =0,
\end{equation*}
even though \(A_{\mathrm I,e}\) has full row rank.

Although Proposition~\ref{prop:equilibrium-formulations-regularity} establishes the nonsingularity of the
assembled Newton system for von Mises limit analysis under Assumption~\ref{ass:reference-loading-regularity},
this global property does not imply local condensability.

\begin{remark}[Additional von Mises evaluations]
	Adding further von Mises constraints at the barycenter, or at any other point of a sub-element, cannot remove
	the hydrostatic mode because \(\operatorname{dev}(p{\mathbf I})=0\) everywhere. Such an augmentation therefore
	leaves the common hydrostatic null mode of \(K_e\) unchanged.
\end{remark}

In this work, this mode is instead eliminated by imposing a two-sided trace bound
\(|\!\operatorname{tr}(\sigma)|\le\bar\tau\) (\(\bar\tau>0\)). To maintain a uniform sub-element construction,
let \(x_{i,c}\) denote the barycenter of \(V_i\), set \(S_{i,c}:=S_i(x_{i,c})\), and evaluate this bound at the
barycenter of every sub-element:
\begin{equation*}
	\left| \operatorname{tr}\!\left( \operatorname{mat}_v(S_{i,c}q_i) \right) \right| \le \bar\tau.
\end{equation*}
Using \(I_v\) from \eqref{eq:svec-identity}, its conic representation is
\begin{equation}
	G_i^{\rm tr}q_i+s_i^{\rm tr}=h_i^{\rm tr}, \qquad s_i^{\rm tr}\in\mathcal Q_2,
	\label{eq:subelement-trace-bound-conic-form}
\end{equation}
where
\begin{equation}
	G_i^{\rm tr} :=
	\begin{bmatrix}
		0_{1\times21}\\
		-I_v^\top S_{i,c}
	\end{bmatrix}, \qquad h_i^{\rm tr}
	:=
	\begin{bmatrix}
		\bar\tau\\0
	\end{bmatrix}, \qquad \widehat G_i^{\rm VM}
	:=
	\begin{bmatrix}
		G_i^{\rm VM}\\
		G_i^{\rm tr}
	\end{bmatrix} \in\mathbb{R}^{26\times21}.
	\label{eq:subelement-trace-regularized-von-mises-operator}
\end{equation}
For macro element \(\mathcal M_e\), define
\begin{equation}
	\widehat{\mathcal G}_e^{\rm VM} := \operatorname{blkdiag}_{r=1}^{4}
	\left(\widehat G_{i(e,r)}^{\rm VM}\right) \in\mathbb{R}^{104\times84}.
	\label{eq:macro-subelement-trace-operator}
\end{equation}

\begin{proposition}[One trace bound per sub-element]
	\label{prop:subelement-trace-regularized-von-mises}
	For every sub-element \(V_i\), the augmented operator \(\widehat G_i^{\rm VM}\) in \eqref{eq:subelement-trace-regularized-von-mises-operator}
	satisfies \(\ker(\widehat G_i^{\rm VM})=\{0\}\). Consequently, in von Mises limit analysis, \(\Omega_e\succ0\), and \(K_e\) is
	nonsingular for every macro element.
\end{proposition}

\begin{proof}
	By Proposition~\ref{prop:von-mises-hydrostatic-kernel}, any \(q_i\in\ker(\widehat G_i^{\rm VM})\) represents a constant
	hydrostatic stress \(p\mathbf I\). The added trace row gives \(3p=0\), and injectivity of the interpolation yields \(q_i=0\).
	The macro-level conclusions follow from \eqref{eq:macro-subelement-trace-operator} and \eqref{eq:limit-analysis-macro-kernel-criterion}.
\end{proof}

Because the hydrostatic null mode is constant throughout \(V_i\), the same proof applies if the trace bound is
evaluated at any single point of the sub-element, including any of its vertices. The barycenter is selected
simply to avoid an arbitrary choice among the four vertices. Additionally, the trace construction augments \(G\)
only by additional rows and therefore cannot enlarge its kernel; the global nonsingularity established in
Proposition~\ref{prop:equilibrium-formulations-regularity} is thus preserved.

The four additional trace constraints per macro element introduce only cone slack and dual variables, which are
eliminated when forming the reduced Newton system \eqref{eq:ipm-reduced-system}. They modify only the local
contributions to \(\Omega_e\) and therefore do not change the dimensions of \(q_e\), \(y_{\mathrm I,e}\), or \(K_e\),
which remains in \(\mathbb{R}^{138\times138}\), nor the block sparsity pattern of the Newton and condensed systems.

The preceding construction is stronger than required for the condensation. Because the internal traction
equations couple the four sub-element hydrostatic modes, one trace bound on a single sub-element of each macro
element is sufficient to satisfy \eqref{eq:local-newton-kernel-condition}. To show this, select
\(r_\star\in\{1,\ldots,4\}\), set \(i_\star:=i(e,r_\star)\), and define
\begin{equation*}
	\mathcal S_{e,c} :=
	\begin{bmatrix}
		0 & \cdots & S_{i_\star,c} & \cdots & 0
	\end{bmatrix} \in\mathbb{R}^{6\times84},
\end{equation*}
so that \(\mathcal S_{e,c}q_e\) is the Voigt stress at the selected point. The trace bound has the same
\(\mathcal Q_2\) representation as \eqref{eq:subelement-trace-bound-conic-form}; its trace block and the
resulting augmented operator are
\begin{equation}
	\mathcal G_e^{\rm tr} :=
	\begin{bmatrix}
		0_{1\times84}\\
		-I_v^\top\mathcal S_{e,c}
	\end{bmatrix}, \qquad \widetilde{\mathcal G}_e^{\rm VM}
	:=
	\begin{bmatrix}
		\mathcal G_e^{\rm VM}\\
		\mathcal G_e^{\rm tr}
	\end{bmatrix} \in\mathbb{R}^{98\times84},
	\label{eq:macro-trace-regularized-von-mises-operator}
\end{equation}

\begin{proposition}[One trace bound per macro element]
	\label{prop:macro-trace-regularized-von-mises}
	For each macro element \(\mathcal M_e\), the operator \(\widetilde{\mathcal G}_e^{\rm VM}\) in
	\eqref{eq:macro-trace-regularized-von-mises-operator} satisfies
	\begin{equation*}
		\dim\ker(\widetilde{\mathcal G}_e^{\rm VM})=3, \qquad \ker(\widetilde{\mathcal G}_e^{\rm VM})
		\cap\ker(A_{\mathrm I,e})=\{0\}.
	\end{equation*}
	Consequently, in von Mises limit analysis, the associated \(\Omega_e\) is positive semidefinite and singular,
	while \(K_e\) is nonsingular.
\end{proposition}

\begin{proof}
	By Proposition~\ref{prop:von-mises-hydrostatic-kernel} and the block-diagonal assembly, \(\ker(\mathcal G_e^{\rm VM})\) is
	parametrized by four independent constant sub-element pressures. The added trace row imposes \(p_{r_\star}=0\), leaving a
	three-dimensional kernel. For \(q_e\in\ker(\widetilde{\mathcal G}_e^{\rm VM})\) satisfying \(A_{\mathrm I,e}q_e=0\),
	Lemma~\ref{lem:hydrostatic-pressure-propagation} makes all four pressures equal and hence zero; injectivity of the interpolation
	then gives \(q_e=0\). The remaining conclusions follow from \eqref{eq:limit-analysis-macro-kernel-criterion}.
\end{proof}

This distinction is particularly relevant for von Mises limit analysis. Stress condensation requires every \(21\times21\) sub-element
stress block to be nonsingular and hence one trace bound per sub-element. Macro-element condensation instead incorporates the internal
traction-continuity multipliers into \(K_e\); the internal traction equations then control the relative hydrostatic modes, so a single
trace bound per macro element makes \(K_e\) nonsingular even though \(\Omega_e\) remains singular.

The stronger sub-element construction in Proposition~\ref{prop:subelement-trace-regularized-von-mises} is nevertheless used in the
implementation because it yields \(\Omega_e\succ0\), thereby supporting both condensations and enabling the structure-aware block-Cholesky
factorization described in the next section. The trace bounds are introduced solely to remove the hydrostatic null mode and are not
intended to model pressure-sensitive yielding.

\begin{remark}[Existence of a nonrestrictive trace bound]
	\label{rem:trace-bound-optimum}
	The von Mises admissible set is unbounded in the hydrostatic direction, so any finite trace bound restricts the feasible set at the
	sub-element level. Nevertheless, suppose that the original discrete limit-analysis problem admits an optimal solution
	\((q^\star,\lambda^\star)\), and define
	\(T^\star := \max_i \left| \operatorname{tr}\!\left( \operatorname{mat}_v(S_{i,c} q_i^\star) \right) \right|.\)
	Because the discretization contains finitely many sub-elements, \(T^\star<\infty\). For any \(\bar\tau>T^\star\), the solution
	\(( q^\star,\lambda^\star)\) remains feasible for the trace-augmented problem. Since the augmented feasible set is contained
	in the original feasible set, the two problems have the same optimal limit multiplier, and \(( q^\star,\lambda^\star)\) remains optimal.
	If the original optimal stress field is unique, it is therefore unchanged.

	The argument does not determine a suitable value of \(\bar\tau\) a priori. The optimal load multiplier \(\lambda^\star\) is uniquely
	determined as the optimal value, whereas the associated optimal stress field need not be unique; a trace bound may therefore select
	among optimal stress fields without changing the multiplier. In the implementation, we use \(\bar\tau=10\sigma_y\), where \(\sigma_y\)
	is the uniaxial yield stress. Although this choice is not guaranteed a priori to be inactive, it did not alter the optimal load
	multiplier, to the reported numerical accuracy, in the example considered in Section~\ref{sec:numerical-results}.
\end{remark}

The kernel criterion in \eqref{eq:limit-analysis-macro-kernel-criterion} also provides a basis for extending the condensation to other
conic yield formulations, such as pressure-sensitive Drucker--Prager~\cite{makrodimopoulos2006lower} and frictional or truncated
Mohr--Coulomb models~\cite{krabbenhoft2008three}. For each representation, the local blocks must incorporate any auxiliary variables,
and the corresponding rank and kernel conditions must be verified. For the pressure-independent Tresca criterion, trace augmentation
can likewise remove the hydrostatic null mode.


\section{GPU Implementation of the IPM with Macro-Element Condensation}
\label{sec:gpu-condensed-solver}

This section presents the GPU implementation of the primal--dual IPM with macro-element condensation, written in CUDA Fortran using
libraries from the CUDA Toolkit~\cite{cuda128}. At each IPM iteration, independent cone-, sub-element-, and macro-local operations are
organized as dense batches, their contributions are assembled directly into the device-resident sparse matrix of the macro-condensed
system \eqref{eq:compact-condensed-system}, and the eliminated directions are recovered independently after each solve. Although the
assembly implementation described in this section focuses on the two base formulations, only a minor extension is needed to incorporate
the scalar rebar contributions of the reinforced-concrete formulation, as detailed in Subsection~\ref{subsec:numerical-reinforced-concrete}.
For the numerical comparisons, the stress and macro-element condensations are also implemented on the CPU. The CPU stress-condensation
implementation uses the same predictor--corrector framework but solves the stress-condensed system \eqref{eq:stress-only-condensed-system}
at each Newton step.

\subsection{Computational Workflow and Data Organization}
\label{subsec:gpu-workflow-data-organization}

The computational workflow separates the one-time setup phase from the numerical operations repeated at each IPM iteration, as summarized
in Algorithm~\ref{alg:ipm-predictor-corrector}. During setup, the sub-element data are grouped by cone family in a field-wise
structure-of-arrays layout and transferred to the device. The macro-local data are stored separately in contiguous batches organized by
macro element, enabling the batched local factorizations and applications of \(K_e^{-1}\) required for condensation,
right-hand-side formation, and recovery.
\begin{algorithm}[!t]
	\caption{Primal--dual interior-point method with macro-element condensation}
	\label{alg:ipm-predictor-corrector}
	\begin{algorithmic}[1]
		\Require Problem data and an initial iterate \(v=(x,s,y,z)\) with \(s,z\in\operatorname{int}\mathcal{K}\)
		\Statex \textbf{Setup:} Precompute the CSR structures and assembly maps,
		allocate the device workspaces, transfer the problem data to the GPU, and initialize the sparse solver
		\While{the convergence criteria are not satisfied}
		\State Compute the Nesterov--Todd scaling data \(W\) and \(\xi\) as described in~\cite{vandenberghe2010cvxopt}
		\label{line:ipm-scaling}
		\State {Form and factorize the macro-local blocks, then assemble and factorize the condensed matrix \(\Sigma\)}
		\label{line:ipm-factorization}
		\State Form the predictor right-hand side, solve \eqref{eq:compact-condensed-system}, and recover
		\(\Delta_{\rm a}\)
		\label{line:ipm-predictor}
		\State Compute \(\alpha_{\rm a}=\alpha_{\max}(\Delta_{\rm a})\) from \eqref{eq:ipm-admissible-step} and
		determine the centering parameter \(\sigma_{\rm cen}\)
		\label{line:ipm-affine-step}
		\State Form the corrector right-hand side, solve \eqref{eq:compact-condensed-system}, and recover \(\Delta_{\rm c}\)
		\label{line:ipm-corrector}
		\State Compute \(\alpha_{\rm c}=\alpha_{\max}(\Delta_{\rm c})\) from \eqref{eq:ipm-admissible-step}
		\label{line:ipm-corrector-step}
		\State Update \(v\gets v+0.99\,\alpha_{\rm c}\Delta_{\rm c}\)
		\label{line:ipm-update}
		\State Update the residuals \eqref{eq:ipm-residuals} and the stopping measures
		\label{line:ipm-residual-update}
		\EndWhile
	\end{algorithmic}
\end{algorithm}

\begin{figure}[!t]
	\centering
	\resizebox{\textwidth}{!}{\definecolor{layoutInk}{HTML}{111827}
\definecolor{layoutText}{HTML}{374151}
\definecolor{layoutMuted}{HTML}{6B7280}
\definecolor{layoutPanel}{HTML}{F8FAFC}
\definecolor{layoutBlue}{HTML}{2563EB}
\definecolor{layoutBluePale}{HTML}{EFF6FF}
\definecolor{layoutOrange}{HTML}{D97706}
\definecolor{layoutOrangePale}{HTML}{FFF7ED}

\begin{tikzpicture}[
    x=1cm,
    y=1cm,
    font=\small,
    panel/.style={draw=layoutInk, line width=0.8pt, rounded corners=2pt,
                  fill=layoutPanel},
    record/.style={draw=layoutMuted, line width=0.65pt, rounded corners=1.5pt,
                   fill=white, minimum height=0.55cm},
    socgroup/.style={draw=layoutBlue, line width=0.8pt, rounded corners=2pt,
                     fill=layoutBluePale},
    psdgroup/.style={draw=layoutOrange, line width=0.8pt, rounded corners=2pt,
                     fill=layoutOrangePale},
    socfield/.style={draw=layoutBlue, line width=0.55pt, rounded corners=1.2pt,
                     fill=white, minimum height=0.42cm},
    psdfield/.style={draw=layoutOrange, line width=0.55pt, rounded corners=1.2pt,
                     fill=white, minimum height=0.42cm},
    batchcell/.style={draw=layoutBlue, line width=0.65pt, fill=layoutBluePale,
                      minimum height=0.62cm},
    matrixcell/.style={draw=layoutInk, line width=0.65pt, rounded corners=1.2pt,
                       fill=white, minimum height=0.52cm},
    csrcell/.style={draw=layoutMuted, line width=0.65pt, fill=white,
                    minimum height=0.55cm}
]

\draw[panel] (0,2.35) rectangle (7.45,8.55);
\draw[panel] (7.75,2.35) rectangle (15.80,8.55);
\draw[panel] (0,0) rectangle (15.80,2.10);

\node[font=\bfseries\large, text=layoutInk] at (3.725,8.18)
    {(a) Conventional CPU layout};
\node[text=layoutMuted] at (3.725,7.78)
    {array of structures (AoS)};

\node[font=\bfseries\large, text=layoutInk] at (11.775,8.18)
    {(b) GPU-oriented layout};
\node[text=layoutMuted] at (11.775,7.78)
    {family-wise SoA and contiguous batches};

\node[anchor=west, font=\bfseries, text=layoutText] at (0.35,7.30)
    {Cone records};

\node[anchor=east, text=layoutText] at (1.20,6.78) {cone \(1\)};
\node[record, anchor=west, minimum width=5.72cm] at (1.33,6.78)
    {\(G_1\;\vert\;h_1\;\vert\;\mathrm{idx}_1\;\vert\;s_1\;\vert\;z_1\;\vert\;W_1
      \;\vert\;\cdots\)};

\node[anchor=east, text=layoutText] at (1.20,6.08) {cone \(2\)};
\node[record, anchor=west, minimum width=5.72cm] at (1.33,6.08)
    {\(G_2\;\vert\;h_2\;\vert\;\mathrm{idx}_2\;\vert\;s_2\;\vert\;z_2\;\vert\;W_2
      \;\vert\;\cdots\)};

\node[text=layoutMuted] at (4.04,5.45) {\(\vdots\)};

\node[anchor=east, text=layoutText] at (1.20,4.84)
    {cone \(j\)};
\node[record, anchor=west, minimum width=5.72cm] at (1.33,4.84)
    {\(G_j\;\vert\;h_j\;\vert\;\mathrm{idx}_j\;\vert\;s_j\;\vert\;z_j\;\vert\;W_j
      \;\vert\;\cdots\)};

\node[anchor=west, font=\bfseries, text=layoutText] at (0.35,4.20)
    {Macro-element records};
\node[record, minimum width=1.30cm] at (1.30,3.47) {\(K_1\)};
\node[record, minimum width=1.30cm] at (3.02,3.47) {\(K_2\)};
\node[text=layoutMuted] at (4.25,3.47) {\(\cdots\)};
\node[record, minimum width=1.72cm] at (5.68,3.47) {\(K_{N_m}\)};

\node[align=center, text=layoutMuted] at (3.725,2.76)
    {All fields associated with one object\\are grouped in the same record.};

\node[anchor=west, font=\bfseries, text=layoutText] at (8.10,7.30)
    {Cone-family arrays};

\draw[socgroup] (8.10,4.42) rectangle (11.62,7.03);
\node[font=\bfseries, text=layoutBlue] at (9.86,6.77)
    {SOC family};
\node[socfield, minimum width=1.36cm] at (8.98,6.22) {\(G_{\rm soc}\)};
\node[socfield, minimum width=1.36cm] at (10.74,6.22) {\(h_{\rm soc}\)};
\node[socfield, minimum width=1.36cm] at (8.98,5.66) {\(\mathrm{idx}_{\rm soc}\)};
\node[socfield, minimum width=1.36cm] at (10.74,5.66) {\(s_{\rm soc}\)};
\node[socfield, minimum width=1.36cm] at (8.98,5.10) {\(z_{\rm soc}\)};
\node[socfield, minimum width=1.36cm] at (10.74,5.10) {\(W_{\rm soc}\)};
\node[text=layoutMuted] at (9.86,4.74) {\(\vdots\)};

\draw[psdgroup] (11.82,4.42) rectangle (15.45,7.03);
\node[font=\bfseries, text=layoutOrange] at (13.635,6.77)
    {PSD family};
\node[psdfield, minimum width=1.40cm] at (12.72,6.22) {\(G_{\rm psd}\)};
\node[psdfield, minimum width=1.40cm] at (14.54,6.22) {\(h_{\rm psd}\)};
\node[psdfield, minimum width=1.40cm] at (12.72,5.66) {\(\mathrm{idx}_{\rm psd}\)};
\node[psdfield, minimum width=1.40cm] at (14.54,5.66) {\(s_{\rm psd}\)};
\node[psdfield, minimum width=1.40cm] at (12.72,5.10) {\(z_{\rm psd}\)};
\node[psdfield, minimum width=1.40cm] at (14.54,5.10) {\(W_{\rm psd}\)};
\node[text=layoutMuted] at (13.635,4.74) {\(\vdots\)};

\node[anchor=west, font=\bfseries, text=layoutText] at (8.10,4.08)
    {Macro-local batch \(\{K_e\}_{e=1}^{N_m}\)};
\node[batchcell, minimum width=1.72cm] at (8.96,3.50) {\(K_1\)};
\node[batchcell, minimum width=1.72cm] at (10.68,3.50) {\(K_2\)};
\node[batchcell, minimum width=1.72cm] at (12.40,3.50) {\(\cdots\)};
\node[batchcell, minimum width=1.72cm] at (14.12,3.50) {\(K_{N_m}\)};

\node[align=center, text=layoutMuted] at (11.775,2.76)
    {Separate fields and equal-size matrices\\are batched contiguously.};

\node[font=\bfseries, text=layoutInk] at (7.90,1.76)
    {Common global sparse representation (CPU and GPU)};

\node[matrixcell, minimum width=0.90cm] at (0.78,1.02) {\(A\)};
\node[matrixcell, minimum width=0.90cm] at (1.78,1.02) {\(P\)};
\node[matrixcell, minimum width=1.82cm] at (3.24,1.02)
    {\(G^\top\Theta G\)};
\node[matrixcell, minimum width=0.90cm] at (4.70,1.02) {\(\Omega\)};
\node[matrixcell, minimum width=0.90cm] at (5.70,1.02) {\(\Sigma\)};

\node[text=layoutMuted] at (7.10,1.02) {\(\longrightarrow\)};
\node[csrcell, minimum width=2.40cm] at (9.225,1.02)
    {row pointers};
\node[csrcell, minimum width=2.40cm] at (11.625,1.02)
    {column indices};
\node[csrcell, minimum width=1.72cm] at (13.685,1.02)
    {values};
\node[text=layoutMuted] at (7.90,0.37)
    {Fixed row and column structures; numerical values updated as required.};

\end{tikzpicture}}
	\caption{Comparison of the conventional CPU and GPU-oriented data layouts. In the CPU array-of-structures
		representation, the data belonging to each cone or macro element are grouped in an individual record. On the
		GPU, cone quantities are separated by field and cone family, and the macro-local data associated
		with \(K_e\) form contiguous batches. The global matrices \(A\), \(P\), \(G^\top\Theta G\), \(\Omega\), and \(\Sigma\) use CSR
		storage in both implementations.}
	\label{fig:gpu-data-layout}
\end{figure}
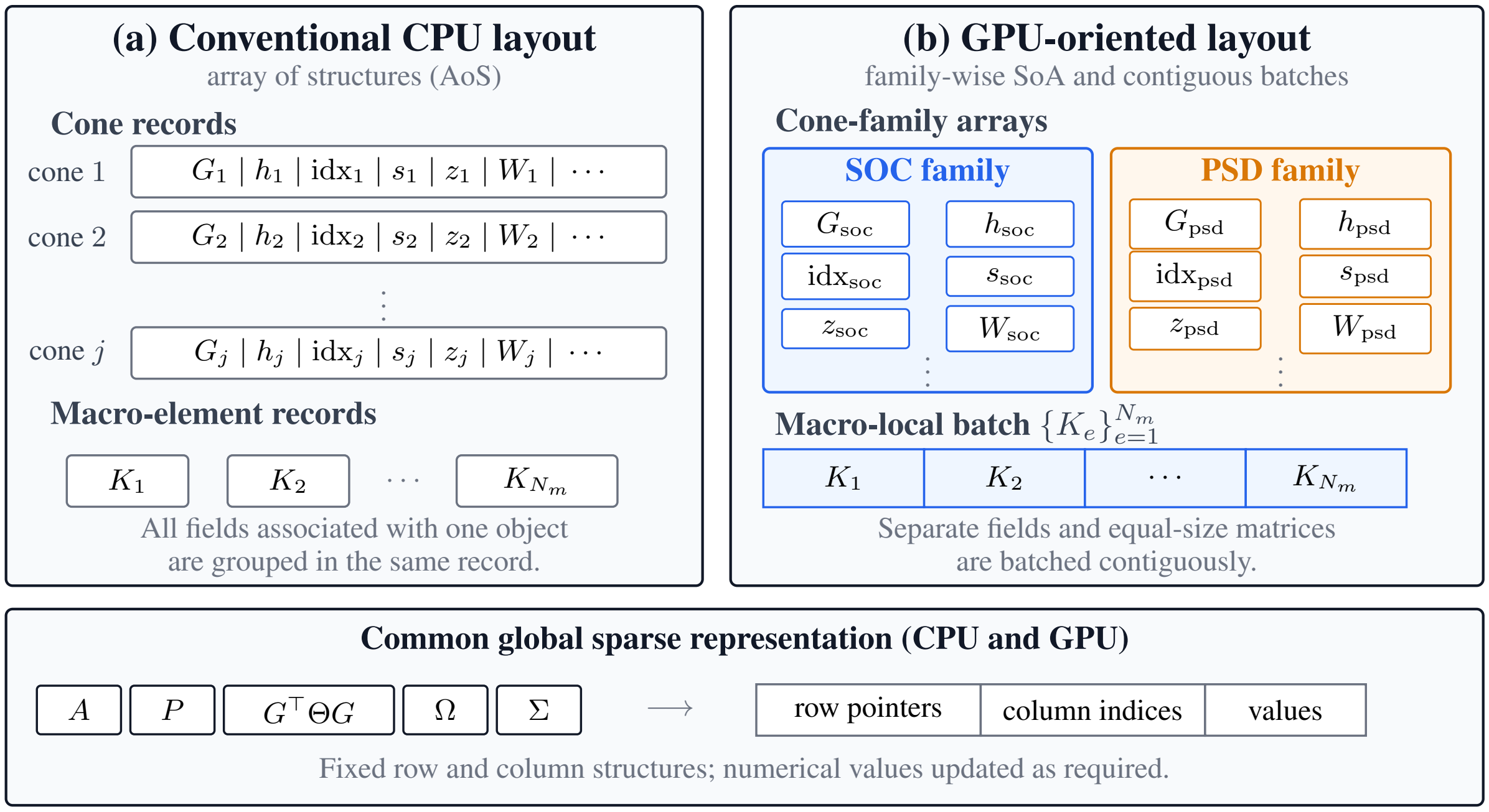

The sparse global representation is common to the CPU and GPU implementations. The matrices \(A\), \(P\),
\(G^\top\Theta G\), \(\Omega\), and \(\Sigma\) are stored in CSR format. Their sparsity patterns are constructed
during setup and remain fixed throughout the IPM iterations. The matrices \(A\) and \(P\) are constant, whereas
the numerical values of \(G^\top\Theta G\), \(\Omega\), and \(\Sigma\) are updated at every iteration. On the
GPU, cuSPARSE is used for CSR matrix additions, including \(\Omega=P+G^\top\Theta G\), and for all sparse
matrix--vector products involving the global CSR matrices, including \(Px\), \(A^\top y\), and \(Ax\) during
residual evaluation. Precomputed scatter maps place each numerical contribution directly into its prescribed
sparse-matrix entries without repeating structural searches inside the IPM loop. Figure~\ref{fig:gpu-data-layout}
compares the CPU and GPU local-data organizations and shows their common global sparse representation.

Let \(N_f\) denote the number of cones in family \(f\), and let \(\bar d_f\) and \(\bar n_f\) denote the
corresponding maximum cone and local-variable dimensions. The family-wise yield operators are stored in the
three-dimensional array \(G_f\in\mathbb{R}^{\bar d_f\times\bar n_f\times N_f}\), while cone vectors such as
\(h_f\), \(s_f\), and \(z_f\) are stored in two-dimensional arrays in \(\mathbb{R}^{\bar d_f\times N_f}\). The
final index enumerates the cones; hence, \(G_f(:,:,j)\) and the vector slices \(h_f(:,j)\), \(s_f(:,j)\), and
\(z_f(:,j)\) associated with cone \(j\) are contiguous, and successive cone slices are stored consecutively.
The arrays \(W_{\rm soc}\), \(\xi_{\rm soc}\), \(W_{\rm psd}\), and \(\xi_{\rm psd}\) analogously collect the
family-specific scaling data.

This layout provides contiguous cone slices and regular strided-batched access. In particular, it enables the matrix--vector products
with \(G\) and \(G^\top\), required for residual evaluation, condensed right-hand-side formation, and cone-direction recovery,
to be computed efficiently using the strided-batched DGEMV routine provided by cuBLAS.

The family-wise layout also supports kernels specialized to the algebra of each cone. SOC operations use
closed-form vector operations, including dot products, scaled sums, square roots, and cone-local norm reductions.
The \(3\times3\) PSD operations require several small dense decompositions: Cholesky factorizations and SVDs for
the Nesterov--Todd scaling, together with symmetric eigenvalue computations for the admissible-step calculation.
For the batches of \(3\times3\) matrices considered here, general-purpose batched routines provided by the CUDA
libraries have high overhead relative to the arithmetic itself, and specialized kernels are therefore substantially
faster~\cite{haidar2017guide,abdelfattah2026efficient}. Accordingly, our implementation uses dedicated kernels
for the \(3\times3\) Cholesky, SVD, and eigenvalue computations.

Within an IPM iteration, \(\Theta\), the macro-local blocks \(K_e\), and \(\Sigma\) are unchanged between the predictor and corrector
solves. Consequently, the macro-local factors and the factorization of \(\Sigma\) are computed once, while right-hand-side formation,
backsolves, and recovery are carried out for each direction. In the main execution path, all data remain GPU-resident throughout the
IPM loop, except for a few scalar quantities, such as \(\alpha\), \(\sigma_{\rm cen}\), \(\mu\), and the residual norms, that are
exchanged with the CPU. Unless stated otherwise, all numerical operations use double precision; the limited use of single precision
is introduced later. Device workspaces are allocated during setup and reused throughout the solve. Dense temporaries are processed in
chunks of macro elements when a single batch would exceed the available workspace; this bounds temporary storage
without changing the algebraic condensation.

Beyond these matrix operations, the reductions required for iteration-wide norms and inner products in both the
SOC and PSD families are performed using the cuBLAS routines \texttt{cublasDnrm2} and \texttt{cublasDdot}. For PSD
quantities, the weights induced by the symmetric vectorization are included. By contrast, the step-length
computations in Lines~\ref{line:ipm-affine-step} and~\ref{line:ipm-corrector-step} use a custom hierarchical
reduction. Denote by \(\overline{\Delta s}_c\) and \(\overline{\Delta z}_c\) the directions obtained through
the standard Nesterov--Todd/Jordan-algebra normalization associated with \(\xi_c\), which maps the current scaled
cone point to the Jordan identity~\cite{nesterov1998primal,vandenberghe2010cvxopt}. Each SOC or PSD cone \(c\)
produces the local violation measure
\[
\rho_c :=
\max\left\{ 0,\, -\lambda_{\min}(\overline{\Delta s}_c),\, -\lambda_{\min}(\overline{\Delta z}_c) \right\}.
\]
Here, \(\lambda_{\min}\) denotes the smallest Jordan eigenvalue in the corresponding cone. For an SOC element
\(u=(u_0,u_v)\), it is \(u_0-\lVert u_v\rVert_2\); evaluating the norm requires a local reduction over the
components of \(u_v\). For a PSD element, \(\lambda_{\min}\) is the ordinary smallest matrix eigenvalue,
evaluated analytically for the \(3\times3\) cones considered here. This cone-local computation precedes the
hierarchical reduction of \(\rho_c\) over all cones summarized in Algorithm~\ref{alg:gpu-step-size-reduction}.
\begin{algorithm}[!htbp]
	\caption{Warp-reduced evaluation of the maximum admissible step}
	\label{alg:gpu-step-size-reduction}
	\begin{algorithmic}[1]
		\Require Cone scaling data and directions \(\Delta s\) and \(\Delta z\)
		\State Initialize the device scalar \(\rho_{\max}\gets0\)
		\ForAll{cones \(c\) in the SOC and PSD families in parallel}
		\State Compute the cone-local measure \(\rho_c\)
		\State Compute \(\rho_{\rm w}\gets\operatorname{warpMax}(\rho_c)\) using register-shuffle operations
		\If{the thread is the warp leader}
		\State \(\operatorname{atomicMax}(\rho_{\max},\rho_{\rm w})\)
		\EndIf
		\EndFor
		\State Set \(\alpha_{\max}\gets\min\{1,\rho_{\max}^{-1}\}\) if
		\(\rho_{\max}>0\), and \(\alpha_{\max}\gets1\) otherwise
	\end{algorithmic}
\end{algorithm}

\subsection{Parallel Formation of the Macro-Local Blocks}
\label{subsec:gpu-macro-local-block-formation}

Each IPM iteration begins with the Nesterov--Todd scaling computation listed in Line~\ref{line:ipm-scaling} of
Algorithm~\ref{alg:ipm-predictor-corrector}. Since the Nesterov--Todd operator is block diagonal over the
elementary cones, its blocks can be computed independently. Separate device kernels process the SOC and PSD
families and write the resulting compact scaling data directly into the family-wise arrays. These data are then
used to form the cone-wise contributions to \(G^\top\Theta G\), from which \(\Omega\) and the
iteration-dependent primal blocks \(\Omega_e\) of \(K_e\) are assembled as part of
Line~\ref{line:ipm-factorization}.

The matrix \(\Theta_c\) associated with an elementary cone \(c\) is not formed explicitly. Instead, the
corresponding cone operator is transformed as
\begin{equation*}
	H_c:=W_c^{-\top}G_c, \qquad G_c^\top\Theta_cG_c =G_c^\top W_c^{-1}W_c^{-\top}G_c =H_c^\top H_c.
\end{equation*}
For each cone family, a dedicated kernel applies the compact representation of \(W_c^{-\top}\) column by column
to \(G_c\) and immediately forms the upper triangle of \(H_c^\top H_c\), without storing \(H_c\) in global memory.

The vertexwise Gram contributions belonging to a sub-element \(V_i\) act on the same coefficient vector \(q_i\)
and are therefore accumulated locally. By contrast, the variables \(q_i\) of distinct sub-elements are disjoint,
so the resulting sub-element blocks occupy nonoverlapping positions in the CSR value array of \(G^\top\Theta G\).
Their destination ranges are precomputed during setup. A GPU work group can consequently process all cones
associated with one sub-element and write the completed block directly to its prescribed CSR positions, without
structural searches, inter-sub-element atomics, or synchronization between sub-elements. This gives parallelism
both across sub-elements and across the entries of each local Gram matrix. Once these values have been assembled,
the primal matrix \(\Omega=P+G^\top\Theta G\) is obtained by sparse CSR addition on the device. The symbolic
structure and workspace of this addition are prepared once during setup.

Precomputed macro-element maps then gather the relevant entries of \(\Omega\) into contiguous batches of macro-local primal data.
Together with the constant internal-equilibrium blocks \(A_{\mathrm I,e}\), prepared during setup, these iteration-dependent
data define \(K_e\) in \eqref{eq:local-saddle-point-block}. The batched factorization and condensation described next exploit
the local block and coupling structure to reduce dense storage.

\subsection{Batched Condensation, Assembly, and Local Recovery}
\label{subsec:gpu-batched-condensation-assembly}

The macro-external coupling matrix \(C_e\) defined in \eqref{eq:macro-external-coupling-block} has nonzero columns only for the
multiplier components associated with the contributing faces of \(\mathcal M_e\), with \(9\) components per face.
Depending on the boundary topology, some macro elements have only \(27\) active face components. For the batched formulation, let
\(\bar C_e\in\mathbb R^{138\times36}\) denote the local coupling matrix, with zero columns corresponding to inactive face components.
The associated connectivity is represented by the Boolean localization matrix
\begin{equation}
	L_e\in\{0,1\}^{36\times n_{\mathrm E}}, \qquad C_e=\bar C_eL_e,
	\label{eq:gpu-boolean-localization-matrix}
\end{equation}
where each active row of \(L_e\) gathers one macro-external multiplier component associated with
\(\mathcal M_e\), while rows introduced only for padding are zero.

With this localization, the assembly of the macro-external matrix
\(\Sigma_{\mathrm E}\) in
\eqref{eq:macro-external-schur-complement} takes the form
\begin{equation*}
	X_e:=K_e^{-1}\bar C_e, \qquad
	\Phi_e:=\bar C_e^\top X_e, \qquad
	\Sigma_{\mathrm E}=\sum_{e=1}^{N_m}L_e^\top\Phi_eL_e,
\end{equation*}
where \(X_e\in\mathbb R^{138\times36}\) and
\(\Phi_e\in\mathbb R^{36\times36}\).

The symmetric-indefinite matrices \(K_e\in\mathbb R^{138\times138}\) could be processed as a dense batch using the general LU routines
\texttt{cublasDgetrfBatched} and \texttt{cublasDgetrsBatched}, since the CUDA libraries do not provide the required batched symmetric-indefinite
factorization. In practice, this route becomes unreliable in the final IPM iterations, when the sub-element primal blocks in \(\Omega_e\) become severely
ill-conditioned. The implementation therefore uses only the structure-aware positive-definite factorization described below.

As established in Section~\ref{sec:local-condensation-block-nonsingularity}, \(\Omega_e\succ0\) for elastoplasticity and, in limit analysis,
for Rankine by Proposition~\ref{prop:rankine-operator-full-rank} and for the von Mises formulation equipped with one trace bound per
sub-element by Proposition~\ref{prop:subelement-trace-regularized-von-mises}. Writing \(\Omega_e=\operatorname{blkdiag}_{r=1}^4(\Omega_{e,r})\),
the Cholesky path applies batched factorizations to the four \(21\times21\) sub-element primal blocks \(\Omega_{e,r}\) and to the
positive-definite matrix \(S_{\mathrm I,e}:=A_{\mathrm I,e}\Omega_e^{-1}A_{\mathrm I,e}^{\top} \in\mathbb R^{54\times54}\); the
corresponding internal-multiplier Schur complement is \(-S_{\mathrm I,e}\). These batched Cholesky factorizations use \texttt{cusolverDnDpotrfBatched},
while the corresponding triangular solves and matrix products use cuBLAS routines. The resulting factors are reused for condensed-matrix formation,
both condensed right-hand sides, and local recovery.

As the IPM approaches convergence, some sub-element primal blocks \(\Omega_{e,r}\) become severely ill-conditioned. Although these
blocks are positive definite in exact arithmetic, their finite-precision representations may exhibit very small or slightly
negative eigenvalues, making the Cholesky path unreliable. Because a fixed diagonal regularization \(\Omega_{e,r}+\delta I\) was not robust
across the changing block scales, we instead use relative minimum-eigenvalue reconditioning \cite{tabeart2020improving}. When the Cholesky
factorization \(\Omega_{e,r}=R_{e,r}^{\top}R_{e,r}\) succeeds, its reliability is assessed using the inexpensive reciprocal-conditioning proxy
\[\widehat{\kappa}_{e,r}^{-1}:=\left(\frac{\min_j |(R_{e,r})_{jj}|}{\max_j |(R_{e,r})_{jj}|}\right)^2.\]
A block is selected for spectral treatment if its Cholesky factorization fails or if \(\widehat{\kappa}_{e,r}^{-1}<\varepsilon_{\mathrm{chol}}\),
where \(\varepsilon_{\mathrm{chol}}=10^{-13}\).

For both CPU and GPU implementations, a selected block is decomposed as \(\Omega_{e,r}=Q_{e,r}\operatorname{diag}(\lambda_j)Q_{e,r}^{\top},\)
and its eigenvalues are replaced by
\(\widetilde{\lambda}_j:=\max\!\left\{\lambda_j,\varepsilon_{\mathrm{eig}}\lambda_{\max} \right\},\) with \(\varepsilon_{\mathrm{eig}}=10^{-13}\)
and \(\lambda_{\max}:=\max_j\lambda_j\). This relative floor replaces roundoff-level nonpositive eigenvalues and bounds the condition number
of the modified block by \(\varepsilon_{\mathrm{eig}}^{-1}\).

In the CPU implementations, \texttt{dpotrf} is called separately for every sub-element block at each IPM iteration, so spectral
reconditioning is applied only to the blocks that trigger the criterion, while all remaining blocks retain their Cholesky factors.
This selective procedure is used in both the stress and macro-element condensations. On the GPU, \texttt{cusolverDnDpotrfBatched}
processes the complete collection of sub-element blocks; if any block triggers the criterion, the entire batch is switched to the
spectral path for the remainder of the IPM solve, thereby preserving uniform batched execution.

The GPU eigendecompositions are computed using \texttt{cusolverDnDsyevjBatched}. From the floored eigenvalues, the inverse-square-root factors
\[F_{e,r}=Q_{e,r}\operatorname{diag}\left(\widetilde{\lambda}_j^{-1/2}\right)\]
are used in the local applications and to form \(S_{\mathrm I,e}\) as a Gram matrix. The resulting \(S_{\mathrm I,e}\) blocks are
factorized by batched Cholesky without additional spectral modification, and all spectral operations and
subsequent local applications remain GPU resident.

In addition to improving numerical robustness, the structure-aware path reduces the storage required by the macro-local operations.
The nonzero stress-coupling portion of \(\bar C_e\) is represented by four \(21\times9\) face blocks aligned with the four blocks of
\(\Omega_e\). Compared with the padded \(138\times36\) dense representation, this reduces the coupling-block storage by approximately
\(85\%\). Likewise, the four \(21\times21\) factors, the \(54\times54\) factor, and the stored \(84\times54\) block \(A_{\mathrm I,e}^{\top}\)
together require approximately \(52\%\) less storage than a full \(138\times138\) dense factor.

The present implementation retains the macro-local factors and \(A_{\mathrm I,e}^{\top}\) blocks on the device to avoid repeated
extraction and factorization within each IPM iteration. An alternative would gather the required blocks from the global CSR matrices into
reusable chunk workspaces and recompute the local factors as needed. With a fixed chunk size, the workspace required for these local data
would be independent of the total number of macro elements, at the cost of repeated block extraction and factorization. This tradeoff
could accommodate larger discretizations when retained macro-local data account for a substantial share of device-memory usage.

Using the retained local factors and compact coupling blocks, \(X_e\) and \(\Phi_e\) are evaluated in reusable chunk workspaces.
Batched matrix products form each complete \(\Phi_e\), whose upper-triangular entries are scattered directly into the macro-external part
of the CSR value array of \(\Sigma\) using precomputed maps \(\pi_e(i,j)\). For pairs associated with active rows of \(L_e\), these maps
give the corresponding global upper-triangular CSR positions; pairs involving padding are assigned valid dummy positions and contribute
zero. Thus, \(\pi_e\) implements the action of \(L_e^\top(\cdot)L_e\) without explicitly storing \(L_e\). Since neighboring macro elements
may contribute to the same entry of \(\Sigma\), the scatter kernel uses hardware-supported atomic addition to accumulate these
contributions without explicit locking. Algorithm~\ref{alg:gpu-condensed-assembly} summarizes this batched assembly using the
structure-aware local factors.

For both the predictor and corrector steps, the macro-external right-hand side \(\widehat\beta_{\mathrm E}\), defined in
\eqref{eq:macro-external-schur-complement}, is assembled using the same batched local solves and sparse accumulation:
\begin{equation*}
	t_e:=K_e^{-1}b_e, \qquad \gamma_e:=\bar C_e^\top t_e,
	\qquad \widehat\beta_{\mathrm E}
	=-b_{y,\mathrm E}+\sum_{e=1}^{N_m}L_e^\top\gamma_e.
\end{equation*}
After each condensed solve, the eliminated directions are recovered by \eqref{eq:local-solution}. For each macro element, the local
right-hand side \(b_e-\bar C_eL_e\Delta y_{\mathrm E}\) is formed independently, and \(w_e\) is computed using the same batched
local factors used for right-hand-side assembly. Combining the recovered vectors \(w_e\) with \(\Delta\chi\) then gives the full Newton
direction \((\Delta x,\Delta y)\).
\begin{algorithm}[!t]
	\caption{Batched macro-local condensation and parallel CSR assembly of \(\Sigma\)}
	\label{alg:gpu-condensed-assembly}
	\begin{algorithmic}[1]
		{\Require Local Cholesky or spectral factor data for
		\(\Omega_{e,r}\), Cholesky factors of \(S_{\mathrm I,e}\),
		blocks \(A_{\mathrm I,e}^{\top}\), stored face blocks representing
		\(\bar C_e\), assembly maps \(\pi_e\), and \(\ell_{\mathrm E}\)
		for limit analysis}
		\State {Initialize \(\Sigma_{\rm val}\) to zero and, for limit analysis, insert \(\ell_{\mathrm E}\) at the
		precomputed load-factor border positions}
		\State {Compute \(X_e\gets K_e^{-1}\bar C_e\) for all \(e\) using batched local solves}
		\State Compute \(\Phi_e\gets \bar C_e^\top X_e\) for all \(e\) using batched matrix products
		\ForAll{macro elements \(e\) and upper-triangular entries \((\Phi_e)_{ij}\), in parallel}
		\State \(k\gets\pi_e(i,j)\)
		\State \(\operatorname{atomicAdd} \bigl({\Sigma_{\rm val}}[k],(\Phi_e)_{ij}\bigr)\)
		\EndFor
	\end{algorithmic}
\end{algorithm}

\subsection{Solution of the Condensed System}
\label{subsec:gpu-condensed-system-solution}

At each IPM iteration, the condensed system \eqref{eq:compact-condensed-system} is solved for the predictor and corrector directions:
\begin{equation*}
	\Sigma\Delta\chi^{\nu}=\widehat\beta^{\nu}, \qquad     \nu\in\{\mathrm a,\mathrm c\}.
\end{equation*}
The CSR pattern of \(\Sigma\) is fixed, while its numerical values are updated once per IPM iteration and shared
by the two solves. Our implementation uses the NVIDIA cuDSS sparse direct solver~\cite{cudss080} as the primary
option for all formulations. For the positive-definite elastoplastic system, a preconditioned conjugate-gradient
(PCG) method~\cite{Saad2003} is also available as a memory-oriented alternative.

\subsubsection{Direct Solution with cuDSS}
\label{subsubsec:gpu-direct-cudss}

cuDSS operates directly on the device-resident CSR and dense arrays. Its symbolic analysis, including ordering
and allocation of the factorization data, is performed once during setup and subsequently reused. At every IPM
iteration, cuDSS performs one numerical factorization of the updated matrix \(\Sigma\), followed by separate
backsolves for \(\widehat\beta^{\mathrm a}\) and \(\widehat\beta^{\mathrm c}\), without host--device transfers.

For the standard elastoplastic system \eqref{eq:elastoplastic-condensed-system}, the coefficient matrix \(\Sigma=\Sigma_{\mathrm E}\) is
positive definite; cuDSS is therefore configured for a symmetric positive-definite factorization. For limit analysis, the
symmetric-indefinite system \eqref{eq:limit-analysis-condensed-system} is solved directly with cuDSS. This avoids the additional
solve with \(\ell_{\mathrm E}\) required by the positive-definite alternative described in Remark~\ref{rem:positive-definite-limit-analysis-solve}.
Avoiding this solve is particularly relevant on GPUs, where dependencies in sparse triangular solves can limit parallelism.

The direct path is adopted as the primary solution strategy to avoid the potentially high Krylov iteration counts
associated with increasingly ill-conditioned IPM systems. Its main limitation is the device memory required by the
sparse factors, which depends on the fill generated by the ordering. The reductions in dimension, structural-nonzero
count, and multiplier-block size established in Subsection~\ref{subsec:sparsity-size-parallelism} therefore extend
the range of problem sizes that can be handled by GPU-resident sparse factorization relative to the original Newton system.

\subsubsection{Memory-Oriented PCG Alternative}
\label{subsubsec:gpu-pcg-alternative}

The elastoplastic condensed system can alternatively be solved using PCG. For this option, the numerical factorization of \(\Sigma\) in
Line~\ref{line:ipm-factorization} of Algorithm~\ref{alg:ipm-predictor-corrector} is replaced by preconditioner construction. The resulting
preconditioner is reused for the predictor and corrector solves. This option is not used for limit analysis here, since the positive-definite
reduction in Remark~\ref{rem:positive-definite-limit-analysis-solve} would require an additional solve with \(\ell_{\mathrm E}\) per IPM iteration.

The structure-exploiting preconditioner proposed here is a macro-element additive-Schwarz operator. Let \(\bar L_e\) denote the active-row restriction
of the localization matrix \(L_e\) in \eqref{eq:gpu-boolean-localization-matrix}, obtained by removing any padding rows, and define
\begin{equation*}
	\Sigma_{e}^{\mathrm{loc}} :=\bar L_e \Sigma_{\mathrm E}\bar L_e^\top, \qquad M^{-1} :=\sum_{e=1}^{N_m} \bar L_e^\top {\Sigma_{e}^{\mathrm{loc}}}^{-1}\bar L_e.
\end{equation*}
Each \(\Sigma_{e}^{\mathrm{loc}}\) is a dense principal block of dimension at most \(36\), corresponding to the four macro-external
faces. Because \(\Sigma_{\mathrm E}\succ0\), every block \(\Sigma_{e}^{\mathrm{loc}}\) is symmetric positive definite.
At every IPM iteration, their entries are gathered from \(\Sigma_{\mathrm E}\) using precomputed macro-element-to-CSR
position maps analogous to \(\pi_e\). The blocks are factorized and inverted in double precision using batched
routines, and the same maps scatter the entries of \({\Sigma_{e}^{\mathrm{loc}}}^{-1}\) into the CSR value array of \(M^{-1}\);
contributions on a face shared by two macro elements are added. Thus, no sparsity search is required during the
IPM iterations. Since the graph of \(\Sigma_{\mathrm E}\) already contains the dense face couplings induced by each macro
element, \(M^{-1}\) uses the same CSR structure as \(\Sigma_{\mathrm E}\) and introduces no additional structural fill.

The PCG iterations terminate when
\begin{equation*}
	\frac{\lVert \widehat\beta^{\nu}-\Sigma_{\mathrm E}\Delta\chi_k^{\nu} \rVert_2}{\lVert\widehat\beta^{\nu}\rVert_2}
	\le \eta_{\nu},
\end{equation*}
where \(\eta_{\rm a}\) and \(\eta_{\rm c}\) denote the relative tolerances for the predictor and corrector solves,
respectively. Following inexact-IPM strategies that progressively tighten the accuracy of the Newton solve as
complementarity decreases~\cite{gondzio2013inexact}, we keep \(\eta_{\rm a}\) fixed and set \(\eta_{\rm c}\) adaptively as
\begin{equation*}
	\eta_{\rm c} = \max\!\left\{ \eta_{\min}, \min\!\left( \eta_{\max}, c_\eta\mu^{\alpha_\eta} \right) \right\},
\end{equation*}
where \(c_\eta>0\) and \(\alpha_\eta>0\) control the scaling and decay of the tolerance, respectively, while
\(\eta_{\min}\) and \(\eta_{\max}\) provide lower and upper bounds.

For comparison in the numerical experiments, the implementation also includes scalar- and block-Jacobi preconditioners. Scalar Jacobi
uses the inverse diagonal of \(\Sigma_{\mathrm E}\), whereas block Jacobi independently inverts nonoverlapping \(27\times27\) diagonal
principal blocks. The order \(27\) is the largest multiple of the nine face-multiplier components that does not exceed the maximum order
of \(32\) supported by \texttt{cublasDmatinvBatched}. All three preconditioners are constructed in double precision and applied in single
precision to reduce their application cost. Scalar Jacobi is applied by diagonal scaling, block Jacobi by batched dense matrix--vector
products, and the macro-element preconditioner by a sparse matrix--vector product. The system matrix, Krylov vectors, inner products, and
residual tests remain in double precision.

Unlike sparse direct factorization, PCG stores no fill-generating global factors; its memory requirement is governed by the CSR matrices,
the fixed-size macro-element batches, and a small number of work vectors. This can make PCG useful when the direct factors exceed the
available device memory. The tradeoff is the repeated sparse matrix--vector products and global reductions, together with an iteration
count that may grow as the IPM systems become ill-conditioned. PCG is therefore presented as a memory-saving alternative rather than as a
uniformly faster replacement for the cuDSS solve.


\section{Numerical Results}
\label{sec:numerical-results}

This section evaluates the proposed macro-element condensation and its CPU and GPU implementations on three 3D benchmark problems, using
the stress condensation of \eqref{eq:stress-only-condensed-system} as a baseline. Both CPU variants share the same OpenMP primal--dual
IPM implementation and multithreaded Intel MKL PARDISO solver \cite{intel_mkl_pardiso}; they differ only in the variables eliminated
before the global solve and are denoted Stress--CPU and Macro--CPU, respectively. The cuDSS-based GPU implementation is denoted Macro--GPU;
for the elastoplastic benchmark, where cuDSS and memory-oriented PCG variants are tested, the solver name is appended to the method label.
Stress--CPU and Macro--CPU are compared to evaluate the additional benefits of macro-element condensation, while Macro--CPU and Macro--GPU
are compared to assess the acceleration achieved by the GPU implementation. MOSEK and CuClarabel are included as general-purpose CPU and
GPU reference solvers, respectively, to assess solution agreement and provide external performance baselines.

The selected benchmarks cover the principal formulations considered in this work: an elastoplastic Cook's membrane problem governed by
the Rankine criterion, the limit analysis of a steel I-section beam under torsion using the von Mises criterion, and the limit analysis
of a reinforced-concrete beam under three-point bending within the extended framework of Section~\ref{sec:structure-exploiting-condensation}.
Each benchmark is evaluated on two meshes to assess performance at different discretization scales.

The numerical experiments were performed using the computational resources of the Toubkal high-performance computing facility at
Mohammed VI Polytechnic University (UM6P). All CPU experiments, including the MOSEK runs, used \(24\) threads on nodes equipped with
Intel Xeon Platinum 8276L processors, while the GPU experiments used a single NVIDIA A100 GPU with \(80\,\mathrm{GB}\) of device memory.
Reported runtimes cover only IPM solver execution and exclude problem generation.

All in-house implementations used the common initialization \(x^{(0)}=y^{(0)}=0\), \(s^{(0)}=10^{3}e\), and \(z^{(0)}=10^{-3}e\), where
\(e\) is the product-cone identity introduced in Subsection~\ref{subsec:ipm-newton-system}; MOSEK and CuClarabel used their default
initialization procedures. All solver configurations used a feasibility tolerance of \(10^{-7}\). The in-house implementations and
CuClarabel used absolute- and relative-gap tolerances of \(10^{-5}\), whereas MOSEK retained its default relative-gap tolerance. Since
MOSEK's homogeneous formulation employs different residual scalings and termination tests, identical nominal gap tolerances do not
necessarily provide comparable solution accuracy. Indeed, for MOSEK on the coarser reinforced-concrete problem, relative-gap tolerances
of \(10^{-5}\) and \(10^{-7}\) produced \(\lambda^\star=0.799\) and \(1.490\), respectively, whereas the default setting gave \(1.514\),
within \(0.21\%\) of the value \(1.517\) obtained by the in-house implementations and CuClarabel.

\subsection{Elastoplastic Cook's membrane}
\label{subsec:numerical-cook-membrane}

The first benchmark is the 3D Cook's membrane, a tapered body clamped on its left face and subjected to a uniform upward
traction on its right face. The prescribed resultant force of \(0.125\,\mathrm{MN}\) corresponds to approximately \(99\%\)
of the collapse load (estimated by lower-bound limit analysis on the finer of the two meshes considered). The material is
modeled as isotropic linear elastic--perfectly plastic, with Young's modulus \(E=10\,\mathrm{GPa}\) and Poisson's ratio
\(\nu=0.499\). Plastic admissibility is governed by the Rankine criterion with equal compressive and tensile strengths,
\(f_c=f_t=20\,\mathrm{MPa}\). The numerical study uses two meshes containing \(67{,}904\) and \(512{,}936\) tetrahedral
elements. Figure~\ref{fig:cook-membrane-benchmark} shows the problem geometry, loading, and boundary conditions together
with the computed normal-stress field \(\sigma_{xx}\) on the coarser mesh.

\begin{figure}[p]
    \centering
    \includegraphics[width=\textwidth]{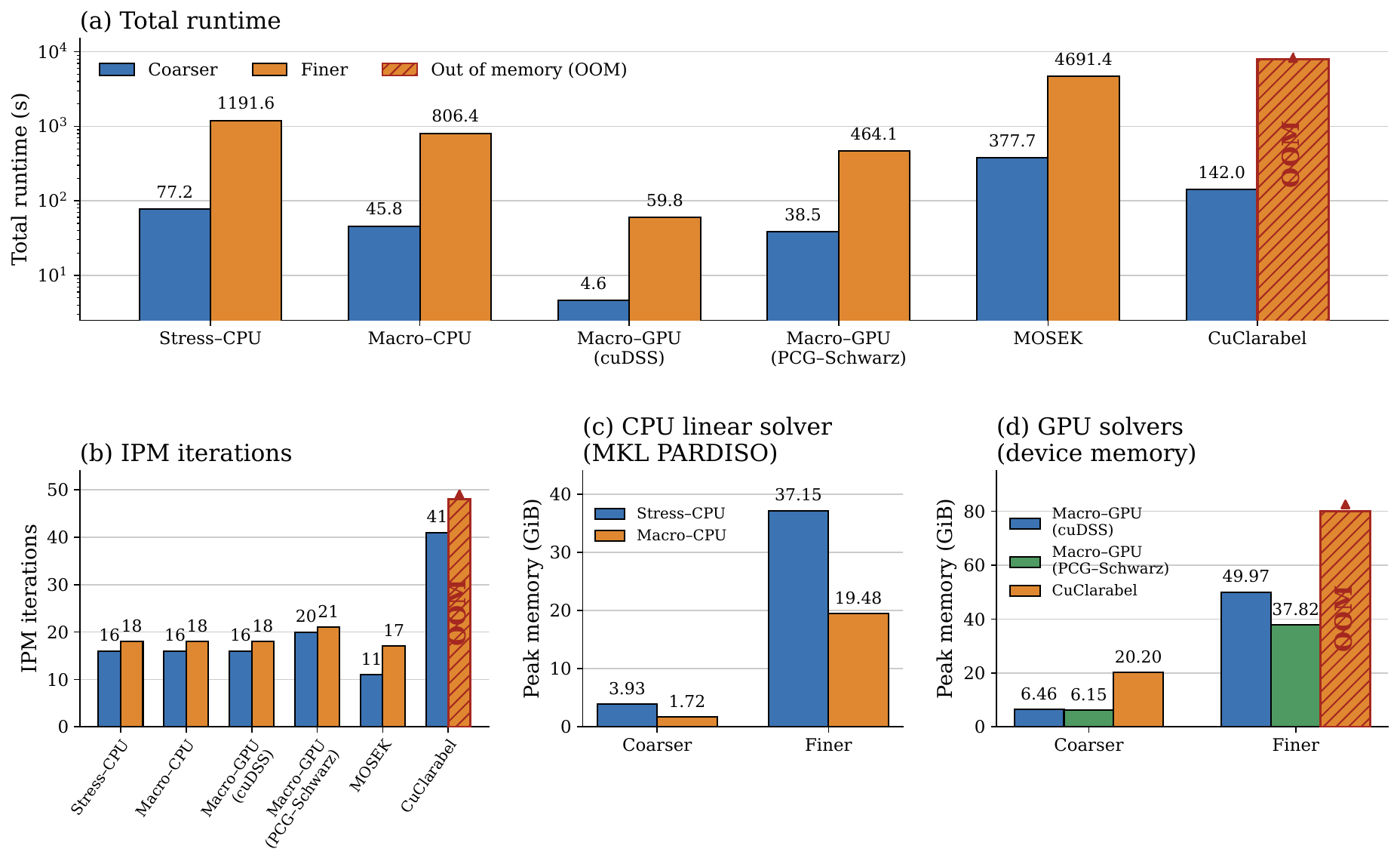}
    \caption{Computational performance for the Cook's membrane benchmark: (a) total runtime, (b) IPM iterations, (c) peak memory used by
	the MKL PARDISO linear solver, and (d) peak device memory of the GPU solvers. Hatched bars mark the finer-mesh CuClarabel out-of-memory
	result and do not represent measured values.}
    \label{fig:cook-performance}
\end{figure}

Table~\ref{tab:cook-condensation} shows that stress condensation reduces the global dimension by a factor of approximately \(2.16\),
but the structural-nonzero count by only \(1.21\), on both meshes. Eliminating the macro-internal multipliers provides a further
reduction of approximately \(4\times\) in both dimension and structural nonzeros. The resulting original-to-macro ratios are close
to the estimates in Equations~\eqref{eq:dimension-reduction-ratios} and~\eqref{eq:nnz-reduction-estimate}. These structural reductions
translate into measurable computational savings. Figure~\ref{fig:cook-performance} shows that, relative to Stress--CPU, Macro--CPU is
\(1.69\times\) and \(1.48\times\) faster and reduces the peak memory used by PARDISO by \(56.4\%\) and \(47.5\%\) on the coarser and
finer meshes, respectively.

\begin{figure}[p]
	\centering
	\begin{subfigure}[c]{0.47\textwidth}
		\centering
		\resizebox{\linewidth}{!}{\begin{tikzpicture}[
    x=1cm,
    y=1cm,
    line cap=round,
    line join=round,
    body edge/.style={draw=black!72,line width=0.85pt},
    dimension/.style={<->,draw=black!75,thin},
    extension/.style={draw=black!45,thin},
    load/.style={->,draw=black,line width=0.9pt}
]
    \coordinate (A) at (0.00,0.00);
    \coordinate (B) at (0.00,4.40);
    \coordinate (C) at (5.30,6.00);
    \coordinate (D) at (5.30,4.40);

    \coordinate (Ap) at (1.10,-0.28);
    \coordinate (Bp) at (1.10,4.12);
    \coordinate (Cp) at (6.40,5.72);
    \coordinate (Dp) at (6.40,4.12);

    \draw[body edge] (A)--(B)--(C)--(D)--cycle;
    \draw[body edge] (Ap)--(Bp)--(Cp)--(Dp)--cycle;
    \draw[body edge] (A)--(Ap) (B)--(Bp) (C)--(Cp) (D)--(Dp);

    \foreach \s in {0.20,0.72,...,4.36}
        \draw[thin,black!65] (-0.28,\s-0.22)--(0.02,\s+0.08);
    \node[anchor=east] at (-1.02,2.20) {clamped};

    \draw[load] (5.50,4.47)--(5.50,5.26);
    \draw[load] (5.82,4.39)--(5.82,5.18);
    \draw[load] (6.14,4.31)--(6.14,5.10);
    \node[above] at (5.82,5.30) {$\bar{\boldsymbol t}$};

    \draw[extension] (B)--(-0.12,4.82);
    \draw[extension] (C)--(5.18,6.42);
    \draw[dimension] (-0.05,4.72)--
        node[above,sloped,fill=white,inner sep=1.5pt]
        {$0.48\,\mathrm{m}$}(5.23,6.31);

    \draw[extension] (A)--(-0.78,0.00);
    \draw[extension] (B)--(-0.78,4.40);
    \draw[dimension] (-0.64,0.00)--
        node[fill=white,inner sep=1.5pt,rotate=90]
        {$0.44\,\mathrm{m}$}(-0.64,4.40);

    \draw[extension] (Dp)--(7.10,4.12);
    \draw[extension] (Cp)--(7.10,5.72);
    \draw[dimension] (6.94,4.12)--
        node[fill=white,inner sep=1.5pt,rotate=90]
        {$0.16\,\mathrm{m}$}(6.94,5.72);

    \draw[dimension] (B)--
        node[below,sloped,fill=white,inner sep=1.2pt]
        {$0.10\,\mathrm{m}$}(Bp);
\end{tikzpicture}}
		\caption{Geometry, loading, and boundary conditions.}
		\label{fig:cook-membrane-geometry}
	\end{subfigure}
	\hfill
	\begin{subfigure}[c]{0.52\textwidth}
		\centering
		\includegraphics[width=0.92\linewidth, trim={35mm 40mm 25mm 72mm}, clip]{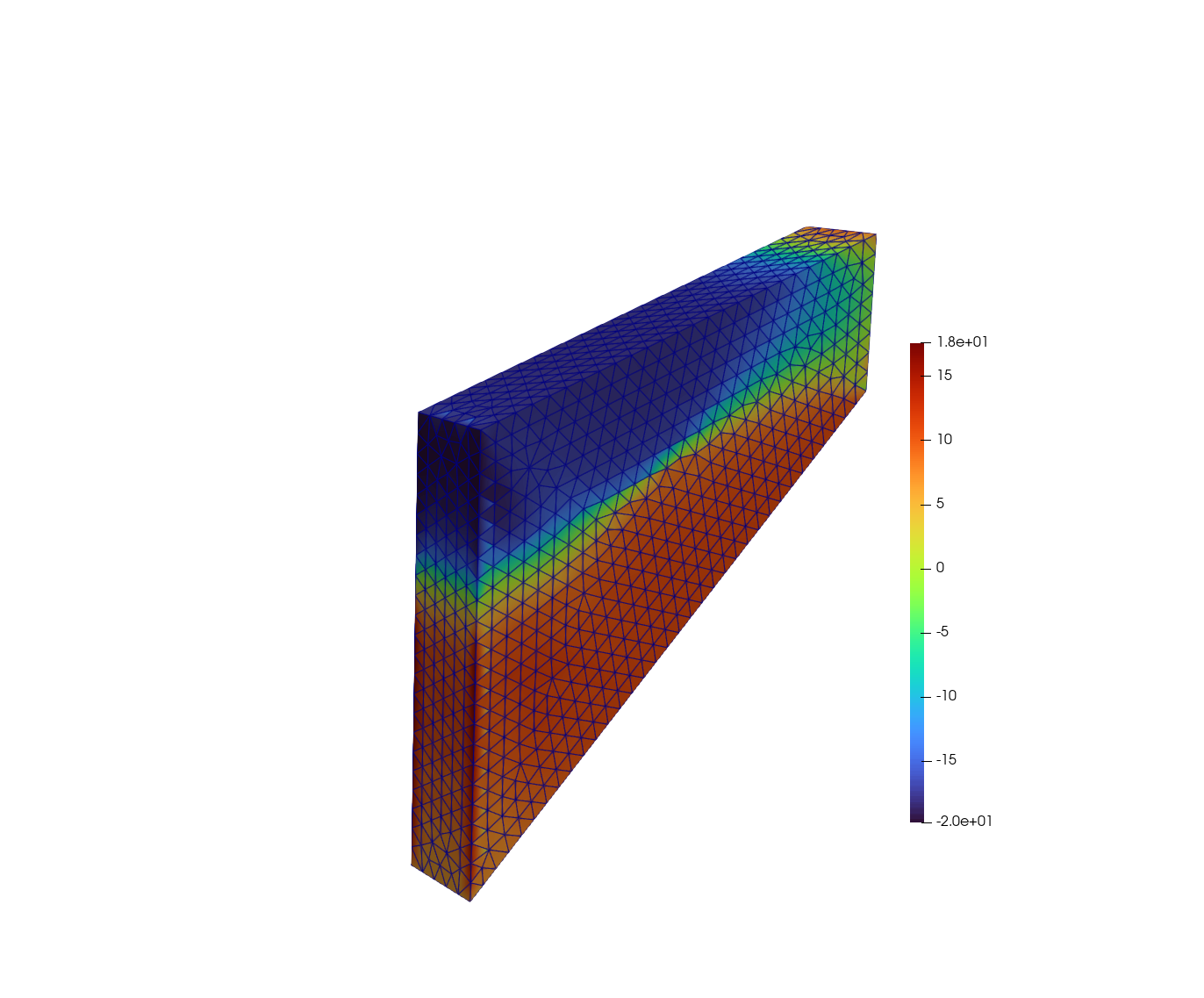}
		\caption{Normal stress \(\sigma_{xx}\), in \(\mathrm{MPa}\), on the coarser mesh.}
		\label{fig:cook-sigma-xx}
	\end{subfigure}
	\caption{Cook's membrane benchmark: geometry and representative stress solution. The geometric schematic is
		not drawn to scale.}
	\label{fig:cook-membrane-benchmark}
\end{figure}
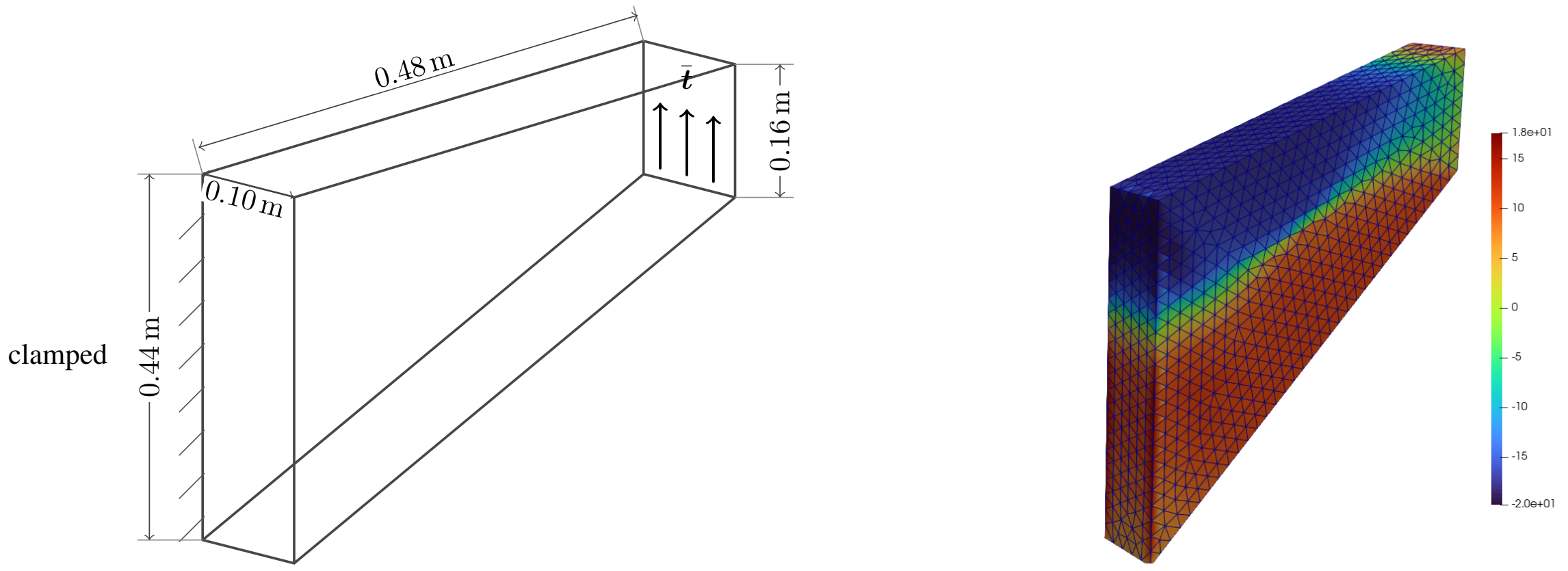
\begin{table}[p]
	\centering
	\caption{Newton-system dimensions and structural-nonzero counts for the Cook's membrane benchmark. The reported ratios are
		\(n_{\mathrm{full}}/n_{\mathrm{cond}}\) and \(\operatorname{nnz}_{\mathrm{full}}/\operatorname{nnz}_{\mathrm{cond}}\), respectively.}
	\label{tab:cook-condensation}
	\small
	\setlength{\tabcolsep}{3pt}
	\begin{tabular*}{\textwidth}{
		@{\extracolsep{\fill}}
		l r r r r c r r r c
		@{}
	}
		\toprule
		& & \multicolumn{4}{c}{Dimension}
		& \multicolumn{4}{c}{Structural nonzeros} \\
		\cmidrule(lr){3-6} \cmidrule(lr){7-10}
		Mesh & Elements
		& $n_{\mathrm{full}}$
		& $n_{\mathrm{stress}}$
		& $n_{\mathrm{cond}}$
		& Ratio
		& $\operatorname{nnz}_{\mathrm{full}}$
		& $\operatorname{nnz}_{\mathrm{stress}}$
		& $\operatorname{nnz}_{\mathrm{cond}}$
		& Ratio \\
		\midrule
		Coarser
		& $67{,}904$
		& $2.66{\times}10^{6}$
		& $1.23{\times}10^{6}$
		& $3.15{\times}10^{5}$
		& $8.43$
		& $4.72{\times}10^{7}$
		& $3.91{\times}10^{7}$
		& $9.75{\times}10^{6}$
		& $4.838$ \\
		Finer
		& $512{,}936$
		& $2.00{\times}10^{7}$
		& $9.27{\times}10^{6}$
		& $2.35{\times}10^{6}$
		& $8.53$
		& $3.57{\times}10^{8}$
		& $2.95{\times}10^{8}$
		& $7.37{\times}10^{7}$
		& $4.843$ \\
		\bottomrule
	\end{tabular*}
\end{table}

\FloatBarrier

The cuDSS-based Macro--GPU implementation is \(9.91\times\) and \(13.48\times\) faster than Macro--CPU on the coarser and finer meshes,
respectively. Against the general-purpose references, it is approximately \(80\times\) faster than MOSEK on both meshes and \(31\times\)
faster than CuClarabel on the coarser mesh. CuClarabel exceeds the available device memory on the finer mesh.

\begin{figure}[!t]
    \centering
    \includegraphics[width=0.90\textwidth]{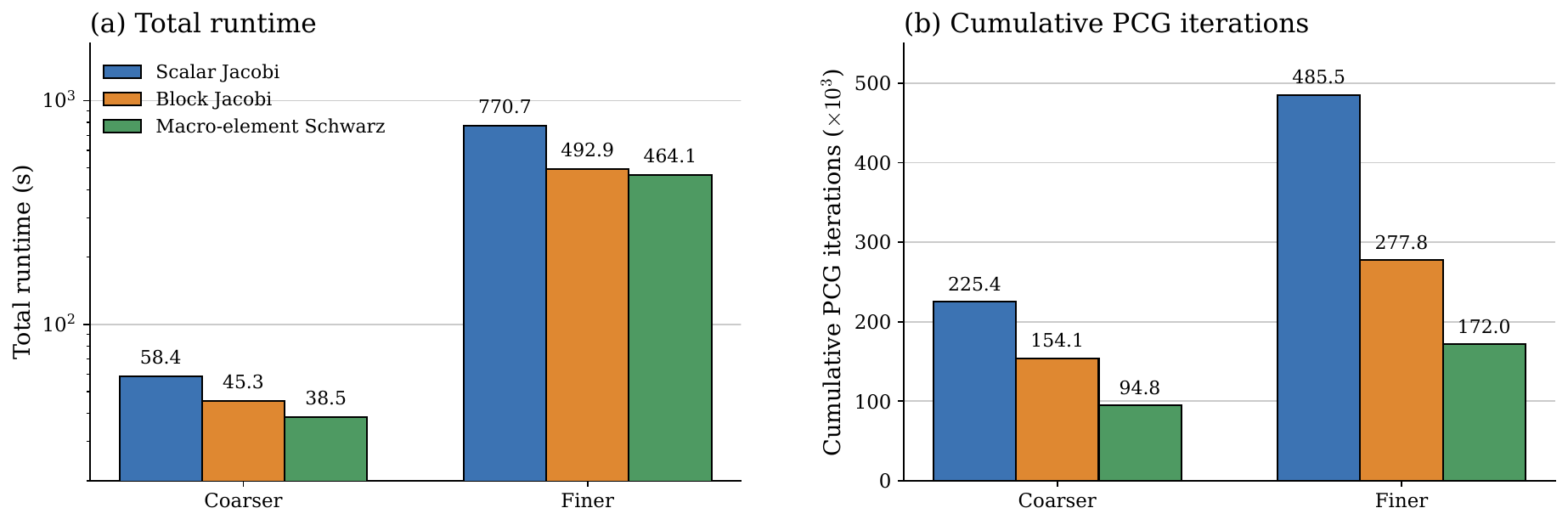}
    \caption{Influence of the PCG preconditioner for the Cook's membrane benchmark: (a) total Macro--GPU (PCG) runtime and
    (b) cumulative PCG iterations over all IPM linear solves. All runs use \(c_\eta=\alpha_\eta=0.5\) and \(\eta_{\min}=10^{-6}\).
    The fixed predictor tolerance and \(\eta_{\max}\) are \(10^{-3}\) for the coarser mesh and \(10^{-4}\) for the finer mesh.}
    \label{fig:cook-pcg-preconditioners}
\end{figure}

The PCG--Schwarz variant of Macro--GPU permits inexact Newton directions through its adaptive linear-solve tolerance and requires more IPM
iterations than the cuDSS-based implementation. It is approximately \(8\times\) slower than the direct variant, but reduces peak device
memory by \(4.8\%\) on the coarser mesh and \(24.3\%\) on the finer mesh. It nevertheless remains faster than Macro--CPU on both meshes,
providing a lower-memory GPU alternative at a higher runtime cost than the direct solve.

Figure~\ref{fig:cook-pcg-preconditioners} compares scalar Jacobi, block Jacobi, and the proposed macro-element additive-Schwarz operator
within the PCG implementation of Macro--GPU. Relative to scalar Jacobi, the macro-element preconditioner reduces the cumulative PCG
iteration count by approximately \(58\)--\(65\%\) and total runtime by \(34\)--\(40\%\) across the two meshes. Compared with block Jacobi,
it reduces the iteration count by approximately \(38\%\) on both meshes, while the runtime savings are more modest: \(15.0\%\) on the
coarser mesh and \(5.8\%\) on the finer mesh.

The smaller runtime gains relative to block Jacobi reflect the different costs of applying the two preconditioners on the GPU. Block Jacobi
comprises independent dense blocks and is applied using efficient batched dense matrix--vector kernels. By contrast, the overlapping
macro-element preconditioner is currently applied through a general sparse matrix--vector product, whose low arithmetic intensity and
irregular memory access limit GPU efficiency. Coloring the macro-element overlap graph would partition the blocks into independent sets
and enable a matrix-free application through batched dense kernels without conflicting memory updates. Even in its current sparse form,
however, the macro-element preconditioner yields the lowest runtime among the three PCG preconditioner variants of Macro--GPU.

\FloatBarrier

\subsection{Limit analysis of an I-section beam under torsion}
\label{subsec:numerical-steel-torsion}

\begin{figure}[p]
	\centering
	\begin{subfigure}[c]{0.49\textwidth}
		\centering
		\includegraphics[width=0.95\linewidth, trim={0mm 42mm 0mm 75mm}, clip]{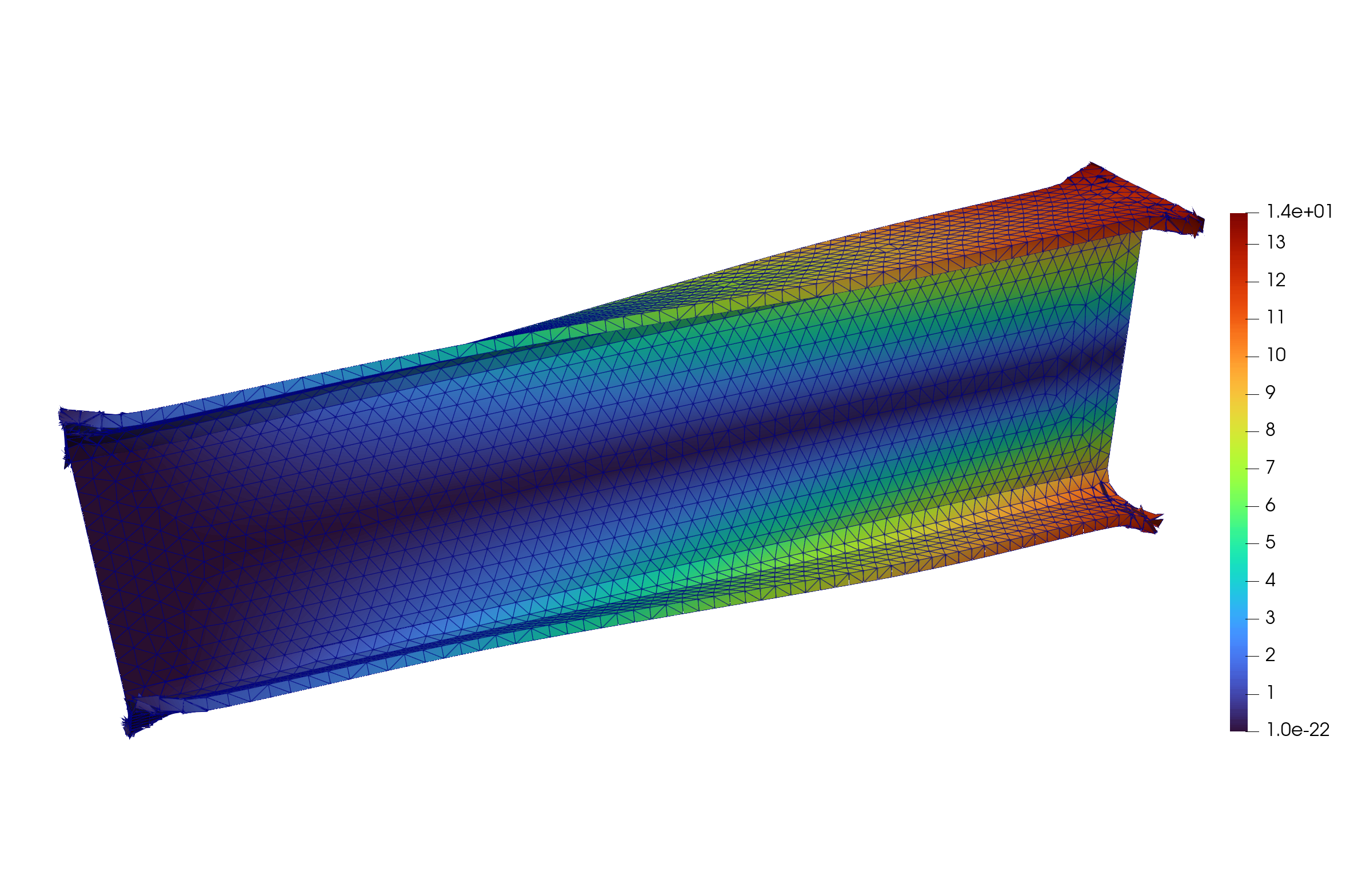}
		\caption{Collapse mechanism.}
		\label{fig:steel-displacements}
	\end{subfigure}
	\hfill
	\begin{subfigure}[c]{0.49\textwidth}
		\centering
		\includegraphics[width=0.95\linewidth, trim={0mm 42mm 0mm 75mm}, clip]{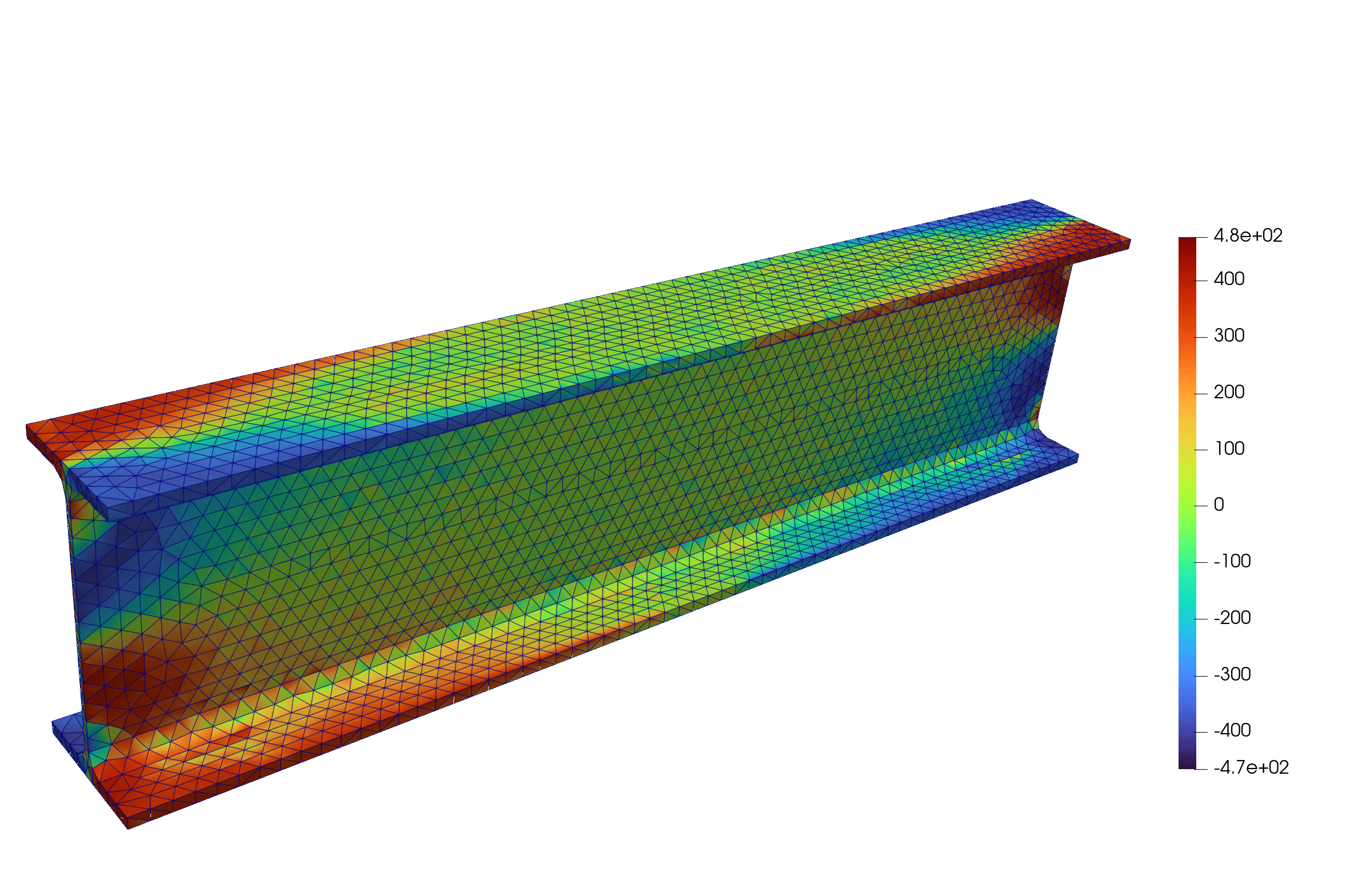}
		\caption{Normal stress \(\sigma_{xx}\), in \(\mathrm{MPa}\).}
		\label{fig:steel-sigma-xx}
	\end{subfigure}
	\caption{IPE~240 beam at torsional collapse on the coarser mesh.}
	\label{fig:steel-solutions}
\end{figure}

We consider a \(1\,\mathrm m\)-long cantilever steel beam with an IPE 240 cross-section. One end of the beam is clamped, and a reference
torque \(M_T=10\,\mathrm{kN\,m}\) is applied at the free end. The load multiplier \(\lambda\) therefore scales the applied torque, so
that the collapse torque is \(\lambda M_T\). The steel is modeled as perfectly plastic with uniaxial yield stress \(\sigma_y=500\,\mathrm{MPa}\),
and plastic admissibility is imposed through the von Mises criterion. To enable both stress and macro-element condensations, the trace
augmentation described in Section~\ref{sec:local-condensation-block-nonsingularity} is applied with \(\bar\tau=10\sigma_y=5000\,\mathrm{MPa}\).
To improve numerical robustness in MOSEK and CuClarabel, the reference torque is divided by \(10^2\) in all solvers, including our
implementations. The computed load multipliers are divided by \(10^2\) for reporting relative to the original reference torque \(M_T\).
Two discretization levels comprising \(61{,}764\) and \(458{,}312\) tetrahedral elements are considered. The computed collapse mechanism
and normal-stress field for the coarser mesh are shown in Figure~\ref{fig:steel-solutions}.

\begin{figure}[p]
	\centering
	\includegraphics[width=0.95\textwidth]{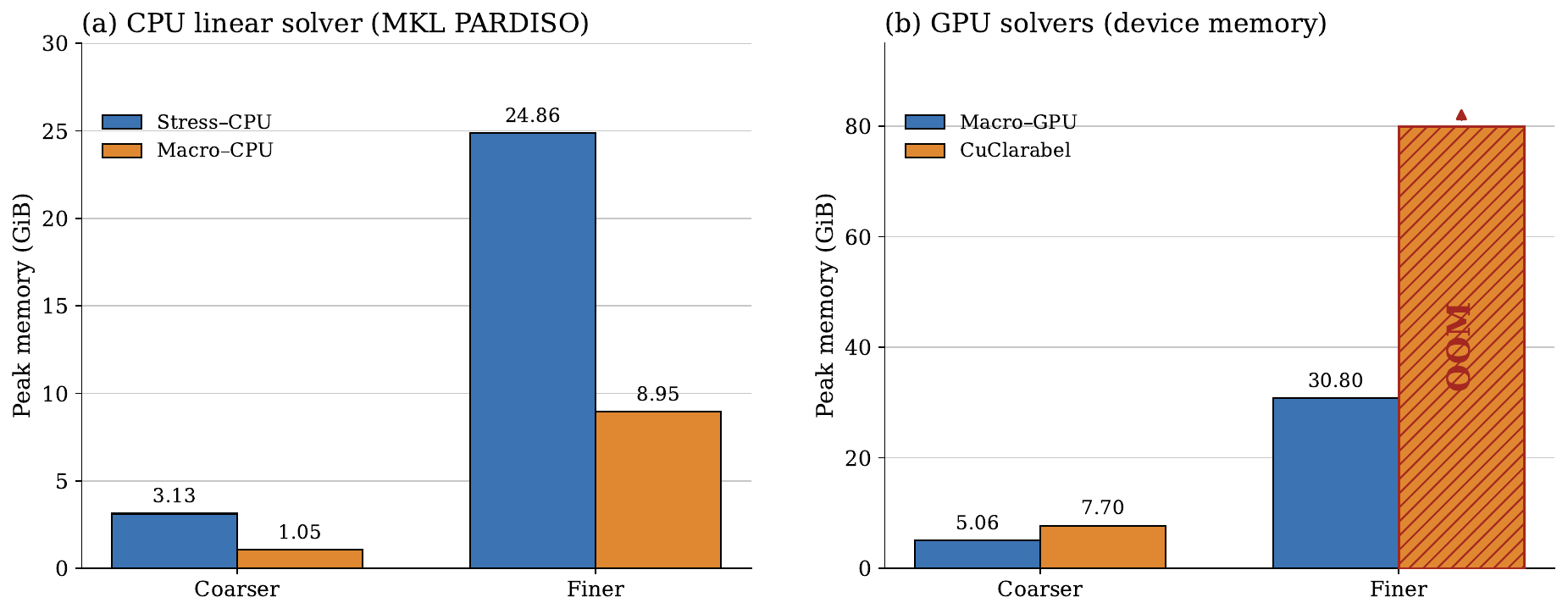}
	\caption{Peak memory usage for the IPE~240 torsion benchmark: (a) MKL PARDISO memory for Stress--CPU and Macro--CPU;
		(b) device memory for Macro--GPU and CuClarabel. OOM denotes exhaustion of the available GPU memory.}
	\label{fig:steel-memory}
\end{figure}

Table~\ref{tab:steel-condensation} shows that, relative to stress condensation, macro-element condensation further reduces the system
dimension by approximately \(3.6\)--\(3.7\times\) and the structural-nonzero count by approximately \(4\times\) on both meshes.
These reductions translate into CPU speedups of \(1.70\times\) and \(1.97\times\) for Macro--CPU over Stress--CPU on the coarser and finer
meshes, respectively. The additional elimination also reduces peak PARDISO memory by approximately \(66.5\%\) and \(64.0\%\), respectively,
as shown in Figure~\ref{fig:steel-memory}.

\begin{table}[p]
	\centering
	\caption{Newton-system dimensions and structural-nonzero counts for the IPE~240 torsion benchmark. The reported ratios are
	\(n_{\mathrm{full}}/n_{\mathrm{cond}}\) and \(\operatorname{nnz}_{\mathrm{full}}/ \operatorname{nnz}_{\mathrm{cond}}\), respectively.}
	\label{tab:steel-condensation}
	\small
	\setlength{\tabcolsep}{3pt}
	\begin{tabular*}{\textwidth}{
		@{\extracolsep{\fill}}
		l r r r r c r r r c
		@{}
	}
		\toprule
		& & \multicolumn{4}{c}{Dimension}
		& \multicolumn{4}{c}{Structural nonzeros} \\
		\cmidrule(lr){3-6} \cmidrule(lr){7-10}
		Mesh & Elements
		& $n_{\mathrm{full}}$
		& $n_{\mathrm{stress}}$
		& $n_{\mathrm{cond}}$
		& Ratio
		& $\operatorname{nnz}_{\mathrm{full}}$
		& $\operatorname{nnz}_{\mathrm{stress}}$
		& $\operatorname{nnz}_{\mathrm{cond}}$
		& Ratio \\
		\midrule
		Coarser
		& $61{,}764$
		& $2.45{\times}10^{6}$
		& $1.16{\times}10^{6}$
		& $3.23{\times}10^{5}$
		& $7.59$
		& $4.10{\times}10^{7}$
		& $3.58{\times}10^{7}$
		& $9.11{\times}10^{6}$
		& $4.50$ \\
		Finer
		& $458{,}312$
		& $1.81{\times}10^{7}$
		& $8.46{\times}10^{6}$
		& $2.28{\times}10^{6}$
		& $7.95$
		& $3.08{\times}10^{8}$
		& $2.65{\times}10^{8}$
		& $6.70{\times}10^{7}$
		& $4.60$ \\
		\bottomrule
	\end{tabular*}
\end{table}
\begin{table}[p]
	\centering
	\caption{Computational performance for the IPE~240 torsion benchmark. Here, \(\lambda^\star\) is the collapse-load multiplier, Iter
	    denotes the number of IPM iterations, and the dagger marks the final iterate returned by MOSEK after numerical stalling.}
	\label{tab:steel-performance}
	\small
	\setlength{\tabcolsep}{8pt}
	\begin{tabular}{l c c c c c c}
		\toprule
		& \multicolumn{3}{c}{Coarser}
		& \multicolumn{3}{c}{Finer} \\
		\cmidrule(lr){2-4} \cmidrule(lr){5-7}
		Method
		& $\lambda^\star$ & Iter & Time (s)
		& $\lambda^\star$ & Iter & Time (s) \\
		\midrule
		Stress--CPU
		& $0.821$ & 19 & 42.1
		& $0.891$ & 23 & 443.3 \\
		Macro--CPU
		& $0.821$ & 19 & 24.8
		& $0.891$ & 23 & 225.2 \\
		Macro--GPU
		& $0.821$ & 20 & 4.1
		& $0.891$ & 23 & 41.7 \\
		\midrule
		MOSEK
		& $0.821$ & 46 & 358.6
		& $0.881^{\dagger}$ & 26 & 2532.0 \\
		CuClarabel
		& $0.820$ & 52 & 107.0
		& -- & OOM & -- \\
		\bottomrule
	\end{tabular}
\end{table}
\FloatBarrier

The three in-house implementations also give matching collapse multipliers to the precision reported in Table~\ref{tab:steel-performance}.
On the coarser mesh, the augmented result agrees with MOSEK and closely matches CuClarabel. Uncondensed Newton solves without trace
augmentation recover the same multipliers on both meshes, confirming that the selected trace bound leaves the computed collapse loads
unchanged within the achieved numerical accuracy.

The cuDSS-based Macro--GPU implementation is \(6.05\times\) and \(5.40\times\) faster than Macro--CPU on the coarser and finer meshes,
respectively. On the coarser mesh, Macro--GPU requires one additional IPM iteration, which reflects the finite-precision effects in the
condensation and linear solves as the local matrices become ill-conditioned near convergence.

Against the general-purpose references, Macro--GPU is \(87.5\times\) faster than MOSEK and \(26.1\times\) faster than CuClarabel on the
coarser mesh. Figure~\ref{fig:steel-memory} also shows that its peak device-memory usage is \(34.3\%\) lower than CuClarabel's on that mesh.
On the finer mesh, Macro--GPU converges within the available device memory, while CuClarabel exhausts that memory and MOSEK returns an
unknown solution status after numerical stalling.

\FloatBarrier

\subsection{Limit analysis of a reinforced-concrete beam under three-point bending}
\label{subsec:numerical-reinforced-concrete}

For reinforced concrete, we adopt the embedded equilibrium formulation of~\cite{ferradi2024steel,ferradi2023elastoplastic}, using the
concrete discretization and notation of Section~\ref{sec:equilibrium-conic-formulation}. In the present discretization, the reinforcement
lines \(L^r\) are partitioned at macro-external faces into \(n_r\) straight segments \(\Delta L_j^r\), each contained in a single macro
element. Let \(\mathcal{SA}\) collect the locally equilibrated stress fields satisfying \(\operatorname{div}\sigma=0\) in each concrete
sub-element \(V_i\) and \(\mathrm d\sigma^r/\mathrm ds=0\) on each rebar segment, where \(\sigma^r\) is the scalar axial rebar stress and
\(s\) is arclength. Each segment therefore carries one constant stress unknown \(q_j^r\), shared across the concrete sub-elements
traversed within its macro element.

Perfect bonding at rebar-segment endpoints identifies the axial rebar displacement with the projection of the interpolated concrete-face
displacement along the bar tangent. For each internal mesh face \(F\) with linear displacement interpolation matrix \(N_F\), let \(I_F^r\)
collect the reinforcement coupling points: the rebar intersection points when \(F\) is macro-external, with \(I_F^r=\varnothing\)
when \(F\) is macro-internal. With interface equilibrium imposed in discrete weak form, the limit-analysis problem then reads
\begin{equation*}
	\begin{aligned}
		\min_{\lambda\in\mathbb{R},\,(\sigma,\sigma^r)\in\mathcal{SA}}\quad & -\lambda,\\
		\text{s.t.}\quad & \sigma \boldsymbol n=\lambda\bar t &&\text{on }\Gamma_t,\\
		&\int_F N_F^\top\left[\!\left[\sigma\right]\!\right]\boldsymbol n\,\mathrm d\Gamma+\sum_{p\in I_F^r}\mathcal A_p^r N_F(p)^\top e_p^r
		\left[\!\left[\sigma^r\right]\!\right]_p=0 && \text{for every internal face }F,\\
		&f(\sigma)\leq0 &&\text{in every }V_i,\\
		&|\sigma^r|\leq f_y &&\text{on }L^r.
	\end{aligned}
\end{equation*}
Here, the last two inequalities impose the Rankine criterion and the rebar yield bound, which admit PSD and SOC representations,
respectively; \(f_y\) denotes the rebar yield stress. At an intersection \(p\), \(\mathcal A_p^r\) and \(e_p^r\) denote the bar
cross-sectional area and oriented unit tangent, respectively; \(\left[\!\left[\sigma^r\right]\!\right]_p\) is the incoming minus
outgoing axial stress along \(e_p^r\). Details of the embedded reinforcement model and its discretization are
given in~\cite{ferradi2023elastoplastic,ferradi2024steel}.

In terms of the concrete stress coefficients \(q\) of \eqref{eq:stress-interpolation} and the rebar stress vector
\(q^r:=(q_1^r,\ldots,q_{n_r}^r)^\top\), the problem above admits an extended conic formulation with equality matrix
\eqref{eq:equality-local-global-structure}, where
\[
x_{\mathrm g}=\begin{bmatrix}q^r\\\lambda\end{bmatrix},\qquad A_{\mathrm g}=\begin{bmatrix}A_r&-\ell\end{bmatrix},\qquad n_{\mathrm g}=n_r+1,
\]
where \(A_r\) is the sparse assembled reinforcement matrix and \(\ell\) is the reference-load vector. The nonzero entries of
\(A_r\) are confined to macro-external equilibrium rows, leaving the concrete macro-internal equations unchanged.

The concrete and rebar yield conditions are separable, so no mixed primal Newton blocks arise between \(q\) and \(q^r\). The SOC bound
on each \(q_j^r\) contributes a positive scalar \(d_j\) to the primal Newton block, giving \(D_r:=\operatorname{diag}(d_1,\ldots,d_{n_r})\succ0\)
and \(\Omega_{\mathrm g}=\operatorname{blkdiag}(D_r,0)\). For Rankine concrete, \(\Omega_q\succ0\); together with \(D_r\succ0\), this
leaves only the \(\lambda\)-direction in \(\ker(\Omega)\). The nonzero load column \(-\ell\) therefore ensures
Assumption~\ref{ass:column-augmentation-regularity} and global Newton nonsingularity.

The macro-external coupling block in \eqref{eq:equality-local-global-structure} is \(B_{\mathrm E}=[\,B_{\mathrm E}^{r}\;-\ell_{\mathrm E}\,]\),
where \(B_{\mathrm E}^{r}\) contains the rebar contributions. With \(\Sigma_{\mathrm E}\) and \(\widehat\beta_{\mathrm E}\) defined in
\eqref{eq:macro-external-schur-complement} and \(b_r\) denoting the component of \(b_x\) associated with \(q^r\), eliminating the rebar
directions from \eqref{eq:compact-condensed-system} gives the final reduced system
\begin{equation*}
	\begin{gathered}
		\begin{bmatrix}
			\widetilde\Sigma_{\mathrm E}&\ell_{\mathrm E}\\
			\ell_{\mathrm E}^\top&0
		\end{bmatrix}
		\begin{bmatrix}\Delta y_{\mathrm E}\\\Delta\lambda\end{bmatrix}
		=\widetilde\beta,\\[1ex]
		\widetilde\Sigma_{\mathrm E}
		=\Sigma_{\mathrm E}+B_{\mathrm E}^{r}D_r^{-1}(B_{\mathrm E}^{r})^{\top},\qquad
		\widetilde\beta=
		\begin{bmatrix}
			\widehat\beta_{\mathrm E}+B_{\mathrm E}^{r}D_r^{-1}b_r\\
			-b_\lambda
		\end{bmatrix}.
	\end{gathered}
\end{equation*}
Since \(D_r\succ0\), the Schur correction is positive semidefinite, and Proposition~\ref{prop:external-schur-positive-definite} gives
\(\widetilde\Sigma_{\mathrm E}\succ0\).

\begin{figure}[!t]
	\centering
	\resizebox{\textwidth}{!}{\begin{tikzpicture}[
    x=1cm,
    y=1cm,
    line cap=round,
    line join=round,
    concrete/.style={draw=black!72,fill=black!3,line width=0.85pt},
    stirrup/.style={draw=black!48,line width=0.65pt},
    longitudinal/.style={draw=black!82,line width=1.25pt},
    dimension/.style={<->,draw=black!75,thin},
    extension/.style={draw=black!45,thin},
    load/.style={->,draw=black,line width=0.9pt},
    annotation/.style={font=\small,align=center}
]
    \coordinate (A) at (0.00,0.00);
    \coordinate (B) at (11.00,0.00);
    \coordinate (C) at (11.00,2.40);
    \coordinate (D) at (0.00,2.40);
    \draw[concrete] (A)--(B)--(C)--(D)--cycle;

    \foreach \x in {0.110,1.188,2.266,3.344,4.422,5.500,6.578,7.656,8.734,9.812,10.890}
        \draw[stirrup] (\x,0.12)--(\x,2.28);

    \draw[longitudinal]
        (0.193,0.65)
        .. controls (0.11,0.50) and (0.06,0.30) .. (0.18,0.12)
        -- (10.82,0.12)
        .. controls (10.94,0.30) and (10.89,0.50) .. (10.81,0.65);
    \foreach \x in {0.110,1.188,2.266,3.344,4.422,5.500,6.578,7.656,8.734,9.812,10.890}
        \fill[black!72] (\x,0.12) circle[radius=0.025];

    \draw[concrete,fill=black!10]
        (5.225,2.40) rectangle (5.775,2.72);
    \draw[black!65,line width=1.5pt] (5.225,2.72)--(5.775,2.72);
    \foreach \x in {5.27,5.38,5.50,5.62,5.73}
        \draw[load] (\x,3.28)--(\x,2.74);
    \node[annotation,anchor=east] at (4.92,3.10)
        {$\bar p=1\,\mathrm{MN/m^2}$};
    \draw[extension] (5.225,2.80)--(5.225,3.83);
    \draw[extension] (5.775,2.80)--(5.775,3.83);
    \draw[dimension] (5.225,3.70)--
        node[above=2pt,inner sep=1.2pt] {$0.5\,\mathrm m$}
        (5.775,3.70);
    \draw[extension] (5.775,2.40)--(6.10,2.40);
    \draw[extension] (5.775,2.72)--(6.10,2.72);
    \draw[dimension] (6.02,2.40)--
        node[right,fill=white,inner sep=1.0pt] {$0.1\,\mathrm m$}
        (6.02,2.72);

    \draw[concrete,fill=black!10] (0.00,-0.32) rectangle (0.55,0.00);
    \draw[concrete,fill=black!10] (10.45,-0.32) rectangle (11.00,0.00);
    \foreach \x in {0.00,0.10,0.20,0.30,0.40}
        \draw[thin,black!65] (\x,-0.34)--(\x+0.13,-0.57);
    \foreach \x in {10.45,10.55,10.65,10.75,10.85}
        \draw[thin,black!65] (\x,-0.34)--(\x+0.13,-0.57);
    \draw[extension] (0.00,-0.32)--(0.00,-0.98);
    \draw[extension] (0.55,-0.32)--(0.55,-0.98);
    \draw[dimension] (0.00,-0.88)--
        node[fill=white,inner sep=1.2pt] {$0.5\,\mathrm m$}
        (0.55,-0.88);
    \draw[extension] (10.45,-0.32)--(10.45,-0.98);
    \draw[extension] (11.00,-0.32)--(11.00,-0.98);
    \draw[dimension] (10.45,-0.88)--
        node[fill=white,inner sep=1.2pt] {$0.5\,\mathrm m$}
        (11.00,-0.88);

    \draw[extension] (0.00,-0.98)--(0.00,-1.42);
    \draw[extension] (11.00,-0.98)--(11.00,-1.42);
    \draw[dimension] (0.00,-1.32)--
        node[fill=white,inner sep=1.5pt] {$10\,\mathrm m$}
        (11.00,-1.32);

    \node[annotation,above right,fill=white,inner sep=1.2pt]
        (stirrup-label) at (7.35,2.84)
        {11 stirrups};
    \draw[->,thin,black!70] (stirrup-label.south west)--(7.656,2.05);
    \node[annotation,fill=white,inner sep=1.2pt]
        (bar-label) at (7.35,-0.35)
        {3 longitudinal bars};
    \draw[->,thin,black!70] (bar-label.north)--(7.08,0.12);

    \node[annotation] at (13.15,3.12) {cross-section};
    \draw[concrete] (12.55,0.00) rectangle (13.75,2.40);
    \draw[stirrup,rounded corners=2.5pt,line width=0.85pt]
        (13.75,2.28)--(12.67,2.28)--(12.67,0.12)--
        (13.63,0.12)--(13.63,2.40);
    \foreach \x in {12.79,13.15,13.51}
        \fill[black!82] (\x,0.12) circle[radius=0.06];

    \draw[extension] (12.55,0.00)--(12.18,0.00);
    \draw[extension] (12.55,2.40)--(12.18,2.40);
    \draw[dimension] (12.28,0.00)--
        node[fill=white,inner sep=1.2pt,rotate=90] {$1\,\mathrm m$}
        (12.28,2.40);
    \draw[extension] (12.55,0.00)--(12.55,-0.40);
    \draw[extension] (13.75,0.00)--(13.75,-0.40);
    \draw[dimension] (12.55,-0.31)--
        node[fill=white,inner sep=1.2pt] {$0.5\,\mathrm m$}
        (13.75,-0.31);
\end{tikzpicture}}
	\caption{Geometry of the reinforced-concrete beam.}
	\label{fig:rc-geometry}
\end{figure}
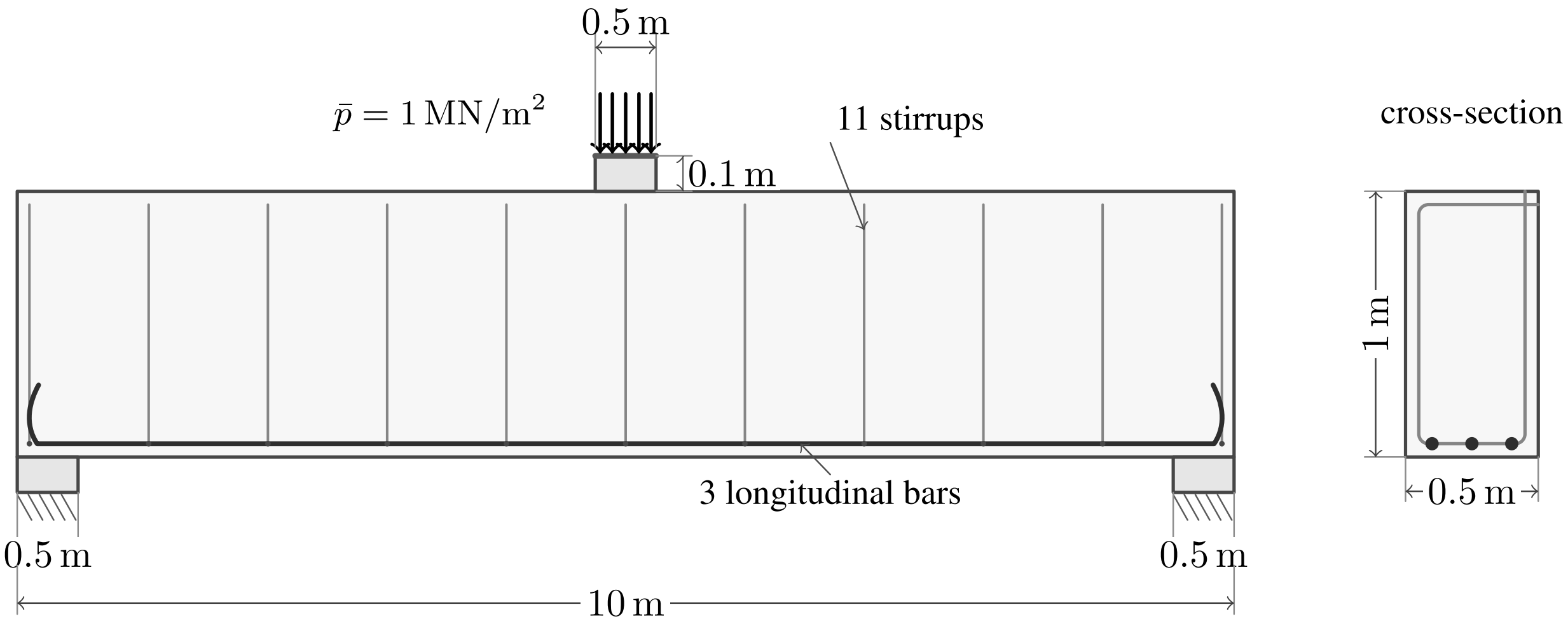

After the global solve, the rebar and concrete macro-local directions are recovered independently as
\begin{equation*}
	\begin{aligned}
		\Delta q^r
		&=D_r^{-1}\!\left(b_r-(B_{\mathrm E}^{r})^{\top}\Delta y_{\mathrm E}\right),\\
		\begin{bmatrix}\Delta q_e\\\Delta y_{\mathrm I,e}\end{bmatrix}
		&=K_e^{-1}\!\left(b_e-C_e\Delta y_{\mathrm E}\right).
	\end{aligned}
\end{equation*}
The concrete recovery follows \eqref{eq:local-solution}, reusing the local factors computed during condensation, while applications
of \(D_r^{-1}\) in the rebar corrections and recovery reduce to scalar reciprocals \(d_j^{-1}\). On the GPU, dedicated kernels
compute and assemble the rebar contributions to \(\widetilde\Sigma_{\mathrm E}\) and \(\widetilde\beta\) independently of the
macro-element operations, complementing the matrix assembly summarized in Algorithm~\ref{alg:gpu-condensed-assembly}.

The same scalar treatment of the rebar directions applies to the stress condensation \eqref{eq:stress-only-condensed-system} used
by Stress--CPU, since \(\Omega_q\succ0\). Eliminating \(\Delta q\) and \(\Delta q^r\) independently leaves a global system in
\((\Delta y_{\mathrm I},\Delta y_{\mathrm E},\Delta\lambda)\).

\begin{figure}[!t]
	\centering
	\begin{subfigure}[c]{0.46\textwidth}
		\centering
		\includegraphics[width=\linewidth, trim={80mm 2mm 220mm 100mm}, clip]{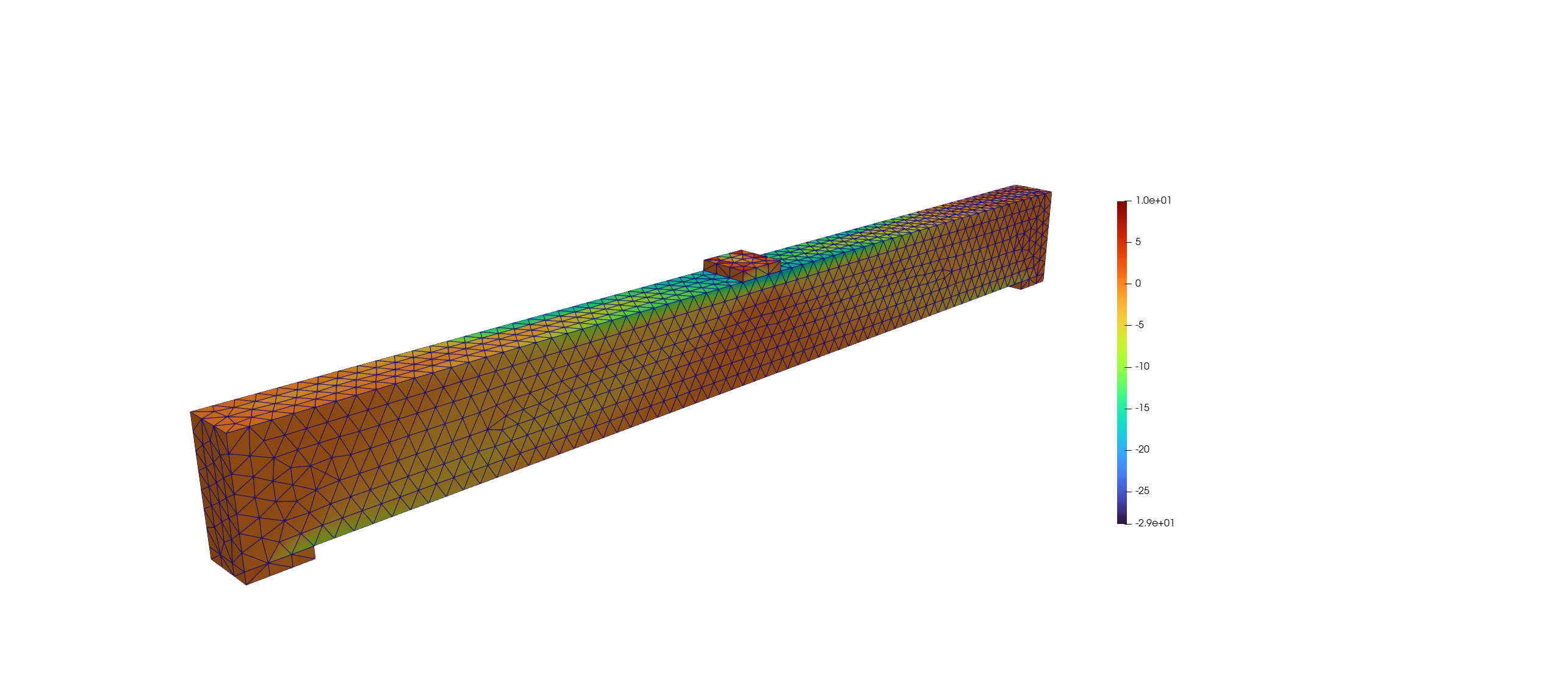}
		\caption{Concrete normal stress \(\sigma_{xx}\), in \(\mathrm{MPa}\).}
		\label{fig:rc-sigma-xx}
	\end{subfigure}
	\hfill
	\begin{subfigure}[c]{0.46\textwidth}
		\centering
		\includegraphics[width=\linewidth, trim={30mm 43mm 100mm 150mm}, clip]{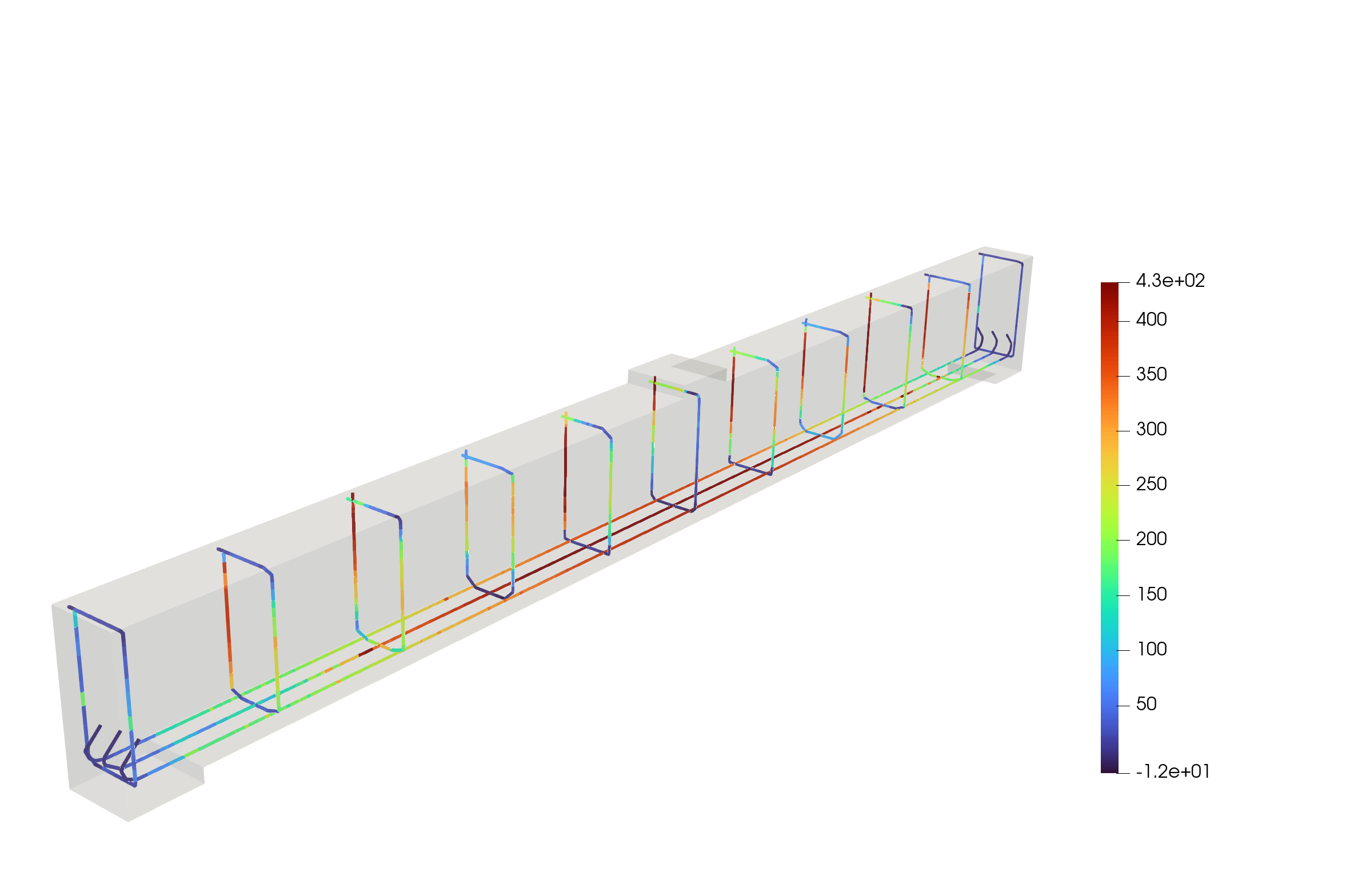}
		\caption{Rebar axial stress, in \(\mathrm{MPa}\).}
		\label{fig:rc-rebars-stress}
	\end{subfigure}
	\par\medskip
	\begin{subfigure}[c]{0.78\textwidth}
		\centering
		\includegraphics[width=\linewidth, trim={0mm 40mm 0mm 120mm}, clip]{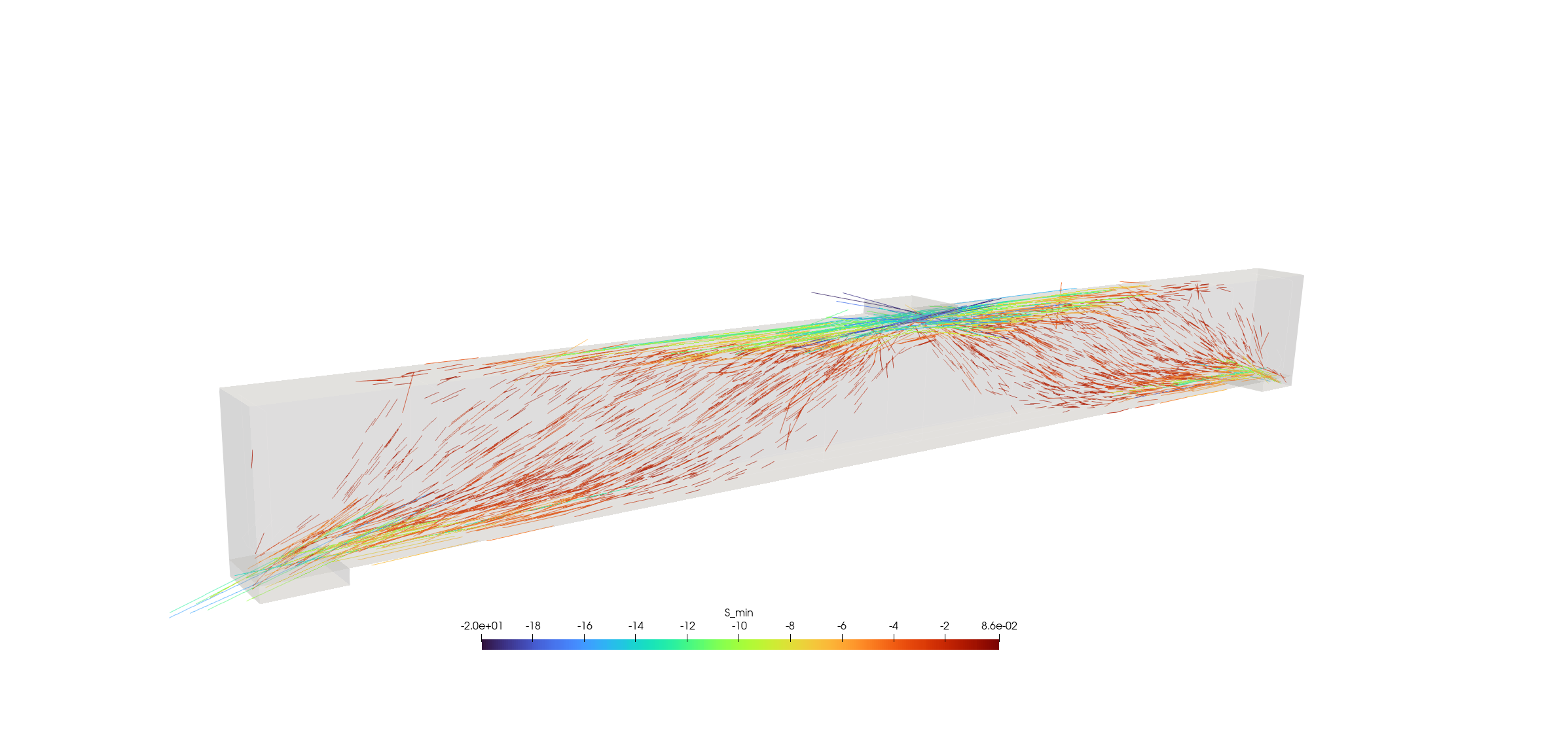}
		\caption{Compression flux.}
		\label{fig:rc-compression-flux}
	\end{subfigure}
	\caption{Reinforced-concrete beam at collapse on the coarser mesh.}
	\label{fig:rc-solutions}
\end{figure}

For the numerical assessment, we consider a \(10\,\mathrm{m}\)-long beam with a rectangular cross-section of width
\(0.5\,\mathrm{m}\) and height \(1\,\mathrm{m}\).
A loading block is placed at midspan on the upper face, and two support blocks are attached to the lower face at
the beam ends; each block is \(0.5\,\mathrm{m}\) long, spans the beam width, and is \(0.1\,\mathrm{m}\) high. A
pressure of \(1\,\mathrm{MN/m^2}\) is applied to the upper face of the loading block, while all displacement
components are prescribed on the lower faces of the support blocks. The concrete obeys the Rankine criterion with
\(f_c=20\,\mathrm{MPa}\) and \(f_t=0.1\,\mathrm{MPa}\), and the reinforcement consists of three longitudinal
HA20 bars and eleven HA10 stirrups spaced at \(0.98\,\mathrm{m}\) center-to-center, with \(0.1\,\mathrm{m}\) end offsets.
The rebar yield stress is \(f_y=435\,\mathrm{MPa}\), and no load is applied directly to the rebar ends.
The coarser configuration comprises \(53{,}816\) concrete tetrahedra and \(1{,}545\) rebar segments, whereas the finer
configuration contains \(441{,}208\) tetrahedra and \(3{,}302\) rebar segments.
Figure~\ref{fig:rc-geometry} summarizes the geometry, loading, boundary conditions, and reinforcement layout, while
Figure~\ref{fig:rc-solutions} presents the concrete normal-stress field \(\sigma_{xx}\), the rebar axial-stress
distribution, and the compression flux at collapse for the coarser configuration.

\begin{table}[!t]
	\centering
	\caption{Newton-system dimensions and structural-nonzero counts for the reinforced-concrete beam benchmark. The reported ratios are
		\(n_{\mathrm{full}}/n_{\mathrm{cond}}\) and \(\operatorname{nnz}_{\mathrm{full}}/ \operatorname{nnz}_{\mathrm{cond}}\), respectively.}
	\label{tab:rc-condensation}
	\small
	\setlength{\tabcolsep}{3pt}
	\begin{tabular*}{\textwidth}{
		@{\extracolsep{\fill}}
		l r r r r c r r r c
		@{}
	}
		\toprule
		& & \multicolumn{4}{c}{Dimension}
		& \multicolumn{4}{c}{Structural nonzeros} \\
		\cmidrule(lr){3-6} \cmidrule(lr){7-10}
		Mesh & Elements
		& $n_{\mathrm{full}}$
		& $n_{\mathrm{stress}}$
		& $n_{\mathrm{cond}}$
		& Ratio
		& $\operatorname{nnz}_{\mathrm{full}}$
		& $\operatorname{nnz}_{\mathrm{stress}}$
		& $\operatorname{nnz}_{\mathrm{cond}}$
		& Ratio \\
		\midrule
		Coarser
		& $53{,}816$
		& $2.12{\times}10^{6}$
		& $9.84{\times}10^{5}$
		& $2.57{\times}10^{5}$
		& $8.22$
		& $3.72{\times}10^{7}$
		& $3.11{\times}10^{7}$
		& $7.82{\times}10^{6}$
		& $4.76$ \\
		Finer
		& $441{,}208$
		& $1.73{\times}10^{7}$
		& $8.01{\times}10^{6}$
		& $2.05{\times}10^{6}$
		& $8.42$
		& $3.06{\times}10^{8}$
		& $2.54{\times}10^{8}$
		& $6.39{\times}10^{7}$
		& $4.80$ \\
		\bottomrule
	\end{tabular*}
\end{table}
\begin{table}[!t]
	\centering
	\caption{Computational performance for the reinforced-concrete beam
		benchmark.}
	\label{tab:rc-performance}
	\small
	\setlength{\tabcolsep}{8pt}
	\begin{tabular}{l c c c c c c}
		\toprule
		& \multicolumn{3}{c}{Coarser}
		& \multicolumn{3}{c}{Finer} \\
		\cmidrule(lr){2-4} \cmidrule(lr){5-7}
		Method
		& $\lambda^\star$ & Iter & Time (s)
		& $\lambda^\star$ & Iter & Time (s) \\
		\midrule
		Stress--CPU
		& $1.517$ & 75 & 199.7
		& $2.350$ & 128 & 3931.8 \\
		Macro--CPU
		& $1.517$ & 75 & 140.8
		& $2.350$ & 129 & 2889.4 \\
		Macro--GPU
		& $1.517$ & 74 & 13.7
		& $2.350$ & 129 & 281.8 \\
		\midrule
		MOSEK
		& $1.514$ & 47 & 457.1
		& $2.318$ & 62 & 6329.8 \\
		CuClarabel
		& $1.517$ & 110 & 205.0
		& -- & OOM & -- \\
		\bottomrule
	\end{tabular}
\end{table}
\begin{figure}[!htbp]
	\centering
	\includegraphics[width=0.88\textwidth]{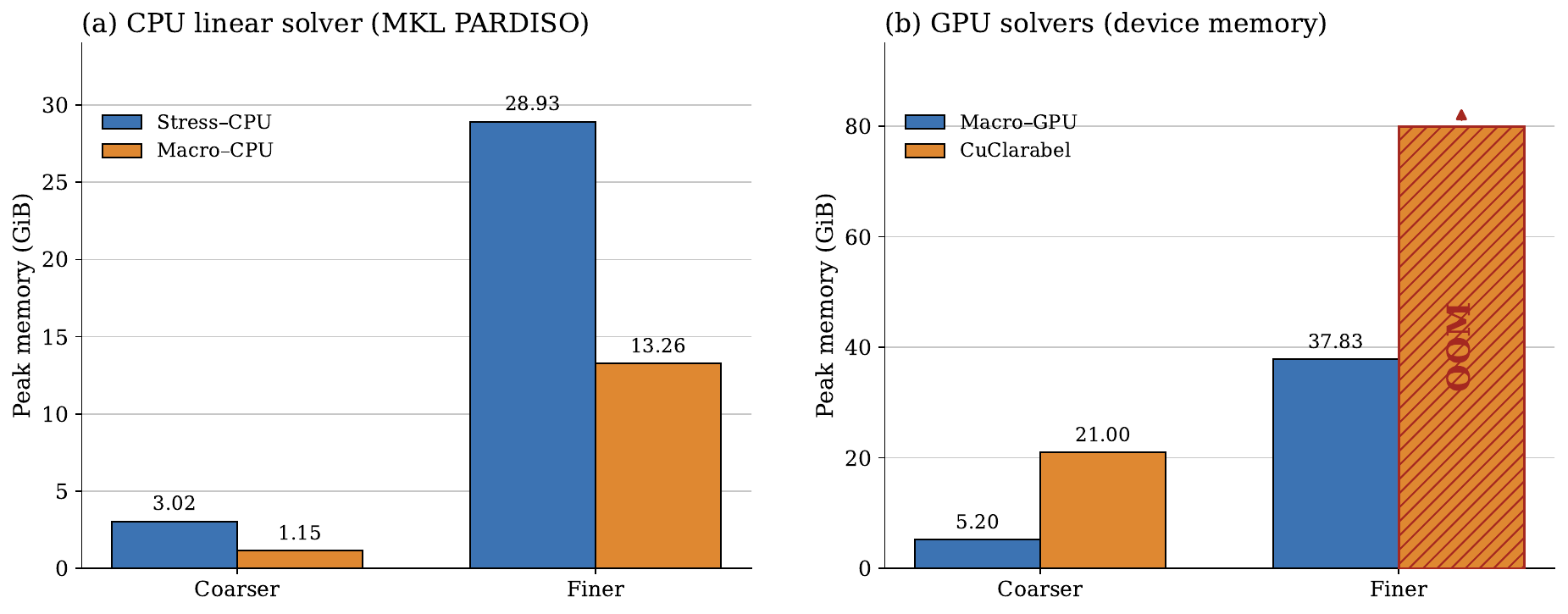}
	\caption{Peak memory usage for the reinforced-concrete beam benchmark: (a) MKL PARDISO memory for Stress--CPU and Macro--CPU;
		(b) device memory for Macro--GPU and CuClarabel. OOM denotes exhaustion of the available GPU memory.}
	\label{fig:rc-memory}
\end{figure}

Table~\ref{tab:rc-condensation} shows further reductions of approximately \(3.8\)--\(3.9\times\) in dimension and \(4\times\) in structural
nonzeros relative to stress condensation, consistent with the preceding benchmarks. The runtimes in Table~\ref{tab:rc-performance} show
that Macro--CPU is \(1.42\times\) and \(1.36\times\) faster than Stress--CPU on the coarser and finer configurations, respectively.
Figure~\ref{fig:rc-memory} shows reductions in peak PARDISO memory of \(61.9\%\) and \(54.2\%\), respectively.

The cuDSS-based Macro--GPU implementation is approximately \(10.3\times\) faster than Macro--CPU on both configurations. Against the
general-purpose references, its speedups over MOSEK are \(33.4\times\) and \(22.5\times\) on the coarser and finer configurations,
respectively, while it is \(15.0\times\) faster than CuClarabel on the coarser configuration. Its peak device-memory usage is also
\(75.2\%\) lower than CuClarabel's on the coarser configuration. On the finer configuration, Macro--GPU completes within the available
device memory, whereas CuClarabel exhausts that capacity.


\section{Conclusions}
\label{sec:conclusions}

This work extended equilibrium macro-element condensation to the iteration-dependent Newton systems arising in
primal--dual interior-point methods for equilibrium-based conic formulations of elastoplasticity and limit analysis.
The analysis established global and macro-local nonsingularity conditions and showed how both formulations can be
solved using positive-definite factorizations. Local nonsingularity holds independently of the yield criterion in
elastoplasticity, while the proposed trace augmentation enables condensation for von Mises limit analysis.
By eliminating macro-internal traction-continuity multipliers, the method reduces the global system relative to
classical stress-only elimination, with structural estimates predicting approximately \(4\times\) reductions in
both dimension and structural nonzeros. For the finer configurations tested, these structural reductions translate
into CPU speedups of up to \(1.97\times\) and reduce peak PARDISO memory by up to \(64.0\%\).

This work also developed a CUDA Fortran implementation that exploits the uniform macro-local block structure through
batched CUDA routines and uses cuDSS for the reduced sparse solve. The GPU implementation achieves speedups of up to
\(13.48\times\) over its CPU counterpart, and up to \(87.5\times\) and approximately \(31\times\) over MOSEK and
CuClarabel, respectively. All three problems are solved at their finest tested mesh size on a single NVIDIA A100 GPU
with \(80\,\mathrm{GB}\) of device memory, while CuClarabel exhausts that capacity in all three cases. For
elastoplasticity, a PCG alternative reduces device-memory requirements at the cost of longer runtimes than the direct GPU solver.

The nonsingularity analysis provides a basis for extending the framework to other yield criteria, shakedown analysis,
and topology optimization, subject to formulation-specific verification of the local and global nonsingularity
conditions. For larger three-dimensional discretizations, aggregating equilibrium super-elements into larger
macro-elements could eliminate additional interface variables, but the resulting growth in local block sizes and
condensed-matrix fill must be balanced against the reduction in the global system. Further studies could characterize
the conditioning of the condensed global matrix and macro-local blocks along the interior-point trajectory.
At the implementation level, the memory required by retained macro-local data could be substantially reduced by extracting
\(A_{\mathrm I,e}^{\top}\) on the fly and recomputing the factors of \(\Omega_e\) and \(S_{\mathrm I,e}\) in temporary batches for
condensation, right-hand-side formation, and recovery, at the cost of additional local computation. To further improve scalability,
a multi-GPU implementation could distribute the macro-local batches and reduced global solve across devices, increasing computational
resources and aggregate device-memory capacity. Hybrid host--device execution could further reduce device-memory requirements and
enable larger discretizations, with its benefits assessed against the associated transfer costs.

Iterative solution methods provide a complementary route to scalability by reducing the memory demands associated with sparse
direct factorization. Their further development is also motivated by the conditioning comparison with classical stress
condensation established in Proposition~\ref{prop:macro-external-conditioning}. For elastoplasticity, the present preconditioned
conjugate-gradient solver could be improved through GPU-efficient preconditioners that better balance iteration count and
application cost. In particular, graph coloring could partition the overlapping macro-element contributions into conflict-free
groups, enabling their parallel application through efficient batched dense kernels. For limit analysis, iterative solution is
more challenging because the stress block of the Newton system lacks the positive-definite compliance contribution present in
elastoplasticity and its conditioning is instead governed by the cone scaling, which may deteriorate near convergence. Further
work should therefore develop suitable iterative solvers and preconditioners that exploit the positive-definite macro-external
core together with the scalar Schur-complement treatment of the load multiplier while remaining robust to this deterioration.
These developments could provide memory-efficient solution paths when direct factorization becomes the principal
scalability bottleneck.


\section*{Declaration of competing interests}
The authors declare no competing interests.

\section*{Data availability}
All data is available from the corresponding author upon reasonable request.



\appendix
\section{Proofs of the equilibrium-operator regularity results}
\label{app:global-newton-regularity}

\subsection{Proof of Proposition~\ref{prop:assembled-equilibrium-full-rank}}
\label{app:proof-assembled-equilibrium-full-rank}

\begin{proof}
	For each macro element, let \(A_{\partial,e}\in\mathbb{R}^{36\times84}\) denote the connector on its four macro-external
	faces and define \(A_e^{\rm all}=[\,A_{\partial,e}^\top\;A_{\mathrm I,e}^\top\,]^\top\) as the complete connector on its four
	macro-external and six macro-internal faces. The classical stability result for the linear four-tetrahedron equilibrium
	macro element states that the left nullspace of the complete connector is precisely the six-dimensional rigid-body
	trace space~\cite[Section~5.5.1]{kempeneers2010pure}; see also~\cite[Section~5.2.1]{maunder2016stability}. The
	point-value connectors used in those analyses and the present coefficients of the linear face basis \(N_F\) are related
	by nonsingular facewise transformations~\cite{kempeneers2010pure}. Consequently, denoting the
	rigid-body traces in the coordinates of \(N_F\) by \(\mathcal R_e\),
	\begin{equation}
		\ker\!\left((A_e^{\rm all})^\top\right)=\mathcal R_e.
		\label{eq:complete-connection-rigid-nullspace}
	\end{equation}
	Now let \(y\in\ker(A_q^\top)\). Restricting \(y\) to each macro element and inserting zero coefficients on its
	external faces contained in \(\Gamma_u\) gives a complete local vector \(y_e^{\rm all}\) satisfying
	\((A_e^{\rm all})^\top y_e^{\rm all}=0\). Equation~\eqref{eq:complete-connection-rigid-nullspace} therefore shows that
	\(y_e^{\rm all}\) is the trace of a rigid-body field \(u_{R,e}\). Consistent assembly makes these traces identical on
	every face shared by two macro elements. Since rigid-body fields that agree on a face are identical, face connectivity
	of the macro-element mesh yields a single global rigid-body field \(u_R\). Its trace vanishes on \(\Gamma_u\), and the
	prescribed-displacement condition in Subsection~\ref{subsec:equilibrium-problem-discretization} thus gives \(u_R=0\).
	Hence \(y=0\), proving \(\ker(A_q^\top)=\{0\}\).
\end{proof}

\subsection{Proof of Proposition~\ref{prop:equilibrium-formulations-regularity}}
\label{app:proof-equilibrium-formulations-regularity}

\begin{proof}
	Proposition~\ref{prop:assembled-equilibrium-full-rank} gives the full row rank of \(A_q\), while
	Corollary~\ref{cor:augmented-equilibrium-full-rank} gives the same property for \(A_\lambda=[\,A_q\;-\ell\,]\), defined in
    \eqref{eq:limit-analysis-equality-matrix}. It therefore remains to verify the kernel-intersection condition in
	Assumption~\ref{ass:global-newton-regularity}.
	
	\medskip
	\noindent\emph{Elastoplasticity.}
	The positive definiteness of \(P_q\), established in \eqref{eq:assembled-compliance}, gives \(\ker(P_q)=\{0\}\). Hence
	the required kernel intersection is trivial, independently of the yield criterion.
	
	\medskip
	\noindent\emph{Rankine limit analysis.}
	Proposition~\ref{prop:rankine-operator-full-rank} shows that every sub-element yield block has full column rank. Since
	the stress part of the assembled yield operator is block diagonal by sub-element, any
	\((q,\lambda_0)\in\ker(G)\cap\ker(A_\lambda)\) satisfies \(q=0\), while the equality equation reduces to \(-\lambda_0\ell=0\).
	Assumption~\ref{ass:reference-loading-regularity} gives \(\ell\ne0\), and therefore \(\lambda_0=0\).
	
	\medskip
	\noindent\emph{Von Mises limit analysis.}
	If \((q,\lambda_0)\in\ker(G)\cap\ker(A_\lambda)\), Proposition~\ref{prop:von-mises-hydrostatic-kernel} implies that
	on every sub-element \(V_i\), \(\operatorname{mat}_v(S_iq_i)=p_i{\mathbf I}\) for some constant \(p_i\). The
	internal-face equations and mesh connectivity then give, by Lemma~\ref{lem:hydrostatic-pressure-propagation}, a single
	constant \(p\) such that \(q=pq_h\). The remaining equality equation becomes \(pA_q q_h-\lambda_0\ell=0\). The linear
	independence of \(\ell\) and \(A_q q_h\) in Assumption~\ref{ass:reference-loading-regularity} therefore gives
	\(p=\lambda_0=0\), and hence \(q=0\). This proves the kernel-intersection condition and completes the proof.
\end{proof}


\bibliographystyle{unsrt}  
\bibliography{references}

\end{document}